\documentclass[12pt,a4paper]{amsart}

\usepackage[utf8]{inputenc}
\usepackage[american, british]{babel}
\usepackage{amsmath, amsthm, amssymb, mathrsfs, graphicx}
\usepackage[dvipsnames]{xcolor}
\usepackage{hyperref}
\usepackage{enumitem}
\usepackage{dsfont}

\newtheorem{theorem}{Theorem}[section]
\newtheorem{proposition}[theorem]{Proposition}
\newtheorem{corollary}[theorem]{Corollary}
\newtheorem{lemma}[theorem]{Lemma}

\theoremstyle{remark}
\newtheorem{definition}[theorem]{Definition}
\newtheorem{claim}{Claim}
\newtheorem{remark}[theorem]{Remark}

\numberwithin{equation}{section}

\DeclareMathOperator{\hs}{Hess} % Hessian
\DeclareMathOperator{\rc}{Ric} % Ricci
\DeclareMathOperator{\tr}{tr} % trace
\DeclareMathOperator{\e}{\mathbf{e}} % basis vector
\DeclareMathOperator{\K}{{\mathbf{K}}} % null flag curv
\DeclareMathOperator{\R}{\mathds{R}} % real nos
\DeclareMathOperator{\N}{\mathds{N}} % natural nos
\DeclareMathOperator{\C}{\mathcal{C}} % cone structure
\DeclareMathOperator{\Cw}{\mathfrak C} % classical Weyl
\DeclareMathOperator{\F}{\mathcal{F}} % smooth real fxns
\DeclareMathOperator{\T}{\mathfrak{T}} %  anisotropic tensors
\DeclareMathOperator{\X}{\mathfrak{X}} % vector filds
\DeclareMathOperator{\vd}{\dot\partial} % vertical derivative
\DeclareMathOperator{\J}{\mathcal{R}} % Jacobi curv operator
\DeclareMathOperator{\W}{\mathcal{W}} % Weyl curv
\DeclareMathOperator{\tnab}{\widetilde{\nabla}} % cnx "til"
\DeclareMathOperator{\hnab}{\hat{\nabla}} % cnx "hat"

\renewcommand{\d}{\mathrm{d}} % differential

\newcommand{\spn}[1]{\mathrm{span}\{{#1}\}} % span
\newcommand{\hor}[1]{\overline{#1}} % "bar"
\newcommand{\til}[1]{\widetilde{#1}} % "til"

\newcommand{\bigslant}[2]{{\raisebox{.1em}{$#1$}/\raisebox{-.1em}{$#2$}}} % quotient

\makeatletter
\newcommand{\pushright}[1]{\ifmeasuring@#1\else\omit\hfill$\displaystyle{#1}$\fi\ignorespaces}
\newcommand{\pushleft}[1]{\ifmeasuring@#1\else\omit$\displaystyle{#1}$\hfill\fi\ignorespaces}
\makeatother

\def\bb#1\eb{\textcolor{blue}{#1}}
\def\br#1\er{\textcolor{red}{#1}}
\def\bv#1\ev{\textcolor{green}{#1}}

\begin{document}
%--------------------------------------------------------------------

%   TITLE
\title[Null flag curvature of conformal pseudo-Finsler metrics]{Null flag curvature of anisotropically conformal pseudo-Finsler metrics}

%   AUTHORs
\author[M.~Inagaki]{Marcelo Inagaki}
\address{Marcelo Inagaki \hfill\break\indent Instituto de Matem\'{a}tica e Estat\'{\i}stica, Universidade de S\~ao Paulo \hfill\break\indent Rua do Mat\~ao, 1010, 05508-090 S\~ao Paulo, S\~ao Paulo, Brazil}
\email{inagaki@ime.usp.br}

\author[M.~\'{A}.~Javaloyes]{Miguel \'Angel Javaloyes}
\address{Miguel \'Angel Javaloyes \textup{(corresponding author)} \hfill\break\indent Departamento de Matem\'{a}ticas, Universidad de Murcia \hfill\break\indent Campus de Espinardo, 30100 Espinardo, Murcia, Spain}
\email{majava@um.es}

\author[P.~Mar\c cal]{Patr\'icia Mar\c cal}
\address{Patr\'icia Mar\c cal \hfill\break\indent Instituto de Matem\'{a}tica e Estat\'{\i}stica, Universidade de S\~ao Paulo \hfill\break\indent Rua do Mat\~ao, 1010, 05508-090 S\~ao Paulo, S\~ao Paulo, Brazil}
\email{pmarcal@usp.br}

%   DATE
\date{\today}

%   ABSTRACT
\begin{abstract}
We consider the Zermelo navigation problem for a time-dependent Finsler	metric, which can be recast as the study of lightlike geodesics of a Lorentz--Finsler metric on a non-relativistic spacetime. To control how such geodesics separate, we introduce a null flag curvature, defined for any	pseudo-Finsler metric with a lightlike cone, and show that it governs the Taylor expansion of the norm of a Jacobi field. We then prove a formula	relating the null flag curvatures of two anisotropically conformal pseudo-Finsler metrics. Applying it to the spacetime associated with a	time-independent Randers metric, which is anisotropically conformal to a stationary spacetime, we obtain compact expressions for the flag curvature of Randers metrics, both in terms of their $(\alpha,\beta)$ data and of their Zermelo data $(h,W)$. These expressions allow us to compute explicitly the	flag curvature when $W$ is conformal, affine or projective, and when $\alpha$ is a K\"ahler metric. We also introduce a Weyl-type curvature, invariant under	anisotropic conformal changes, whose vanishing characterizes Finsler manifolds of scalar flag curvature, and we obtain the flag curvature of $F + \d f$, showing in particular that $F$ has scalar flag curvature if and only if	$F + \d f$ does.
\end{abstract}

%   SETUP
\maketitle
\bibliographystyle{abbrv}

%   SECTIONS
%--------------------------------------------------------------------
\section{Introduction}
%--------------------------------------------------------------------

General Relativity rests on the idea that space and time merge into a fabric in which measurements depend on the observer. In particular, spacelike distances are observer-dependent, since the notion of rest space changes from one observer to another. We consider instead the situation in which distances are independent of the observer but vary with time. The first issue is how to measure the length of a curve in this setting. A natural choice, suited to many familiar scenarios, is to regard the unit vectors of the metric as the prescribed (or maximal) velocities of a moving object, so that the arc-length parameter coincides, up to an additive constant, with the elapsed time. Note that, in contrast with General Relativity, time is here absolute.

This setting is precisely that of the Zermelo navigation problem, in which the admissible velocity is prescribed in every direction and one seeks the trajectory that minimizes the travel time. It is well known that this problem can be recast in terms of a spacetime (see \cite{cap11i, cap24w,  gib09s, jav20d}). We consider the non-relativistic spacetime $\R \times M$ endowed with a family of Finsler metrics $\{F_t\}$, one for each instant $t$, so that, for every $(t,p)$, the set of unit vectors is a strongly convex hypersurface enclosing the zero vector in its bounded region. For a trajectory $s\mapsto (t(s),x(s))$ whose first coordinate is the time, we have
%--------------------------------------------------------------------
\begin{equation*}
t(s) = t(0) + \int_0^s F_{t(\sigma)}\big( \dot x(\sigma) \big) \d\sigma , \text{ hence } \; \dot t(s) = F_{t(s)}\big( \dot x(s) \big) .
\end{equation*}
%--------------------------------------------------------------------
Consequently, a curve $\gamma(s)=(t(s),x(s))$ is admissible if and only if its velocity lies in the subset
%--------------------------------------------------------------------
\begin{equation*}
\C_{(t,p)} = \{(\tau,v) \in \R\times T_p M : \tau = F_t(v) , \, v \neq 0 \} .
\end{equation*}
%--------------------------------------------------------------------
These subsets may be regarded as the lightcones of a (Finsler) spacetime; in general, we call them a cone structure. More precisely, the Lorentz--Finsler metric
%--------------------------------------------------------------------
\begin{equation}
\label{LFintro}
L(\tau,v) = \tau^2 - F_t(v)^2
\end{equation}
%--------------------------------------------------------------------
has $\C_{(t,p)}$ as its lightlike cone, and the admissible curves are exactly the lightlike curves of $L$. We next identify the geodesics of these time-dependent metrics. Since length is identified with elapsed time, the Fermat principle for Finsler spacetimes \cite{per06f}, applied with the arrival observer $\alpha(s) = (s,q)$, shows that the time-minimizing admissible curves are lightlike pre-geodesics.

As in Riemannian and Finsler geometry, a basic question concerns how geodesics departing from a common point separate. In Finsler geometry, the flag curvature controls this behavior. For lightlike geodesics, however, the tangent plane to a variation through lightlike curves is degenerate with respect to the osculating metric of the curve $\gamma$, denoted $g_{\dot\gamma}$, so the flag curvature cannot be computed. Inspired by the null sectional curvature introduced by Harris in Lorentzian geometry \cite{har82t}, we define a null flag curvature and show that it measures the separation of lightlike geodesics. This is done following the strategy used in the Finslerian case \cite[\S~5.5]{bao12i}: we obtain the Taylor expansion around zero of the $g_{\dot\gamma}$-norm of a Jacobi field $J$ along a lightlike geodesic $\gamma$ (Lemma~\ref{JacobiTaylor}) and the null flag curvature shows up as the coefficient of order 4, up to multiplying by $\tfrac{-3}{g_{\dot\gamma(0)}(J^{\prime}(0),J^{\prime}(0))}$.

The null flag curvature is well defined in any pseudo-Finsler space with a lightlike cone, not necessarily one arising from a Zermelo problem. Moreover, lightlike geodesics are invariant, up to reparametrization, under anisotropic conformal changes: if $L$ and $\til{L} = \varphi L$ are pseudo-Finsler metrics defined on $A\subset TM$, with $\varphi$ a function on $A$, then they have the same lightlike geodesics up to reparametrization (see \cite{jav20a}). It is therefore natural to ask whether the null flag curvatures of a lightlike pre-geodesic $\gamma$, computed with $L$ and with $\til{L}$, are related. We answer this question affirmatively in Theorem~\ref{AnisotropicConformalFlagCurv}. Since the anisotropic dependence of $\varphi$ makes it extremely difficult to relate the Chern connections of $L$ and $\til{L}$, we rely instead on the Taylor expansion obtained in Section~\ref{nullflag}.

Our most notable application of this relation concerns the flag cur\-va\-ture of Randers metrics. Given a time-independent Randers metric $F$, the Lorentz--Finsler metric \eqref{LFintro} is anisotropically conformal to a stationary spacetime. Computing the null flag curvature of this stationary spacetime, we obtain expressions for the flag curvature of $F$ in terms of both the $(\alpha,\beta)$ data (Theorem~\ref{Randerssub}) and the Zermelo data $(h,W)$ (Theorem~\ref{ZermeloFlagCurvature}). The formula for Zermelo data is especially noteworthy: to the best of our knowledge, no such formula is available for an arbitrary wind $W$ (see \cite{bao04z,hua15f} for special winds, such as Killing or conformal ones), and it opens a new path towards the understanding of special Randers metrics. Formulas for $(\alpha,\beta)$ data have appeared before \cite[Lemma~2.1]{lis09r}, but ours is considerably more compact. Moreover, these expressions allow us to compute the flag curvature of Randers metrics in several relevant cases: when the vector field $W$ in the Zermelo data is conformal, affine or projective (Corollaries~\ref{conformalW} and \ref{projectiveW}), and when $\alpha$ is a K\"ahler metric and $\beta$ is built from its potential function $\Phi$ and its complex structure $J$ (Corollary~\ref{Kahler}).

Recall that the classification of Randers metrics of constant flag curvature in \cite{bao04z} was achieved by reformulating earlier results in terms of Zermelo data, and that many of the computations there were carried out with computer assistance. Our method, which exploits a relation with conformal anisotropy, allows us to obtain explicit formulas for the flag curvature in terms of Zermelo data. We hope that these formulas will help advance the classification of special classes of Randers metrics, such as those of scalar flag curvature. In this regard, the work \cite{che12n} is particularly relevant: its authors give an example of a Randers metric of scalar flag curvature, which we recover in Corollary~\ref{Kahler}. For further references on Randers metrics of scalar flag curvature, see \cite{chen20c,she13c}.

The conformal relationship has further consequences. We define a Weyl-type curvature that is invariant under anisotropic conformal changes (Definition~\ref{weyldef}) and whose vanishing characterizes Finsler manifolds of scalar flag curvature (Proposition~\ref{nullW}). We also compute the flag curvature of the Finsler metric $F + \d f$, where $F$ is a Finsler metric and $f$ is a smooth function such that $F + \d f > 0$. This is possible because the spacetimes \eqref{LFintro} associated with $F$ and with $F + \d f$ are anisotropically conformal. In particular, the formula shows that $F$ has scalar flag curvature if and only if $F + \d f$ does.

The paper is organized as follows. In Section~\ref{preliminaries} we recall some preliminaries on Finsler geometry, focusing on anisotropic tensor calculus, which is the main tool throughout the paper. In Section~\ref{Zermelo} we briefly review the Zermelo problem and its relation with Lorentz--Finsler spacetimes. In Section~\ref{nullflag} we define the null flag curvature and show that it measures how lightlike geodesics separate (Lemma~\ref{JacobiTaylor}). In Section~\ref{conformalflag} we obtain the formula relating the null flag curvatures of two anisotropically conformal pseudo-Finsler metrics (Theorem~\ref{AnisotropicConformalFlagCurv}). In Section~\ref{weylsection} we define the Weyl-type curvature, prove that it is an anisotropic conformal invariant (Corollary~\ref{weylinv}) and show that its vanishing characterizes Finsler manifolds of scalar flag curvature. In Section~\ref{associatespacetime} we prove that, when $F$ is independent of time, the null flag curvature $\K_{(1,v)}(\tau,w)$ of $L$ in \eqref{LFintro} for a flagpole $(1,v)$ with $F(v)=1$ coincides with the flag curvature $K_v(w)$ of $F$ for the flagpole $v$. In Section~\ref{samecone} we combine this fact with the anisotropic conformal relation between the spacetimes associated with $F$ and $F + \d f$ to obtain a formula for the flag curvature of $F + \d f$ (Proposition~\ref{F+dfProp}). In Sections~\ref{Randersalphabeta} and \ref{RandersZermelo} we obtain formulas for the flag curvature of a Randers metric, using that its associated spacetime is anisotropically conformal to a classical stationary spacetime, both in terms of $(\alpha,\beta)$ data (Theorem~\ref{Randerssub}) and of Zermelo data $(h,W)$ (Theorem~\ref{ZermeloFlagCurvature}); in both cases we compute the flag curvature explicitly for specific data, for instance when $W$ is projective or conformal, or when $\alpha$ comes from a K\"ahler manifold. Finally, in Appendices~\ref{stationarysection} and \ref{stdstationarysection} we compute the null curvature of a stationary spacetime for two different expressions of the stationary metric, respectively, which is needed for the formulas for Randers metrics.
%--------------------------------------------------------------------
\section{Preliminaries on Finsler geometry}
\label{preliminaries}
%--------------------------------------------------------------------

%--------------------------------------------------------------------
\subsection{Finsler metric and indicatrix}
\label{Indicatrix}
%--------------------------------------------------------------------

Finsler metrics arise effortlessly when generalizing the concept of path length in a smooth manifold $M$. Indeed, to measure length in $M$, the most natural approach is to assign a positive value $F(v)$ to each vector $v$ tangent to a regular curve, thus setting the length of $\gamma: [a,b] \rightarrow M$ as
%--------------------------------------------------------------------
\begin{equation}
\label{Length}
\mathscr{L}(\gamma) = \int_a^b F(\dot\gamma (s)) \, \d s .
\end{equation}
%--------------------------------------------------------------------
For this notion to be independent of the parametrization, $F$ must satisfy positive homogeneity. Being positively, rather than absolutely, homogeneous implies that a parametrized curve and a reparametrization reversing its direction may yield different results. The (oriented constant speed) critical curves for the length structure \eqref{Length} are called {\em geodesics}. To ensure geodesics exist and locally minimize length, an additional (and stronger) condition is usually required; namely, for each nonzero tangent vector $v$, the fundamental tensor
%--------------------------------------------------------------------
\begin{equation}
\label{FundTensorF}
g_v (u,w) = \left. \frac{\partial^2}{\partial t \partial s}  \frac{1}{2} F^2(v + tu + sw) \right\vert_{t = s = 0}
\end{equation}
%--------------------------------------------------------------------
must be a positive-definite symmetric bilinear form on $T_{\pi(v)} M$. Under these hypotheses, a smooth function $F$ on the slit tangent bundle $TM \setminus \mathbf{0}$ will be called a {\em Finsler metric}.

Alternatively, a Finsler metric $F$ is also determined by its indicatrix
%--------------------------------------------------------------------
\begin{equation*}
\Sigma = \{ v \in TM : F(v) = 1 \} .
\end{equation*}
%--------------------------------------------------------------------
In fact, if $\Sigma_p$ is a hypersurface of $T_p M$  such that each half-ray departing from the origin intersects it exactly once, then, for each $v \in T_p M \setminus 0$, there exists $\lambda > 0$ such that $\lambda v \in \Sigma_p$, in which case, $F(v) = 1/\lambda$. When $\Sigma_p$  encloses the origin in its bounded region and has positive curvature with respect to an arbitrary Euclidean metric in the tangent space, it always determines a Minkowski norm on $T_pM$. If, in addition, $\Sigma = \cup_{p\in M}\Sigma_p$ is a smooth hypersurface of $TM$ that is transverse to the vertical spaces $T_{v}(T_pM)$, it defines a Finsler metric. Positive curvature is equivalent to the fundamental tensor \eqref{FundTensorF} being positive-definite. For more details, see e.g. \cite{jav14d}.

%--------------------------------------------------------------------
\subsection{Pseudo-Finsler metric and cone structure}
%--------------------------------------------------------------------

We will need to consider  more general families of Finsler metrics by weakening the positive-definiteness of \eqref{FundTensorF}. In general, these families are known as pseudo-Finsler metrics and have a special role in the extensions of General Relativity to the Finslerian realm, e.g. \cite{asa12f, pfe12f}.

A subset $A$ of $TM \setminus \mathbf{0}$ is said to be {\em conic} when $\lambda v \in A$, $\forall v \in A$, and $\lambda > 0$. Given a conic open subset $A$ of $TM \setminus \mathbf{0}$ such that $\pi(A) = M$, a smooth function $L: A \subset TM \setminus \mathbf{0} \to \mathbb{R}$ is a {\em (conic) pseudo-Finsler metric} if:
%--------------------------------------------------------------------
\begin{enumerate}[label=(\roman*)]
\item it is positively homogeneous of degree $2$, i.e. $L(\lambda v) = \lambda^2 L(v)$, $\forall v \in A$, $\lambda > 0$; and,

\item for each $v \in A$, the {\em fundamental tensor} defined as
%--------------------------------------------------------------------
\begin{equation}
\label{FundTensorL}
g_v(u,w) = \left. \frac{\partial^2}{\partial t \partial s} \frac{1}{2} L(v + tu + sw) \right\vert_{t = s = 0}
\end{equation}
%--------------------------------------------------------------------
is a nondegenerate symmetric bilinear form on $T_{\pi(v)}M$.
\end{enumerate}
%--------------------------------------------------------------------
Whenever the domain of $L$ is $TM \setminus \mathbf{0}$, the word ``conic'' is completely omitted from its name, and $L$ is extended continuously to $TM$.

By homogeneity, a (conic) pseudo-Finsler metric $L$ satisfies
%--------------------------------------------------------------------
\begin{subequations}
\begin{equation}
\label{LhomogeneousEq}
L(v) = g_v(v,v) , \text{ and }
\end{equation}
\begin{equation}
\label{dLhomogeneousEq}
\d L_v(w) = 2 g_v(v,w) ,
\end{equation}
\end{subequations}
%--------------------------------------------------------------------
for each $v \in A$ and all $w \in T_{\pi(v)}M$. Another simple consequence of homogeneity is the {\em cone structure} or {\em lightcone} associated with $L$, i.e. the conic subset of $A$ defined by  $\C = \{ v \in A : L(v) = 0 \}$. Since $0$ is a regular value of $L$, the cone structure (when non-empty) is a smooth hypersurface that is transverse to the vertical fibers $T_v(T_{\pi(v)}M)$; see \cite[Lemma~2.3]{jav20a}. Moreover, by \eqref{dLhomogeneousEq}, the tangent space to $\C$ at $v$ is
%--------------------------------------------------------------------
\begin{equation}
\label{ConeTangentSp}
T_v\C = \{ w \in T_{\pi(v)}M : g_v(v,w) = 0 \} = \ker(\ell_v) ,
\end{equation}
%--------------------------------------------------------------------
where $\ell_v = g_v(v,\cdot)$, and $\ell: A \subset TM\setminus \mathbf{0} \to TM^{\ast} \setminus \mathbf{0} : v \mapsto \ell_v$ is the Legendre transform of $L$. Here, as usual, we identify the vertical space $T_v(T_{\pi(v)}M)$ with $T_{\pi(v)}M$, so $T_v\C$ denotes the tangent space at $v$ of the fiber $\C\cap T_{\pi(v)}M$. The orthogonal subspace to $\ker(\ell_v)$ in $(T_{\pi(v)}M,g_v)$ is
%--------------------------------------------------------------------
\begin{equation}
\label{ConeOrt}
\ker(\ell_v)^{\perp} = \spn{v} .
\end{equation}
%--------------------------------------------------------------------
Indeed, $\ker(\ell_v)^{\perp} = \{ u \in T_{\pi(v)}M : g_v(u,w) = 0 , \forall w \in \ker(\ell_v) \}$. So $\spn{v} \subset \ker(\ell_v)^{\perp}$ is immediate, and equality follows from dimension comparison: $\dim\ker(\ell_v)^{\perp} = 1$; e.g. \cite[Lemma~2.22]{one83s}.

%--------------------------------------------------------------------
\subsection{\texorpdfstring{$A$}{A}-Anisotropic tensors}
%--------------------------------------------------------------------

Let $A$ be an open subset of $TM \setminus \mathbf{0}$ with $\pi(A) = M$. First, recall that an {\em $A$-anisotropic vector field} in $M$ is a smooth map ${\mathcal X}: A \to TM$ such that ${\mathcal X}_v := {\mathcal X}(v) \in T_{\pi(v)}M$, and the space of $A$-anisotropic vector fields is denoted by $\T^1_0 (M_A)$. Analogously, an {\em $A$-anisotropic one-form} in $M$ is a smooth map $\vartheta: A \to TM^\ast$ such that $\vartheta_v := \vartheta(v) \in T_{\pi(v)} M^\ast$, and the space of $A$-anisotropic one-forms is denoted by $\T^0_1 (M_A)$.

%--------------------------------------------------------------------
\begin{definition}
Given $r,s \in \N \cup \{0\}$, with $r + s > 0$, an {\em $A$-anisotropic tensor $T$ of type $(r,s)$} is an $\F(A)$-multilinear map 
%--------------------------------------------------------------------
\begin{equation*}
T: \T^0_1 (M_A)^r \times \T^1_0 (M_A)^s \to \F(A) ,
\end{equation*}
%--------------------------------------------------------------------
where $\F(A)$ is the set of smooth real functions on $A$, $f: A \to \R$.
\end{definition}
%--------------------------------------------------------------------
The space of $A$-anisotropic tensors of type $(r,s)$ is denoted as $\T^r_s (M_A)$ and, by convention, $\T^0_0 (M_A) \equiv \F(A)$. For instance, the fundamental tensor $g$ introduced by \eqref{FundTensorL} is an $A$-anisotropic tensor of type $(0,2)$.

The $\F(A)$-multilinearity implies that for every $v \in A$, $T$ determines a multilinear map
%--------------------------------------------------------------------
\begin{equation*}
T_v: (T_{\pi(v)}M^\ast)^r \times (T_{\pi(v)}M)^s \to \R .
\end{equation*}
%--------------------------------------------------------------------
Moreover, since any $A$-anisotropic vector field or one-form can be locally expressed as an $\F(A)$-linear combination of the usual vector fields or one-forms, whose spaces are denoted respectively by $\X(M)$ and $\X^\ast(M)$, it follows that the map
%--------------------------------------------------------------------
\begin{equation}
\label{TensorReduction}
T: \X^\ast(M)^r \times \X(M)^s \to \F(A)
\end{equation}
%--------------------------------------------------------------------
suffices to determine the $A$-anisotropic tensor $T\in \T^r_s (M_A)$.

As in classical tensor calculus, an $\F(A)$-multilinear map
%--------------------------------------------------------------------
\begin{equation}
\label{TensorEndo}
T: \T^1_0 (M_A)^s \to \T^1_0 (M_A)
\end{equation}
%--------------------------------------------------------------------
corresponds to an $A$-anisotropic tensor of type $(1,s)$,
%--------------------------------------------------------------------
\begin{equation*}
\hor{T}: \T^0_1 (M_A) \times \T^1_0 (M_A)^s \to \F(A) ,
\end{equation*}
%--------------------------------------------------------------------
by means of
%--------------------------------------------------------------------
\begin{equation}
\label{Evaluation}
\hor{T}(\theta, X_1, \ldots, X_s) = \theta( T (X_1, \ldots, X_s) ) ,
\end{equation}
%--------------------------------------------------------------------
for arbitrary $(\theta, X_1, \ldots, X_s) \in \X^*(M) \times \X(M)^s$, in view of \eqref{TensorReduction}.

A vector field $V$ defined on an open subset $\Omega \subset M$ is said to be {\em $A$-admissible} if $V_p\in A$ for every $p\in \Omega$. In such a case, for each $T \in \T^r_s(M_A)$, one can consider a classical tensor field
%--------------------------------------------------------------------
\begin{equation*}
T_V: \X^\ast(\Omega)^r \times \X(\Omega)^s \to \F(\Omega) ,
\end{equation*}
%--------------------------------------------------------------------
defined as
%--------------------------------------------------------------------
\begin{equation*}
T_V (\theta^1, \ldots, \theta^r, X_1, \ldots, X_s) (p) = T_{V_p}(\theta^1, \ldots, \theta^r, X_1, \ldots, X_s) ,
\end{equation*}
%--------------------------------------------------------------------
where $\X(\Omega)$ and $\X^\ast(\Omega)$ denote, respectively, the vector fields and the one-forms on $\Omega$, and $\F(\Omega)$ is the set of smooth real functions on $\Omega$.

Finally, since $A$-anisotropic tensors are defined in the vertical fibers of $TM$, it is necessary to define derivatives on this vertical bundle.
%--------------------------------------------------------------------
\begin{definition}
\label{VertDerivative}
Given an $A$-anisotropic tensor $T \in \T^{r}_s (M_A)$, its {\em vertical derivative} is the tensor $\vd T \in \T^{r}_{s+1} (M_A)$ given by
%--------------------------------------------------------------------
\begin{multline*}
(\vd T)_v (\theta^1, \ldots, \theta^r, X_1, \ldots, X_s, Y) = \\
= \left. \frac{\d}{\d t} T_{v+tY_{\pi(v)}}(\theta^1, \ldots, \theta^r, X_1, \ldots, X_s)\right\vert_{t=0} ,
\end{multline*}
%--------------------------------------------------------------------
for any $v \in A$ and $(\theta^1, \ldots, \theta^r, X_1, \ldots, X_s, Y) \in \X^\ast (M)^{r} \times \X (M)^{s+1}$. By use of \eqref{Evaluation}, an analogous definition is made for the vertical derivative of anisotropic tensors of the type \eqref{TensorEndo}.
\end{definition}
%--------------------------------------------------------------------

%--------------------------------------------------------------------
\subsection{\texorpdfstring{$A$}{A}-anisotropic connections}
%--------------------------------------------------------------------

The notion of connection extends to the anisotropic realm through, basically, the same properties as classical connections. Perhaps the main difference is that $\nabla_X Y \in \T^1_0(M_A)$, $\forall X,Y \in \X(M)$. In other words, $\nabla_X Y$ is an $A$-anisotropic vector field and, therefore, $A\ni v \mapsto \nabla^v_X Y \in TM$ is a smooth map with $\nabla^v_X Y \in T_{\pi(v)}M$.

In this context, the {\em torsion} of $\nabla$ is the $A$-anisotropic tensor defined by
%--------------------------------------------------------------------
\begin{equation}
\label{Torsion}
\mathcal{T}_v (X,Y) = \nabla^v_X Y - \nabla^v_Y X - [X,Y] , \; \forall X,Y \in \X(M) ;
\end{equation}
%--------------------------------------------------------------------
recall~\eqref{TensorReduction} and~\eqref{Evaluation}. As usual, an $A$-anisotropic connection is {\em torsion-free} if $\mathcal{T}=0$. Moreover, given an $A$-admissible vector field $V$ over some open neighborhood $\Omega \subset M$, there is an affine connection $\nabla^V$ on $\Omega$, defined as
%--------------------------------------------------------------------
\begin{equation*}
(\nabla^V_X Y)_p = \nabla^{V_p}_X Y , \; \forall X,Y\in \X(\Omega) .
\end{equation*}
%--------------------------------------------------------------------

An anisotropic connection $\nabla$ induces, for each vector field $X \in \X (M)$, an {\em anisotropic tensor derivation} $\nabla_X$ in the space of anisotropic tensors $\T^r_s (M_A)$; see \cite[\S~2.2]{jav19a} for the general definition. Briefly, this can be done as in the classical case by a product rule, but one key point is the derivation of $A$-anisotropic functions. It turns out that, for any function $h \in \F(A)$, $\nabla_X h \in \F (A)$ is given by
%--------------------------------------------------------------------
\begin{equation}
\label{DerivAFunct}
\nabla_X h (v) = X(h(V))(\pi(v)) - (\vd h)_v(\nabla_X^v V) ,
\end{equation}
%--------------------------------------------------------------------
where $V$ is any $A$-admissible vector field extending $v$, i.e. $V_{\pi(v)}=v$. The expression in \eqref{DerivAFunct} does not depend on the choice of $V$; refer to \cite[Lemma~9]{jav19a}. Furthermore, given any vector field $X \in \X(M)$, an $A$-admissible vector field $V \in \X(\Omega)$ in some open subset $\Omega \subset M$, and $T \in \T^r_s (M_A)$, the derivative tensor $\nabla_X T$ can be computed using the affine connection $\nabla^V$, the classical tensor $T_V$, and the vertical derivative of $T$; namely,
%--------------------------------------------------------------------
\begin{multline}
\label{DerivATen}
(\nabla_X T)_v (\theta^1, \ldots, \theta^r, X_1, \ldots, X_s) =\\
= \nabla_X^{V} (T_V)_{\pi(v)} (\theta^1, \ldots, \theta^r, X_1, \ldots, X_s) \\
- (\vd T)_v (\theta^1, \ldots, \theta^r, X_1, \ldots, X_s, \nabla^V_X V) ,
\end{multline}
%--------------------------------------------------------------------
where $V$ is again an $A$-admissible extension of $v\in A$.

%--------------------------------------------------------------------
\subsection{Curvature tensor}
%--------------------------------------------------------------------

First of all, observe that, by means of $\F(A)$-linearity, the tensor derivation $\nabla_{\mathcal X}$ extends to anisotropic vector fields $\mathcal{X} \in \T_0^1(M_A)$. So each $A$-anisotropic connection $\nabla$ gives rise to a {\em cur\-vature tensor} $R: \T^1_0 (M_A)^3 \to \T^1_0 (M_A)$ by
%--------------------------------------------------------------------
\begin{equation}
\label{CurvTensor}
R_v(X,Y)Z = \nabla^v_X (\nabla_Y Z) - \nabla^v_Y (\nabla_X Z) - \nabla^v_{[X,Y]} Z ,
\end{equation}
%--------------------------------------------------------------------
for all $v \in A$ and $X,Y,Z \in \X(M)$. One can verify that $R$ is indeed an $\F(A)$-multilinear map and, hence, an anisotropic tensor of the type described in \eqref{TensorEndo}, which is also antisymmetric in $X$ and $Y$. Moreover, this curvature tensor can be computed using any $A$-admissible vector field $V$ that extends $v$ in an open subset $\Omega \subset M$. To this end, one has to turn to the vertical derivative of $\nabla$, called the {\em Chern tensor}, that is, the anisotropic tensor $P: \T^1_0 (M_A)^3 \to \T^1_0 (M_A)$ defined by
%--------------------------------------------------------------------
\begin{equation}
\label{ChernTensor}
P_v (X,Y,Z) = \left. \frac{\d}{\d t} \left( \nabla^{v+tZ_{\pi(v)}}_X Y \right) \right\vert_{t=0} ,
\end{equation}
%--------------------------------------------------------------------
for all $v \in A$ and $X,Y,Z \in \X(M)$. Then,
%--------------------------------------------------------------------
\begin{equation*}
R_v(X,Y)Z = (R^V(X,Y)Z - P_V(Y,Z,\nabla^V_X V) + P_V(X,Z,\nabla^V_Y V))_{\pi(v)} \, ,
\end{equation*}
%--------------------------------------------------------------------
where $X,Y,Z \in \X(\Omega)$ and $R^V$ is the curvature tensor of the affine connection $\nabla^V$.

%--------------------------------------------------------------------
\subsection{Chern connection and flag curvature }
\label{ChernFlag}
%--------------------------------------------------------------------

The {\em Chern connection} associated with a pseudo-Finsler metric $L$ is the unique $A$-anisotropic connection $\nabla$ that  is torsion-free (see \eqref{Torsion}) and satisfies $\nabla g=0$, where $g$ is the fundamental tensor \eqref{FundTensorL} and $\nabla g$ is the anisotropic derivative tensor computed by \eqref{DerivATen}; refer to \cite{mat80z}, \cite[(7.20) and (7.21)]{she01d} or \cite[Proposition~2.3]{jav14c}.

For every smooth vector field $V$ on some open subset $\Omega \subset M$, the affine connection $\nabla^{V}$ associated with the Chern connection satisfies:
%--------------------------------------------------------------------
\begin{enumerate}[label=(\roman*)]
\item $\nabla^V_X Y - \nabla^V_Y X = [X,Y]$, $\forall X,Y \in \X(\Omega)$;
\item $X g_V(Y,Z) = g_V(\nabla^V_X Y,Z) + g_V(Y,\nabla^V_X Z) + 2C_V(\nabla^V_X V,Y,Z)$, $\forall X,Y,Z \in \X(\Omega)$.
\end{enumerate}
%--------------------------------------------------------------------
Here, $C$ is the {\em Cartan tensor}, defined for each $v \in A$ as
%--------------------------------------------------------------------
\begin{equation*}
C_v (u,w,y) = \frac{1}{2} \vd g_v(u,w,y) = \left.\frac{\d}{\d t} \frac{1}{2} g_{v + ty} (u,w) \right\vert_{t=0} , \; \forall u,w,y \in T_{\pi(v)}M .
\end{equation*}
%--------------------------------------------------------------------
Furthermore, this connection is determined by a Koszul type formula:
%--------------------------------------------------------------------
\begin{multline}
\label{KoszulQ}
2 g_v(\nabla^V_X Y,Z) = \\ = \left( X(g_V(Y,Z)) - Z(g_V(X,Y)) + Y(g_V(Z,X)) \right)_{\pi(v)} \\
+ g_v([X,Y],Z) + g_v([Z,X],Y) - g_v([Y,Z],X) \\
+ 2\left( - C_v(Y,Z,\nabla^V_X V) - C_v(Z,X,\nabla^V_Y V) + C_v(X,Y,\nabla^V_Z V) \right),
\end{multline}
%--------------------------------------------------------------------
where $V$ is an $A$-admissible local extension of $v\in A$ and $X,Y,Z$ are arbitrary vector fields.

The curvature tensor $R$ associated with the Chern connection extends some symmetries of the usual Riemann curvature tensor. For the current work, we will need the following:
%--------------------------------------------------------------------
\begin{subequations}
\begin{equation}
\label{symmetry1}
g_v(R_v(u,w)v,v) = 0 ;
\end{equation}
\begin{equation}
\label{symmetry2}
g_v(R_v(w,v)v,u) = -g_v(R_v(w,v)u,v) ;
\end{equation}
\begin{equation}
\label{symmetry3}
g_v(R_v(u,v)w,v) = g_v(R_v(w,v)u,v) ;
\end{equation}
\end{subequations}
%--------------------------------------------------------------------
for any $v \in A$ and all $u,w \in T_{\pi(v)} M$. These identities are well-known and, in the form they are presented here, can be easily derived from \cite[Proposition~3.1]{jav20c}.

Finally, the curvature tensor provides one of the most fundamental geometric invariants of pseudo-Finsler geometry, the {\em flag curvature}, which measures how geodesics with initial vector in a certain plane (the flag) spread apart from a fixed geodesic with initial velocity $v$ (the flagpole). The formal definition is given by
%--------------------------------------------------------------------
\begin{equation}
\label{FlagCurv}
K_v(w) = \frac{g_v(R_v(w,v)v,w)}{g_v(v,v) g_v(w,w) - g_v(v,w)^2} ,
\end{equation}
%--------------------------------------------------------------------
for $v\in A$ and $w\in T_{\pi(v)}M$ satisfying $g_v(v,v) g_v(w,w) - g_v(v,w)^2 \neq 0$, i.e. the plane (or flag) $\Pi=\spn{v,w}$ is $g_v$-nondegenerate.

%--------------------------------------------------------------------
\subsection{Covariant derivative, geodesics and Jacobi fields}
\label{GeoJacobi}
%--------------------------------------------------------------------

Given a coordinate system $(\Omega,\varphi)$ on $M$, the {\em Christoffel symbols} of $\nabla$ are the functions $\Gamma_{\;ij}^k: T\Omega \cap A \to \R$ determined by
%--------------------------------------------------------------------
\begin{equation*}
\nabla^v_{\partial_i}\partial_j=\Gamma_{\;ij}^k(v)\,\partial_k ,
\end{equation*}
%--------------------------------------------------------------------
with use of the Einstein summation convention on the right-hand side. From the Christoffel symbols stems a definition of covariant derivative along any curve on the manifold, by fixing a reference vector field along the said curve. Denote the {\em covariant derivative} along $\gamma:[a,b] \to M$ with respect to an admissible $V \in \X(\gamma)$ as $D^V_\gamma: \X(\gamma) \to \X(\gamma)$. When the curve $\gamma$ has image in the chart $(\Omega,\varphi)$, then
%--------------------------------------------------------------------
\begin{equation*}
D^V_\gamma X = (\dot X^k + \dot\gamma^i X^j \Gamma^k_{\;ij} \circ V) \partial_k ,
\end{equation*}
%--------------------------------------------------------------------
where $X = X^i \partial_i$ and $\dot\gamma = \dot\gamma^i \partial_i$ in the fixed chart.

This covariant derivative is positively homogeneous of degree zero in $V$, i.e. $D_{\gamma}^{\lambda V} = D_{\gamma}^V$, $\forall \lambda > 0$. Moreover, it allows us to characterize the geodesics of the pseudo-Finsler space $(M,L)$. Precisely, a smooth admissible curve $\gamma: [a,b] \to M$ is a {\em geodesic} of $L$ when $D^{\dot\gamma}_\gamma \dot\gamma = 0$. Then a vector field $J$ along a geodesic $\gamma$ is the variational vector field of a variation by geodesics if and only if it satisfies
%--------------------------------------------------------------------
\begin{equation}
\label{JacobiEq}
J^{\prime\prime} + R_{\dot\gamma}(J,\dot\gamma)\dot\gamma = 0 ,
\end{equation}
%--------------------------------------------------------------------
where $\prime$ denotes $D^{\dot\gamma}_\gamma$. A vector field that satisfies the equation \eqref{JacobiEq} is called a {\em Jacobi field}.

%--------------------------------------------------------------------
\subsection{Lorentz--Finsler metric and lightcone}
\label{LorentzFinsler}
%--------------------------------------------------------------------

Let $M^n$ be a smooth manifold. A pseu\-do-Finsler metric $L: A \to \R$ is said to be a {\em Lorentz--Finsler metric} if
%--------------------------------------------------------------------
\begin{enumerate}[label=(\roman*)]
\item $L>0$; 
\item its fundamental tensor \eqref{FundTensorL} has index $n-1$;
\item $L$ extends smoothly to the boundary of $A$ in $TM \setminus \mathbf{0}$, where $L=0$ and the fundamental tensor $g$ is nondegenerate.
\end{enumerate}
%--------------------------------------------------------------------
We assume $L$ smoothly extends to an open cone $\mathcal{A}\subset TM \setminus \mathbf{0}$ with $\bar{A}\setminus \mathbf{0} \subset \mathcal{A}$. In these conditions, the boundary $\C$ of $A$ in $TM\setminus \mathbf{0}$ is the (non-empty) lightcone of $L$ on $\mathcal{A}$. It is a smooth hypersurface of $TM$, transversal to the vertical spaces $T_{v}(T_{\pi(v)}M)$, and $\C_p = \C \cap T_p M$ is a {\em strong cone} for every $p\in M$. To be precise, $\C_p$ is conic (i.e. if $v \in \C_p$, then $\lambda v \in \C_p$ for all $\lambda>0$), and there is an affine hyperplane $\Pi_p$ that cuts all the half-rays of $\C_p$ in a strongly convex hypersurface of $\Pi_p$; see \cite{jav20d} for more details. See also \cite{edw18r,hee27e,min19c,pfe19f} and references therein for other aspects of Finsler spacetimes.
%--------------------------------------------------------------------
\section{Zermelo navigation problem}
\label{Zermelo}
%--------------------------------------------------------------------

Consider a moving object on a manifold $M$. Suppose its maximum speed is fixed for each direction at each point, and that we want to compute the trajectory between two given points that minimizes the travel time. In addition, suppose the set of possible maximum speeds at a given point $p \in M$ forms a hypersurface $\Sigma_p \subset T_p M$ that encloses the zero vector and has positive curvature with respect to some Euclidean metric. Under smoothness and transversality assumptions, there is a (positive-definite) Finsler metric $F$ whose indicatrix at $p\in M$ is the hypersurface $\Sigma_p$; recall \S~\ref{Indicatrix}.

Since the velocity of a (regular) trajectory $c:[a,b] \rightarrow M$ is determined at every instant, we can compute how the curve $c$ is traveled with respect to absolute time. In other words, assume that the parameter $s$ depends on the absolute time $t$, namely, it can be expressed as $s(t)$. The $F$-unit vector in the direction $\dot c(s(t))$ is $\frac{1}{F(\dot c(s(t)))} \dot c(s(t))\in\Sigma_{c(s(t))}$. So the velocity satisfies
%--------------------------------------------------------------------
\begin{equation}
\label{cUnit}
\frac{\d}{\d t} c(s(t))=\frac{1}{F(\dot c(s(t)))} \dot c(s(t)) .
\end{equation}
%--------------------------------------------------------------------
On the other hand,
%--------------------------------------------------------------------
\begin{equation}
\label{cChain}
\frac{\d}{\d t} c(s(t))=\dot s(t) \dot c(s(t)) .
\end{equation}
%--------------------------------------------------------------------
From \eqref{cUnit} and \eqref{cChain}, it follows that
%--------------------------------------------------------------------
\begin{equation*}
 \dot s(t) = \frac{1}{F(\dot c(s(t)))} ,
\end{equation*}
%--------------------------------------------------------------------
and the inverse function $t(s)$ satisfies
%--------------------------------------------------------------------
\begin{equation*}
\dot t(s) = F(\dot c(s)) .
\end{equation*}
%--------------------------------------------------------------------
From the above identity,
%--------------------------------------------------------------------
\begin{equation*}
t(b) - t(a) = \int_a^b \dot{t}(s) \d s = \int_a^b F(\dot{c}(s)) \d s ,
\end{equation*}
%--------------------------------------------------------------------
that is, the elapsed time to travel along $c$ is the length of the curve measured with respect to the Finsler metric $F$. Therefore, the time-minimizing paths are also length-minimizing and, specifically, geodesics for the Finsler metric $F$. In fact, as the geodesics of a Finsler metric are locally minimizing, sufficiently small segments of a geodesic will be global minima of the length between its endpoints.

Now, what happens in the Zermelo problem when the maximum speeds vary with time? In this case, for each time $t \in \R$, there is a Finsler metric $F_t$, and the time dependence prevents us from calculating the absolute time as the length with respect to a (positive-definite) Finsler metric. Instead, given a trajectory $c(s)$, we have 
%--------------------------------------------------------------------
\begin{equation}
\label{TimeDependence}
\dot{t}(s) = F_{t(s)}(\dot{c}(s)) ,
\end{equation}
%--------------------------------------------------------------------
which is a differential equation in the variable $t(s)$. To address this issue, we introduce a non-relativistic spacetime $\R \times M$, where the first component $\R$ encodes the time coordinate $t$.

Define on $\R \times M$ the Lorentz--Finsler metric
%--------------------------------------------------------------------
\begin{equation}
\label{LFmetric}
L(\tau, v) = \tau^2 - F_t(v)^2 ,
\end{equation}
%--------------------------------------------------------------------
where $(\tau,v)$ must satisfy $ \tau  > F_t(v)$. Observe that when $v=0$, $L$ is smooth only if $F_t$ is Riemannian (i.e. if $L$ is Lorentzian), but this is not important for our purposes as we are interested in the lightlike directions. In any case, it is possible to modify $L$ in a neighborhood of the straight line $s\mapsto (s,0)$ to have a smooth Lorentz--Finsler metric (see \cite[\S 5]{jav20d} for more details).  An object travels the trajectory $c:[a,b] \rightarrow M$ at maximum speed if and only if it satisfies \eqref{TimeDependence}, which is equivalent to
%--------------------------------------------------------------------
\begin{equation*}
L(\dot{t}(s), \dot{c}(s)) = 0 ,
\end{equation*}
%--------------------------------------------------------------------
that is, the curve $s \mapsto (t(s), c(s))$ is lightlike. Therefore, finding the trajectory $c:[a,b] \rightarrow M$ between two points $p$ and $q$ of $M$ traveling at maximum speed that minimizes the arrival time, amounts to finding the lightlike curve $s \mapsto (t(s), c(s))$ connecting $(t(a), p)$ to the vertical line $s \mapsto (s, q)$ that minimizes the $t$-coordinate at the endpoint among all lightlike curves between $(t(a), p)$ and the line $s \mapsto (s, q)$. Thanks to Fermat's principle for Finsler spacetimes obtained in \cite{per06f}, we conclude the solutions are lightlike pre-geodesics of \eqref{LFmetric}, i.e. lightlike geodesics of $L$ up to parametrization. Thus, we have the following.

%--------------------------------------------------------------------
\begin{proposition}
\label{ZermeloLightlike}
The solutions of the time-dependent generalized Zermelo problem are given by the (projections of) lightlike pre-geodesics of the Lorentz--Finsler metric \eqref{LFmetric}.
\end{proposition}
%--------------------------------------------------------------------
%--------------------------------------------------------------------
\section{Null flag curvature}
\label{nullflag}
%--------------------------------------------------------------------

Recall from \S~\ref{ChernFlag} that the flag curvature of a pseudo-Finsler metric is defined only on flags $\Pi = \spn{v,w}$ such that the plane $\Pi$ is nondegenerate with respect to $g_v$, where $v$ is the flagpole. However, taking into account Proposition~\ref{ZermeloLightlike}, we need to consider lightlike geodesics of a Lorentz--Finsler metric to solve the time-dependent Zermelo problem, and variations of lightlike geodesics will span degenerate planes which are flags with respect to the velocity of the initial lightlike geodesic. As a consequence, to measure how these geodesics spread apart, we cannot use the classical flag curvature. A suitable definition for this setting is obtained following the idea for the sectional curvature of degenerate planes in the Lorentzian case \cite{har82t,pal07f,alb10g}.

%--------------------------------------------------------------------
\begin{definition}
\label{NullFlgaCurvature}
Given a pseudo-Finsler metric $L: A \subset TM \setminus \mathbf{0} \to \R$, a {\em null flag} is a degenerate plane $\Pi = \spn{v,w}$ such that $g_v(v,w)=g_v(v,v)=0$ and $g_v(w,w)\not=0$; its flagpole is the null vector $v$. Then we define the {\em null flag curvature} of the flag $\Pi$ with respect to $v$ as
%--------------------------------------------------------------------
\begin{equation}
\K_v(w) = \frac{g_v(R_v(w,v)v,w)}{g_v(w,w)} .
\end{equation}
%--------------------------------------------------------------------
\end{definition}
%--------------------------------------------------------------------

Requiring that $w$ is $g_v$-orthogonal to $v$ ensures that the degenerate plane $\Pi$ is tangent to the lightcone of $L$ at $v$. From the curvature properties \eqref{symmetry1}--\eqref{symmetry3}, $K_v(w)$ depends only on the flagpole $v$ and the degenerate plane $\Pi$; namely,
%--------------------------------------------------------------------
\begin{equation}
\label{NullFlgaCurvWellDef}
\K_v(w) = \K_v(u) , \; \forall u \in \spn{v,w} \text{ with } g_v(u,u) \neq 0 .
\end{equation}
%--------------------------------------------------------------------
Unlike the classical flag curvature, the null flag curvature is positive homogeneous of degree two in $v$, i.e. $\K_{\lambda v}(w) = \lambda^2 \K_v(w)$, $\forall \lambda > 0$. But similarly to the classical definition, the null flag curvature also measures how geodesics spread apart. To verify this property, consider a Jacobi field $J$ along a lightlike geodesic $\gamma: [0,\delta] \to M$; recall \S~\ref{GeoJacobi}. As in the case of the flag curvature, the null flag curvature will be, up to some normalization, the coefficient of order $4$ of the Taylor expansion of $g_{\dot\gamma}(J,J)$. Indeed, since $\gamma$ is a geodesic, the fundamental tensor $g_{\dot\gamma}$ is compatible with $D_{\gamma}^{\dot\gamma}$, because the Cartan tensor $C_{\dot\gamma}(D_{\gamma}^{\dot\gamma} \dot\gamma, \cdot, \cdot)$ vanishes, i.e. for arbitrary vector fields $X, Y$ along $\gamma$,
%--------------------------------------------------------------------
\begin{equation}
\label{Compatibility}
\frac{\d}{\d t} g_{\dot\gamma}(X,Y) = g_{\dot\gamma}(X^{\prime},Y) + g_{\dot\gamma}(X,Y^{\prime}) ,
\end{equation}
%--------------------------------------------------------------------
where $\prime$ denotes the covariant derivative $D^{\dot\gamma}_\gamma$. So we can repeat the computations of \cite[\S~5.5]{bao12i} adapted to the pseudo-Finsler case.

%--------------------------------------------------------------------
\begin{lemma}
\label{JacobiTaylor}
Let $J$ be a Jacobi field along a geodesic $\gamma: [0,\delta] \to M$. If $\dot\gamma(0) = v$, $J(0) = 0$, and $J^{\prime}(0) = w$, then
%--------------------------------------------------------------------
\begin{equation}
\label{ExpansionJ}
g_{\dot\gamma(t)}(J(t),J(t)) = g_v(w,w) t^2 - \frac{1}{3} g_v(R_v(w,v)v,w) t^4 + \mathcal{O}(t^5) .
\end{equation}
%--------------------------------------------------------------------
\end{lemma}
%--------------------------------------------------------------------
%--------------------------------------------------------------------
\begin{proof}
Consider the real function $f = g_{\dot\gamma}(J,J)$. In view of \eqref{Compatibility}, and the Jacobi equation \eqref{JacobiEq}, we obtain that
%--------------------------------------------------------------------
\begin{align*}
\dot f &= 2g_{\dot\gamma}(J^{\prime},J) ,\\
\ddot f &= 2g_{\dot\gamma}(J^{\prime\prime}, J) + 2g_{\dot\gamma}(J^{\prime},J^{\prime}) ,\\
f^{(3)} &= 2g_{\dot\gamma}(J^{(3)}, J) + 6g_{\dot\gamma}(J^{\prime\prime},J^{\prime}) \\
    &= 2g_{\dot\gamma}(J^{(3)}, J) - 6g_{\dot\gamma}(R_{\dot\gamma} (J,\dot\gamma)\dot\gamma,J^{\prime}) ,\\
f^{(4)} &= 2g_{\dot\gamma}(J^{(4)}, J) + 8g_{\dot\gamma}(J^{(3)}, J^{\prime}) + 6g_{\dot\gamma}(J^{\prime\prime}, J^{\prime\prime}) \\
    &= 2g_{\dot\gamma}(J^{(4)}, J) - 8g_{\dot\gamma}(D^{\dot\gamma}_\gamma (R_{\dot\gamma} (J,\dot\gamma)\dot\gamma), J^{\prime}) \\
    &\pushright{+ 6g_{\dot\gamma}(R_{\dot\gamma} (J,\dot\gamma)\dot\gamma, R_{\dot\gamma} (J,\dot\gamma)\dot\gamma) .}
\end{align*}
%--------------------------------------------------------------------
By the anisotropic tensor derivation \eqref{DerivATen}, and the geodesic equation $D^{\dot\gamma}_\gamma\dot\gamma=0$, it follows that $D^{\dot\gamma}_\gamma (R_{\dot\gamma} (J,\dot\gamma)\dot\gamma) =(\nabla^{\dot\gamma}_{\dot\gamma} R)_{\dot\gamma} (J,\dot\gamma)\dot\gamma + R_{\dot\gamma}(J^\prime, \dot\gamma)\dot\gamma$. Hence,
%--------------------------------------------------------------------
\begin{align*}
\dot f(0) &= 0 ,\\
\ddot f(0) &= 2g_{v}(w,w) ,\\
f^{(3)}(0) &= 0 , \\
f^{(4)}(0) &= -8g_{v}(R_{v} (w,v)v, w) ,
\end{align*}
%--------------------------------------------------------------------
which implies \eqref{ExpansionJ}.
\end{proof}
%--------------------------------------------------------------------

Let us conclude this section with another consequence of the metric compatibility of $g_{\dot\gamma}$ with $D_{\gamma}^{\dot\gamma}$, that will be subsequently necessary.

%--------------------------------------------------------------------
\begin{lemma}
\label{AffJangle}
If $J$ is a Jacobi field along a geodesic $\gamma$ with initial conditions $\dot\gamma(0) = v$, $J(0) = u$, and $J^{\prime}(0) = w$, then
%--------------------------------------------------------------------
\begin{equation}
\label{Jangle}
g_{\dot\gamma(t)}(\dot\gamma(t),J(t)) = g_v(v,u) + t g_v(v,w) .
\end{equation}
%--------------------------------------------------------------------
\end{lemma}
%--------------------------------------------------------------------
%--------------------------------------------------------------------
\begin{proof}
From \eqref{Compatibility}, and the Jacobi equation \eqref{JacobiEq}, the real function $f = g_{\dot\gamma}(\dot\gamma,J)$ satisfies
%--------------------------------------------------------------------
\begin{align*}
\dot f &= g_{\dot\gamma}(\dot\gamma , J^{\prime}) , \\
\ddot f &= g_{\dot\gamma}(\dot\gamma , J^{\prime\prime}) = - g_{\dot\gamma}(\dot\gamma , R_{\dot\gamma}(J , \dot\gamma) \dot\gamma) .
\end{align*}
%--------------------------------------------------------------------
By the symmetry \eqref{symmetry1} of the curvature tensor, the second derivative of $f$ vanishes, thus \eqref{Jangle}.
\end{proof}
%--------------------------------------------------------------------
%--------------------------------------------------------------------
\section{Null flag curvature of anisotropically conformal metrics}
\label{conformalflag}
%--------------------------------------------------------------------

Two pseudo-Finsler metrics $L,\til{L}$ defined on $A \subset TM \setminus \mathbf{0}$ are {\em anisotropically conformal} when there exists a smooth, never-zero, and positive homogeneous of degree zero function $\varphi : A \subset TM \setminus \mathbf{0} \to \R$ such that $\til{L} (v) = \varphi(v) L(v)$. No generality is lost by assuming $\varphi>0$. Indeed, $-\til{L}$ has the same Chern connection and the same null flag curvature as $\til{L}$, and \eqref{ConformalFlagCurvEq} below is invariant under $\varphi\mapsto-\varphi$.

We settle upon the notation $\sim$ for the geometric quantities associated with the metric $\til{L}$, conformal to $L$ by the anisotropic factor $\varphi$. Their fundamental tensors can be related by a straightforward calculation from the definition~\eqref{FundTensorL} in view of the chain rule.

%--------------------------------------------------------------------
\begin{lemma}
\label{ConformalFundTensorLemma}
If $g$ and $\til{g}$ are the fundamental tensors of $(M,L)$ and $(M,\varphi L)$, respectively, then
%--------------------------------------------------------------------
\begin{multline}
\label{ConformalFundTensorEq}
\til{g}_v (u,w) = \varphi (v) g_v(u,w) + (\vd\varphi)_v(u) g_v(v,w) + (\vd\varphi)_v(w) g_v(v,u) \\
+ \frac{1}{2} (\hs \varphi)_v (u,w) L(v) ,
\end{multline}
%--------------------------------------------------------------------
where $(\vd\varphi)$ denotes the vertical derivative of $\varphi$ and $(\hs \varphi)_v$ the vertical Hessian of $\varphi$ at $v$. In particular, when $v$ is lightlike for $L$ and $g_v(v,u)=g_v(v,w)=0$, \eqref{ConformalFundTensorEq} reduces to
%--------------------------------------------------------------------
\begin{equation}
\label{v-Orthogonalities}
\til{g}_v (u,w) = \varphi(v) g_v(u,w) .
\end{equation}
%--------------------------------------------------------------------
\end{lemma}
%--------------------------------------------------------------------

It is well-known that two pseudo-Finsler metrics are anisotropically conformal if and only if they share the same lightlike cone; see \cite[Theorem~3.11]{jav20d} or \cite[Theorem~2.5]{jav20a}. Furthermore, the lightlike geodesics of two anisotropically conformal pseudo-Finsler metrics coincide up to parametrization;  e.g. \cite[Proposition~3.4]{jav20a}. Therefore, we can compare their null flag curvatures, by considering variations of lightlike geodesics; recall Lemma~\ref{JacobiTaylor}.

First, we need to specify the reparametrization obtained in \cite{jav20a}. For simplicity, we assume the conformal factor $\varphi$ is positive and we fix zero as the initial parameter of curves throughout the text. Let $\til{\gamma}: [0,\varepsilon] \to M$ be a lightlike geodesic for $\til{L}$. Then $\gamma = \til{\gamma} \circ \psi$ is a lightlike geodesic for $L$, where $\psi: [0,\delta] \to [0,\varepsilon]$ is a solution of the differential equation $\dot\psi(t) = \varphi(\dot{\til{\gamma}}(\psi(t)))$, or
%--------------------------------------------------------------------
\begin{equation}
\label{PsiPhi}
\dot\psi(t) = \varphi(\dot\gamma (t)) ,
\end{equation}
%--------------------------------------------------------------------
because $\varphi$ is positive homogeneous of degree zero; see \cite[Proposition~3.4]{jav20a} for the proof, taking into account our different notation. If $\sigma := \psi^{-1}$, then $\til{\gamma}= \gamma \circ \sigma$ and, by \eqref{PsiPhi},
%--------------------------------------------------------------------
\begin{equation}
\label{SigmaPhi}
\dot\sigma(s)=\frac{1}{\dot\psi(\sigma(s))} = \frac{1}{\varphi(\dot\gamma(\sigma(s)))} = \frac{1}{\varphi(\dot{\til{\gamma}}(s))} ,
\end{equation}
%--------------------------------------------------------------------
where we use the positive homogeneity of $\varphi$ once again. As a consequence,
%--------------------------------------------------------------------
\begin{equation}
\label{GammaTil}
\dot{\til{\gamma}}(s) = \dot\sigma(s) \dot\gamma(\sigma (s)) = \frac{1}{\varphi(\dot\gamma(\sigma(s)))} \dot\gamma(\sigma(s)) = \frac{1}{\varphi(\dot{\til{\gamma}}(s))} \dot\gamma(\sigma(s)) .
\end{equation}
%--------------------------------------------------------------------
In fact, the opposite also follows from \cite[Proposition~3.4]{jav20a}. If $\gamma$ is a lightlike geodesic for $L$, then there exists $\sigma$ satisfying \eqref{SigmaPhi} such that $\til{\gamma} = \gamma\circ\sigma$ is a lightlike geodesic for $\til{L}$.

Next, we must relate a Jacobi field for $L$ along $\gamma$ with a Jacobi field for $\til{L}$ along the reparametrized geodesic $\til{\gamma} = \gamma \circ \sigma$.

%--------------------------------------------------------------------
\begin{lemma}
\label{JTilLemma}
Let $D$ and $\til{D}$ denote the covariant derivatives associated with the Chern connections of $L$ and $\til{L} = \varphi L$, respectively. Let $J$ be a Jacobi field along a lightlike geodesic $\gamma$ for $(M,L)$ satisfying $J(0) =0$ and $J^{\prime}(0) = D_{\gamma}^{\dot\gamma}J (0)$ $g_{\dot\gamma(0)}$-orthogonal to $\gamma$. There exists a (smooth) function $h : [0,\varepsilon] \to \R$ with $h(0)=0$ such that
%--------------------------------------------------------------------
\begin{equation}
\label{JTilEq}
\til{J}(s) = J(\sigma(s)) + h(s) \dot{\til{\gamma}}(s)
\end{equation}
%--------------------------------------------------------------------
is a Jacobi field along the lightlike geodesic $\til{\gamma} = \gamma \circ \sigma$ for $(M,\til{L})$, and it satisfies
%--------------------------------------------------------------------
\begin{equation}
\label{JTilDeriv0}
\til{J}^{\prime}(0) = \til{D}_{\til{\gamma}}^{\dot{\til{\gamma}}} \til{J}(0) = \dot\sigma(0) J^{\prime}(0) + \dot h(0) \dot{\til{\gamma}}(0) .
\end{equation}
%--------------------------------------------------------------------
\end{lemma}
%--------------------------------------------------------------------
%--------------------------------------------------------------------
\begin{proof}
Assume $J$ is a Jacobi field along $\gamma$ in $(M,L)$ such that $J(0) = 0$ and $g_v(v,w)=0$, where $v = \dot\gamma(0)$ and $w = J^{\prime}(0)$. Consider a variation of $\gamma$ by lightlike geodesics with variational vector field $J$ (which is possible because $g_v(v,w)=0$) and then reparametrize these geodesics to be lightlike geodesics for $\til{L}$. As this reparametrization depends smoothly on the parameter of the variation, it is straightforward to check that its variational vector field $\til{J}$ satisfies \eqref{JTilEq}.

By \cite[Remark~3.3]{jav20a}, if $\til{V} = V \circ \sigma$ and $\til{W} = W \circ \sigma$, then
%--------------------------------------------------------------------
\begin{equation}
\label{DerivRmk}
\til{D}^{\til{V}}_{\til{\gamma}} \til{W}(s) = \dot\sigma(s) \til{D}^{V}_{\gamma} W(\sigma(s)) .
\end{equation}
%--------------------------------------------------------------------
Consider $V(t) = \dot\gamma(t)$ and $W(t) = J(t)$. From~\eqref{GammaTil} and~\eqref{JTilEq}, we have
%--------------------------------------------------------------------
\begin{align*}
\til{V}(s) &= \dot\gamma(\sigma(s)) = \varphi(\dot{\til{\gamma}}(s)) \dot{\til{\gamma}}(s) , \\
\til{W}(s) &= J(\sigma(s)) = \til{J}(s) - h(s) \dot{\til{\gamma}}(s) .
\end{align*}
%--------------------------------------------------------------------
Taking into account that $\til{D}^{\dot{\til{\gamma}}}_{\til{\gamma}} \dot{\til{\gamma}} = 0$, because $\til{\gamma}$ is geodesic in $(M,\til{L})$, \eqref{DerivRmk} yields that
%--------------------------------------------------------------------
\begin{equation}
\label{JTilDeriv}
\begin{aligned}
\til{D}^{\dot{\til{\gamma}}}_{\til{\gamma}} \til{J}(s) &= \til{D}^{\varphi(\dot{\til{\gamma}})\dot{\til{\gamma}}}_{\til{\gamma}} \til{J}(s) = \til{D}^{\dot\gamma\circ\sigma}_{\til{\gamma}} (J\circ\sigma)(s) + \dot h(s) \dot{\til{\gamma}}(s) \\
&= \dot\sigma(s) \til{D}^{\dot\gamma}_{\gamma} J(\sigma(s)) + \dot h(s) \dot{\til{\gamma}}(s) .
\end{aligned}
\end{equation}
%--------------------------------------------------------------------
Since $J$ is the variational vector field for a variation of $\gamma$ by curves, we can apply \cite[Proposition~3.2]{jav14c} to get $\til{D}_\gamma^{\dot\gamma} J(0) = \til{D}_\beta^{\dot\gamma} \dot\gamma (0)$, where $\beta$ denotes the transversal curve of the variation through $p = \gamma(0)$. As the variation fixes the starting point, we have $\dot\beta(0) = 0$ and, by direct computation, $\til{D}_\beta^{\dot\gamma} \dot\gamma(0) = D_\beta^{\dot\gamma} \dot\gamma(0)$ (it does not depend on the covariant derivative). Using again \cite[Proposition~3.2]{jav14c}, $D_\beta^{\dot\gamma} \dot\gamma (0) = D_\gamma^{\dot\gamma} J(0)$ for the other covariant derivative. Putting all this together, we conclude that
%--------------------------------------------------------------------
\begin{equation*}
\til{D}_\gamma^{\dot\gamma} J(0) = \til{D}_\beta^{\dot\gamma} \dot\gamma(0) = D_\beta^{\dot\gamma} \dot\gamma(0) = D_\gamma^{\dot\gamma} J(0) .
\end{equation*}
%--------------------------------------------------------------------
Substituting the last identity into \eqref{JTilDeriv} for $s=0$ gives  \eqref{JTilDeriv0}.
\end{proof}
%--------------------------------------------------------------------

We are now ready for the curvature comparison. We will use the notation $\nabla^2\varphi$ for the anisotropic derivative $\nabla T\in \T^0_2 (M_A)$, where $T = \nabla\varphi\in \T^0_1 (M_A)$, and the evaluation $\nabla^2 \varphi (v)$ denotes $\nabla T (v,v)$.

%--------------------------------------------------------------------
\begin{theorem}
\label{AnisotropicConformalFlagCurv}
Let $(M,L)$ and $(M,\til{L})$ be two anisotropically conformal pseudo-Finsler spaces with $\til{L}=\varphi L$ and let $\spn{v,w}$ be a null flag with flagpole $v$. If $\K_v(w)$ is the null flag curvature for $L$ and $\til{\K}_v(w)$ is the null flag curvature for $\til{L}$, then
%--------------------------------------------------------------------
\begin{equation}
\label{ConformalFlagCurvEq}
\til{\K}_v(w) = \K_v(w) - \frac{1}{2} \frac{\nabla^2_v \varphi(v)}{\varphi(v)} + \frac{3}{4} \frac{(\nabla_v \varphi (v))^2}{\varphi(v)^2} ,
\end{equation}
%--------------------------------------------------------------------
where $\nabla$ is the anisotropic derivative associated with the Chern connection of $L$.
\end{theorem}
%--------------------------------------------------------------------
%--------------------------------------------------------------------
\begin{proof}
Bear in mind that the null flag curvatures $\K_v(w)$ and $\til{\K}_v(w)$ are only well-defined if $v$ is a lightlike vector for both $L$ and $\til{L}$, and $w$ is a nondegenerate vector such that $g_v(v, w) = 0 = \til{g}_v(v,w)$; recall Definition~\ref{NullFlgaCurvature}. Indeed, it can be seen from Lemma~\ref{ConformalFundTensorLemma}, that $L, \til{L}$ share the same lightlike cone, and $g_v(v,w) = 0$ if and only if $\til{g}_v (v,w) = 0$.

Let $\gamma$ be a geodesic for $L$ such that $\dot\gamma(0) = v$ is lightlike, and let $J$ be the Jacobi field along $\gamma$ satisfying $J(0) = 0$ and $J^{\prime}(0) = w$, with $g_v(v,w) = 0$. From Lemma~\ref{AffJangle}, it follows that $g_{\dot \gamma}(\dot\gamma, J) = 0$. Lemma~\ref{ConformalFundTensorLemma} then implies that
%--------------------------------------------------------------------
\begin{equation}
\label{RelationForExpansion}
\til{g}_{\dot\gamma} (J,J) = \varphi(\dot\gamma) g_{\dot\gamma}(J,J) .
\end{equation}
%--------------------------------------------------------------------

The proof of the theorem will be based on comparing the fourth order terms of the Taylor expansions on both sides of \eqref{RelationForExpansion}. To do this, we present some preliminary calculations in the format of claims.

%--------------------------------------------------------------------
%   CLAIM A
%--------------------------------------------------------------------
\begin{claim}
\label{ClaimA}
The Taylor expansions for $g_{\dot\gamma} (J,J)$ and $\til{g}_{\dot\gamma} (J,J)\circ\sigma$ up to order $4$ around $0$ are
%--------------------------------------------------------------------
\begin{subequations}
\begin{equation*}
g_{\dot\gamma} (J,J)(t) = \underbrace{g_{v} (w,w)}_{\mathrm{a}_2} t^2 \underbrace{- \frac{1}{3} g_{v} (R_{v} (w,v)v,w)}_{\mathrm{a}_4} t^4 + \mathcal{O} (t^5) ,
\end{equation*}
\begin{equation*}
\til{g}_{\dot\gamma} (J,J) (\sigma(s)) = \underbrace{\frac{1}{\varphi(v)^2} \til{g}_v (w,w)}_{\til{\mathrm{a}}_2} s^2 \underbrace{- \frac{1}{3\varphi(v)^4} \til{g}_v (\til{R}_v(w,v)v,w)}_{\til{\mathrm{a}}_4} s^4 + \mathcal O(s^5) .
\end{equation*}
\end{subequations}
%--------------------------------------------------------------------
\end{claim}
%--------------------------------------------------------------------
%--------------------------------------------------------------------
\begin{proof}
The image of $\gamma$ is a geodesic for $\til{L}$ up to parametrization; say $\til{\gamma}$ is such lightlike geodesic for $\til{L}$ with $\dot{\til{\gamma}}(0) = \til{v} = \tfrac{v}{\varphi (v)}$ by \eqref{GammaTil}. Consider the Jacobi field $\til{J}$ along $\til{\gamma}$ given by \eqref{JTilEq}. By Lemma~\ref{JTilLemma}, $\til{J}$ is determined by the initial conditions $\til{J}(0) = 0$ and $\til{J}^{\prime}(0) = \til{w} = \frac{1}{\varphi(v)} w +\frac{\dot{h}(0)}{\varphi(v)} v$. The Taylor expansions around $0$ for $g_{\dot\gamma} (J,J)$ and $\til{g}_{\dot{\til{\gamma}}} (\til{J},\til{J})$ are, as proved in Lemma~\ref{JacobiTaylor},
%--------------------------------------------------------------------
\begin{subequations}
\begin{equation}
\label{ExpansionJacobiNorm}
g_{\dot\gamma} (J,J) (t) = \underbrace{g_{v} (w,w)}_{\mathrm{a}_2} t^2 \underbrace{- \frac{1}{3} g_{v} (R_{v} (w,v)v,w)}_{\mathrm{a}_4} t^4 + \mathcal{O} (t^5) ,
\end{equation}
\begin{equation}
\label{ExpansionJacobiNormTil}
\til{g}_{\dot{\til\gamma}} (\til{J}, \til{J}) (s)
= \til{g}_{\til{v}} (\til{w}, \til{w}) s^2 - \frac{1}{3}\til{g}_{\til{v}} (\til{R}_{\til{v}} (\til{w},\til{v})\til{v} , \til{w}) s^4 + \mathcal O(s^5) .
\end{equation}
\end{subequations}
%--------------------------------------------------------------------
It follows from positive $0$-homogeneity and multilinearity of the fundamental and curvature tensors that the coefficients on the right-hand side of \eqref{ExpansionJacobiNormTil} are
%--------------------------------------------------------------------
\begin{subequations}
\begin{equation*}
\til{g}_{\til{v}}(\til{w},\til{w}) = \frac{1}{\varphi(v)^2} \til{g}_{v}(w,w) = \til{\mathrm{a}}_2 ,
\end{equation*}
\begin{equation*}
- \frac{1}{3} \til{g}_{\til{v}}(\til{R}_{\til{v}}(\til{w},\til{v})\til{v},\til{w}) = - \frac{1}{3\varphi(v)^4} \til{g}_{v}(\til{R}_{v}(w,v)v,w) = \til{\mathrm{a}}_4 .
\end{equation*}
\end{subequations}
%--------------------------------------------------------------------
Moreover, $\til{g}_{\til{v}}(\til{v},\til{w}) = \tfrac{1}{\varphi(v)^2} \til{g}_v(v,w) \overset{\eqref{v-Orthogonalities}}{=} \tfrac{1}{\varphi(v)} g_v(v,w) = 0$, and it follows from Lemma~\ref{AffJangle} again that $\til{g}_{\dot{\til{\gamma}}}(\dot{\til{\gamma}},\til{J}) = 0$. Then, equation~\eqref{JTilEq} implies that $\til{g}_{\dot{\til{\gamma}}}(\dot{\til{\gamma}},J\circ\sigma) = 0$. Combining this with identities \eqref{GammaTil} and \eqref{JTilEq} once again, the left-hand side of \eqref{ExpansionJacobiNormTil} is
%--------------------------------------------------------------------
\begin{equation*}
\til{g}_{\dot{\til{\gamma}}}(\til{J},\til{J}) = \til{g}_{\dot\gamma}(J,J)\circ\sigma ,
\end{equation*}
%--------------------------------------------------------------------
which completes the proof of Claim~\ref{ClaimA}.
\phantom\qedhere
\end{proof}
%--------------------------------------------------------------------

%--------------------------------------------------------------------
%   CLAIM B
%--------------------------------------------------------------------
\begin{claim}
\label{ClaimB}
The Taylor expansion of $\varphi(\dot\gamma)$ around $0$ is
%--------------------------------------------------------------------
\begin{equation*}
\varphi(\dot\gamma(t)) = \underbrace{\varphi(v)}_{\mathrm{b}_0} + \underbrace{\nabla_{v} \varphi(v)}_{\mathrm{b}_1} t + \underbrace{\frac{1}{2} \nabla^2_{v} \varphi(v)}_{\mathrm{b}_2} t^2 + \mathcal O (t^3) .
\end{equation*}
%--------------------------------------------------------------------
\end{claim}
%--------------------------------------------------------------------
%--------------------------------------------------------------------
\begin{proof}
By \eqref{PsiPhi} and straightforward differentiation, we have that
%--------------------------------------------------------------------
\begin{equation}
\label{DerivativesPsi}
\dot\psi (t) = \varphi(\dot\gamma(t)) , \;
\ddot\psi (t) = \nabla_{\dot\gamma (t)} \varphi (\dot\gamma(t)) , \;
\psi^{(3)}(t) = \nabla^2_{\dot\gamma (t)} \varphi (\dot\gamma (t)) ,
\end{equation}
%--------------------------------------------------------------------
where it is essential  that $\gamma$ be a geodesic; recall \eqref{DerivAFunct}. Since $\psi(0)=0$, the Taylor expansion of $\psi$ up to order $3$ around $0$ is
%--------------------------------------------------------------------
\begin{equation*}
\psi(t) = \varphi(v) t + \frac{\nabla_{v} \varphi(v)}{2} t^2 + \frac{\nabla^2_{v} \varphi(v)}{6} t^3 + \mathcal O(t^4) ,
\end{equation*}
%--------------------------------------------------------------------
and the derivative of the above equation proves Claim~\ref{ClaimB}.
\phantom\qedhere
\end{proof}
%--------------------------------------------------------------------

%--------------------------------------------------------------------
%   CLAIM C
%--------------------------------------------------------------------
\begin{claim}
\label{ClaimC}
 The Taylor expansion of $\sigma$ around $0$ is
%--------------------------------------------------------------------
\begin{multline*}
\sigma(s) = \underbrace{\frac{1}{\varphi(v)}}_{\mathrm{c}_1} s \underbrace{-\frac{\nabla_{v} \varphi(v)}{2\varphi(v)^3}}_{\mathrm{c}_2} s^2 \underbrace{- \frac{1}{6\varphi(v)^5} (\varphi(v) \nabla^2_{v} \varphi(v) - 3 (\nabla_{v} \varphi(v))^2)}_{\mathrm{c}_3} s^3 \\
+ \mathcal O(s^4) .
\end{multline*}
%--------------------------------------------------------------------
\end{claim}
%--------------------------------------------------------------------
%--------------------------------------------------------------------
\begin{proof}
Remember that $\sigma$ is the inverse of $\psi$. Then, $\sigma(0) = 0$ and the derivatives of $\sigma$ are
%--------------------------------------------------------------------
\begin{equation*}
\dot\sigma = \frac{1}{\dot\psi(\sigma)} , \;
\ddot\sigma = - \frac{\ddot\psi(\sigma)}{\dot\psi(\sigma)^3} , \;
\sigma^{(3)} = - \frac{1}{\dot\psi(\sigma)^5} \left( \dot\psi (\sigma) \psi^{(3)}(\sigma) - 3 \ddot\psi(\sigma)^2 \right) .
\end{equation*}
%--------------------------------------------------------------------
So Claim~\ref{ClaimC} follows by evaluating these derivatives at $0$, using \eqref{DerivativesPsi}.
\phantom\qedhere
\end{proof}
%--------------------------------------------------------------------

Claim~\ref{ClaimA} directly gives the Taylor expansion for the left-hand side of \eqref{RelationForExpansion}. It only remains to compute the expansion on the right-hand side of \eqref{RelationForExpansion}.

From Claims~\ref{ClaimA} and~\ref{ClaimB}, we obtain that
%--------------------------------------------------------------------
\begin{equation*}
\varphi(\dot\gamma) g_{\dot\gamma}(J,J) = \mathrm{a}_2 \mathrm{b}_0 t^2 + \mathrm{a}_2 \mathrm{b}_1 t^3 + (\mathrm{a}_2 \mathrm{b}_2 + \mathrm{a}_4 \mathrm{b}_0) t^4 + \mathcal O(t^5) .
\end{equation*}
%--------------------------------------------------------------------
From Claim~\ref{ClaimC}, we get that
%--------------------------------------------------------------------
\begin{align*}
\sigma(s)^2 &= \mathrm{c}_1^2 s^2 + 2 \mathrm{c}_1 \mathrm{c}_2 s^3 + (2 \mathrm{c}_1 \mathrm{c}_3 + \mathrm{c}_2^2) s^4 + \mathcal{O}(s^5) , \\
\sigma(s)^3 &= \mathrm{c}_1^3 s^3 + 3 \mathrm{c}_1^2 \mathrm{c}_2 s^4 + \mathcal{O}(s^5) , \\
\sigma(s)^4 &= \mathrm{c}_1^4 s^4 + \mathcal{O}(s^5) .
\end{align*}
%--------------------------------------------------------------------
Thus,
%--------------------------------------------------------------------
\begin{multline}
\label{ExpansionJacobiNormAlongSigma}
(\varphi(\dot\gamma) g_{\dot\gamma}(J,J)) \circ \sigma = \mathrm{a}_2 \mathrm{c}_1 s^2 + \mathrm{a}_2 (2 \mathrm{c}_2 + \mathrm{b}_1 \mathrm{c}_1^3 ) s^3 \\
+ \left( \mathrm{a}_2 ( 2 \mathrm{c}_3 + \mathrm{b}_0 \mathrm{c}_2^2 + 3 \mathrm{b}_1 \mathrm{c}_1^2 \mathrm{c}_2 + \mathrm{b}_2 \mathrm{c}_1^4 ) + \mathrm{a}_4\mathrm{c}_1^3 \right) s^4 + \mathcal O(s^5) ,
\end{multline}
%--------------------------------------------------------------------
where the identity $\mathrm{b}_0 \mathrm{c}_1 = 1$ has been used. The equality $\til{g}_v(w,w) = \varphi(v) g_v(w,w)$, obtained from \eqref{v-Orthogonalities}, indicates that $\mathrm{a}_2 \mathrm{c}_1 = \til{\mathrm{a}}_2$. Besides, it is trivial to verify that $2 \mathrm{c}_2 + \mathrm{b}_1 \mathrm{c}_1^3 = 0$.

Comparison of the $4^{\text{th}}$ order terms in the equality \eqref{RelationForExpansion}, whose expansions are given by Claim~\ref{ClaimA} and \eqref{ExpansionJacobiNormAlongSigma}, implies that
%--------------------------------------------------------------------
\begin{equation*}
\begin{aligned}
\til{\mathrm{a}}_4 &= \mathrm{a}_2 (2 \mathrm{c}_3 + \mathrm{b}_0 \mathrm{c}_2^2 + 3 \mathrm{b}_1 \mathrm{c}_1^2 \mathrm{c}_2 + \mathrm{b}_2 \mathrm{c}_1^4) + \mathrm{a}_4 \mathrm{c}_1^3 \\
&= \mathrm{a}_2 \left( \frac{1}{3}\frac{\mathrm{b}_2}{\mathrm{b}_0^4} - \frac{1}{4}\frac{\mathrm{b}_1^2}{\mathrm{b}_0^5} \right) + \frac{\mathrm{a}_4}{\mathrm{b}_0^3} ,
\end{aligned}
\end{equation*}
%--------------------------------------------------------------------
where we have used that $\mathrm{c}_1 = \frac{1}{\mathrm{b}_0}$, $\mathrm{c}_2 = -\frac{1}{2}\frac{\mathrm{b}_1}{\mathrm{b}_0^3}$, and $\mathrm{c}_3 = \frac{1}{2}\frac{\mathrm{b}_1^2}{\mathrm{b}_0^5} -\frac{1}{3}\frac{\mathrm{b}_2}{\mathrm{b}_0^4}$. Finally, replacing $\til{\mathrm{a}}_4$, $\mathrm{a}_2$, $\mathrm{a}_4$, $\mathrm{b}_0$, $\mathrm{b}_1$ and $\mathrm{b}_2$ by their respective values, this simplifies to
%--------------------------------------------------------------------
\begin{multline}
\label{ProofEq}
\frac{1}{\varphi(v)} \til{g}_v (\til{R}_v(w,v)v,w) = \\
= g_v(R_v(w,v)v,w) - g_v(w,w) \left( \frac{1}{2} \frac{\nabla^2_v \varphi (v)}{\varphi(v)} - \frac{3}{4} \frac{(\nabla_v \varphi (v))^2}{\varphi(v)^2} \right) .
\end{multline}
%--------------------------------------------------------------------
Observe that \eqref{ProofEq} holds for every $w \in T_v\C$, degenerate or not, since the argument above only uses that $g_v(v,w)=0$. Using the definition of the null flag curvature, and once again the equality $\til{g}_v(w,w) = \varphi(v) g_v(w,w)$, the above equation yields \eqref{ConformalFlagCurvEq}.
\end{proof}
%--------------------------------------------------------------------

Notice the null flag curvatures of anisotropically conformal pseudo-Finsler metrics differ by a quantity that depends only on the lightlike flagpole; namely,
%--------------------------------------------------------------------
\begin{equation*}
\til{\K}_{v}(w) - \K_{v}(w) = - \frac{1}{2} \frac{\nabla^2_v \varphi(v)}{\varphi(v)} + \frac{3}{4} \frac{(\nabla_v \varphi(v))^2}{\varphi(v)^2} .
\end{equation*}
%--------------------------------------------------------------------
In addition, we can extend the result to the Jacobi operator:
%--------------------------------------------------------------------
\begin{equation}
\label{JacobiOperator}
\J_v(w) = R_v (w,v)v , \; \forall v \in A , w \in T_{\pi(v)}M.
\end{equation}
%--------------------------------------------------------------------

%--------------------------------------------------------------------
\begin{corollary}
Let $M$ be a smooth manifold with $n = \dim M \geq 3$, and let $L,\til{L}$ be two anisotropically conformal pseudo-Finsler metrics on $M$ with $\til{L} = \varphi L$. Let $\J_v(w), \til{\J}_v(w)$ denote the Jacobi operators for $L,\til{L}$, respectively. If $v$ is lightlike (for either) and $w$ satisfies $g_v(v,w) = 0$ (or equivalently, $\til{g}_v(v,w) = 0$), then there exists $\mu(v,w) \in \R$ such that
%--------------------------------------------------------------------
\begin{equation}
\label{ConformalJacobiOp}
\til{\mathcal{R}}_v(w) - \mathcal{R}_v(w) = \left( - \tfrac{1}{2} \tfrac{\nabla^2_v \varphi (v)}{\varphi(v)}  + \tfrac{3}{4} \tfrac{(\nabla_v \varphi (v))^2}{\varphi(v)^2} \right) w + \mu(v,w) v .
\end{equation}
%--------------------------------------------------------------------
\end{corollary}
%--------------------------------------------------------------------
%--------------------------------------------------------------------
\begin{proof}
By~\eqref{symmetry2} and~\eqref{symmetry3}, $\til{g}_v(\til{R}_v(w,v)v,u)$ and $g_v(R_v(w,v)v,u)$ are symmetric in $w$ and $u$. When $v \in \C$, we get from~\eqref{ConformalFundTensorEq} and~\eqref{symmetry1} that $\varphi(v) g_v(\til{R}_v(w,v)v,v) = \til{g}_v(\til{R}_v(w,v)v,v) = 0$, because $(\vd\varphi)_v(v) = 0$ by positive-homogeneity. Since $\varphi(v) \neq 0$, we obtain $g_v(\til{R}_v(w,v)v,v) = 0$. From the last identity and \eqref{v-Orthogonalities}, if $u \in T_v\C$ (recall \eqref{ConeTangentSp}), then
%--------------------------------------------------------------------
\begin{equation*}
\til{g}_v (\til{R}_v(w,v)v,u) = \varphi(v)g_v(\til{R}_v(w,v)v,u) .
\end{equation*}
%--------------------------------------------------------------------
So $g_v(\til{R}_v(w,v)v,u)$ is also symmetric in $w,u \in T_v\C$.

Now, it follows from \eqref{ProofEq} in the proof of Theorem~\ref{AnisotropicConformalFlagCurv} that
%--------------------------------------------------------------------
\begin{equation*}
g_v \left(\til{R}_v(w,v)v - R_v(w,v)v + \left( \tfrac{(\nabla_v^2 \varphi(v))}{2\varphi(v)} - \tfrac{3(\nabla_v \varphi(v))^2}{4\varphi(v)^2} \right) w , w \right) = 0 ,
\end{equation*}
%--------------------------------------------------------------------
for any $w \in T_v\C$, since \eqref{ProofEq} holds for all such $w$; see the observation after \eqref{ProofEq}. By the polarization identity, we obtain
%--------------------------------------------------------------------
\begin{equation*}
g_v \left(\til{R}_v(w,v)v - R_v(w,v)v +  \left( \tfrac{(\nabla_v^2 \varphi(v))}{2\varphi(v)} - \tfrac{3(\nabla_v \varphi(v))^2}{4\varphi(v)^2} \right) w , u \right) = 0 ,
\end{equation*}
%--------------------------------------------------------------------
for all $u \in T_v\C$. This means the right-hand side vector of the above product is $g_v$-orthogonal to all the elements of $T_v\C$ (including $v$), so it must belong to the radical of $(T_v\C,g_v)$, which is $\spn{v}$; recall~\eqref{ConeOrt}. Hence \eqref{ConformalJacobiOp}.
\end{proof}
%--------------------------------------------------------------------
%--------------------------------------------------------------------
\section{Null Weyl curvature}
\label{weylsection}
%--------------------------------------------------------------------

The Ricci curvature is defined as the contraction of the Jacobi operator \eqref{JacobiOperator}; i.e. $\rc(v) := \tr(\J_v)$, $\forall v \in A$. We will see that the traceless part of $\J_v$ is a conformal invariant when $v$ is lightlike and $w$ is tangent to the lightcone.

%--------------------------------------------------------------------
\begin{definition}
\label{weyldef}
Assume, throughout the section, that $(M,L)$ is a pseudo-Finsler manifold with $n = \dim M \geq 3$. We define the {\em null Weyl curvature} by
%--------------------------------------------------------------------
\begin{equation}
\W_v(w) := \J_v(w) - \tfrac{1}{n-2}\rc(v)w ,
\end{equation}
%--------------------------------------------------------------------
for any $v \in \C$, and $w \in T_{\pi(v)}M$. Compare this notion with that of Weyl shadow flag curvature of a causal structure in \cite[\S~2.3.9]{mak18d}.
\end{definition}
%--------------------------------------------------------------------

First observe that, with use of symmetry~\eqref{symmetry1}, $\W_v(w)\in T_v\C$ for all $v$ lightlike and $w\in T_v\C$. Given an element $w\in T_v\C$, we will denote by $[w]$ the class of $w$ in $\bigslant{T_v\C}{\spn{v}}$; recall~\eqref{ConeOrt}. Since $[\W_{v}(w+\lambda v)] = [\W_v(w)]$, the null Weyl curvature can pass to the quotient as a map
%--------------------------------------------------------------------
\begin{equation}
\label{Weylmap}
[\W_v]: \bigslant{T_v\C}{\spn{v}} \to \bigslant{T_v\C}{\spn{v}} .
\end{equation}
%--------------------------------------------------------------------
For the same reasons, for any $v \in \C$, we can also consider the map
%--------------------------------------------------------------------
\begin{equation}
\label{Ricmap}
[\J_v]: \bigslant{T_v\C}{\spn{v}} \to \bigslant{T_v\C}{\spn{v}} .
\end{equation}
%--------------------------------------------------------------------
As $\spn{v}$ is the radical of $(T_v\C,g_v)$, it turns out that $g_v$ restricts to $\bigslant{T_v\C}{\spn{v}}$ as a well-defined and nondegenerate scalar product.

%--------------------------------------------------------------------
\begin{lemma}
\label{ricquotient}
Let $(M,L)$ be a pseudo-Finsler manifold. If $v$ is lightlike, then $\rc(v)$ can be computed as the trace of the map \eqref{Ricmap} using the nondegenerate scalar product $g_v$ on $\bigslant{T_v\C}{\spn{v}}$.
\end{lemma}
%--------------------------------------------------------------------
%--------------------------------------------------------------------
\begin{proof}
Since $g_v$ is nondegenerate on $\bigslant{T_v\C}{\spn{v}}$, it is possible to consider an orthonormal basis $\{[\e_{i}]\}_{i=1}^{n-2}$ of $\left( \bigslant{T_v\C}{\spn{v}},g_v \right)$. In this case, $\{\e_{i}\}_{i=1}^{n-2}$ is a $g_v$-orthonormal system and it can be extended to a $g_v$-orthonormal basis $\{\e_{i}\}_{i=1, \ldots, n}$ of $(T_{\pi(v)}M,g_v)$ in such a way that $\e_{n-1} = v+u$, and $\e_n = v-u$, where $u$ is a $g_v$-lightlike vector with $g_v(v,u) = \tfrac{1}{2}$. If $\varepsilon_i = g_v(\e_i,\e_i)$, then
%--------------------------------------------------------------------
\begin{align*}
\rc(v) =& \tr(\J_v) \\
=& \sum_{i=1}^{n-2} \varepsilon_i g_{v}(R_{v}(\e_{i}, v)v, \e_{i}) \\
&+ \varepsilon_{n-1} g_{v}(R_{v}(u+v, v)v, u+v) + \varepsilon_{n} g_{v}(R_{v}(v-u, v)v, v-u) .
\end{align*}
%--------------------------------------------------------------------
The first summand on the right is just the trace of \eqref{Ricmap}. By the multilinearity of $R_v(X,Y)$, its antisymmetry in $X,Y$, and \eqref{symmetry1}, we have that
%--------------------------------------------------------------------
\begin{multline*}
\varepsilon_{n-1} g_{v}(R_{v}(u+v, v)v, u+v)+\varepsilon_{n} g_{v}(R_{v}(v-u, v)v, v-u) = \\
= (\varepsilon_{n-1}+\varepsilon_n)g_{v}(R_{v}(u, v)v, u) .
\end{multline*}
%--------------------------------------------------------------------
Since $\varepsilon_{n-1} + \varepsilon_n = 1-1 = 0$, we conclude the desired result.
\end{proof}
%--------------------------------------------------------------------

%--------------------------------------------------------------------
\begin{corollary}
\label{weylinv}
The null Weyl curvature map \eqref{Weylmap} is invariant under anisotropic conformal changes.
\end{corollary}
%--------------------------------------------------------------------
%--------------------------------------------------------------------
\begin{proof}
Let $L,\til{L}$ be anisotropically conformal pseudo-Finsler metrics. Taking the trace of \eqref{ConformalJacobiOp}, we get
%--------------------------------------------------------------------
\begin{equation}
\label{TraceRic}
\til{\rc}(v) - \rc(v) = (n-2) \left( - \tfrac{1}{2} \tfrac{\nabla^2_v \varphi (v)}{\varphi(v)}  + \tfrac{3}{4} \tfrac{(\nabla_v \varphi (v))^2}{\varphi(v)^2} \right).
\end{equation}
%--------------------------------------------------------------------
Then, from \eqref{ConformalJacobiOp}, \eqref{TraceRic}, and the definition of the Weyl curvature as the map \eqref{Weylmap}, we conclude that $[\til{\W}_v(w)]=[\W_v(w)]$, as required.
\end{proof}
%--------------------------------------------------------------------

%--------------------------------------------------------------------
\begin{corollary}
Let $(M,L)$ and $(M,\til{L})$ be two anisotropically conformal pseudo-Finsler spaces. If $v\in\C$ and $w \in T_{\pi(v)}M$ satisfies $g_v(v,w)=0$ with $g_v(w,w)\not=0$, then
%--------------------------------------------------------------------
\begin{equation}
\label{Ricci_comparison}
\tfrac{1}{n-2} \left( \til{\rc}(v) - \rc(v) \right) = \til{\K}_{v}(w) - \K_{v}(w) .
\end{equation}
%--------------------------------------------------------------------
\end{corollary}
%--------------------------------------------------------------------
%--------------------------------------------------------------------
\begin{proof}
It  follows directly from \eqref{TraceRic} and \eqref{ConformalFlagCurvEq}.
\end{proof}
%--------------------------------------------------------------------

%--------------------------------------------------------------------
\begin{proposition}
\label{nullW}
The null Weyl curvature map \eqref{Weylmap} vanishes for $v \in \C$ if and only if the null flag curvature is scalar; i.e. there exists some scalar $\lambda$ depending on the lightlike vector $v$ such that $\K_v(w) = \lambda(v)$, for all nondegenerate $w\in T_v\C$, which is equivalent to $[\J_v(w)] = \lambda(v) [w]$.
\end{proposition}
%--------------------------------------------------------------------
%--------------------------------------------------------------------
\begin{proof}
The forward implication is immediate from the definitions. Assume now that the null flag curvature is scalar with $[\J_v(w)] = \lambda(v) [w]$. Then, taking the trace of the last identity in $(\bigslant{T_v\C}{\spn{v}},g_v)$, one gets that $\rc(v)=(n-2)\lambda(v)$; recall Lemma~\ref{ricquotient}. As a consequence,
%--------------------------------------------------------------------
\begin{equation*}
[\W_v(w)] = [\J_v(w)] -\tfrac{1}{n-2}\rc(v)[w] = \lambda(v)[w] -\tfrac{1}{n-2}(n-2)\lambda(v)[w] = 0 .
\end{equation*}
%--------------------------------------------------------------------
\end{proof}
%--------------------------------------------------------------------

%--------------------------------------------------------------------
\begin{remark}
Observe that the Weyl tensor in a semi-Riemannian manifold $(M,g)$ is defined as
%--------------------------------------------------------------------
\begin{align*}
\Cw(X, Y, Z, W) =& g(R(X, Y) W, Z) \\
&-\tfrac{1}{n-2} \Big( \rc(X,Z)g(Y,W) - \rc(X,W)g(Y,Z) \\
&\qquad + g(X,Z)\rc(Y,W) - g(X,W)\rc(Y,Z) \Big) \\
&+ \tfrac{S}{(n-1)(n-2)} \Big( g(X,Z)g(Y,W) - g(X,W)g(Y,Z) \Big) ,
\end{align*}
%--------------------------------------------------------------------
where $R$, $\rc$, and $S$ are, respectively, the curvature tensor, the Ricci tensor, and the scalar curvature of $(M,g)$; see \cite[(1.116)]{bes87e}, which uses the opposite sign convention for the curvature. In particular, if $Y$ is lightlike and $g(Y, X) = g(Y, Z) = 0$, then we have
%--------------------------------------------------------------------
\begin{align*}
\Cw(X, Y, Z, Y) &= g(R(X,Y)Y,Z) - \tfrac{1}{n-2} g(X,Z) \rc(Y,Y) \\
&= g(\W_Y(X),Z) . 
\end{align*}
%--------------------------------------------------------------------
Therefore, our null Weyl curvature is just a certain projection of the classical Weyl tensor.
\end{remark}
%--------------------------------------------------------------------
%--------------------------------------------------------------------
\section{Null flag curvature of the associated Lorentz--Finsler metric}
\label{associatespacetime}
%--------------------------------------------------------------------

Reconsider the Zermelo problem studied in Section~\ref{Zermelo}, but for a Finsler metric $F$ that does not depend on time. Now, the Lorentz--Finsler metric in \eqref{LFmetric} takes the form
%--------------------------------------------------------------------
\begin{equation}
\label{LFmetric2}
L(\tau, v) = \tau^2 - F(v)^2 ,
\end{equation}
%--------------------------------------------------------------------
which is a Lorentz--Finsler metric as in \S \ref{LorentzFinsler} up to the smoothness of $L$ in the direction $(\tau,0)$ for $\tau>0$. Nonetheless, \eqref{LFmetric2} determines the open subset $A=\{(t,\tau,v)\in \R\times \R\times TM: L(\tau,v)>0, \tau>0 \}$ such that $L$ extends smoothly to its boundary, where $L=0$ and the fundamental tensor is nondegenerate. As we will compute the curvature in the lightlike directions of the boundary, the possible lack of smoothness on $(t,\tau,0)$ is not a problem. Moreover, the boundary $\C$ is a cone structure as defined in \S \ref{LorentzFinsler} and it is possible to modify $L$ in a small neighborhood of the half-ray $\tau\mapsto (t,\tau,0)$ to obtain a proper Lorentz--Finsler metric; see \cite[Theorem~5.6]{jav20d}.

In the following, as the pseudo-Finsler metric $L$ in \eqref{LFmetric2} does not depend on the first coordinate $t$, we will consider this metric defined on $\R\times TM$, disregarding $t$.

%--------------------------------------------------------------------
\begin{lemma}
\label{connections}
Let $\hnab$ and $\nabla$ be respectively the Chern connections of $L$ and $F$, where $L$ satisfies \eqref{LFmetric2}. If $c_1,c_2\in\R$, $X,Y\in \mathfrak{X}(M)$, $V$ is any vector field of some open subset $\Omega\subset M$ without zeroes, and $f$ is any positive function on $\Omega$, then
%--------------------------------------------------------------------
\begin{equation*}
\hnab^{(f,V)}_{(c_1,X)}(c_2,Y)=(0,\nabla^V_XY) .
\end{equation*}
%--------------------------------------------------------------------
\end{lemma}
%--------------------------------------------------------------------
%--------------------------------------------------------------------
\begin{proof}
First observe that by the positive $0$-homogeneity of $\hnab^{(f,V)}$ we can assume that $f=1$. Next, let us point out that
%--------------------------------------------------------------------
\begin{subequations}
\begin{equation*}
[(c_1,X),(c_2,Y)]=(0,[X,Y]) ;
\end{equation*}
\begin{equation*}
g^L_{(1,V)}((c_1,X),(c_2,Y))=c_1c_2-g^F_V(X,Y) ;
\end{equation*}
\begin{equation*}
C^L_{(1,V)}((c_1,X),(c_2,Y),(c_3,Z)))=-C^F_V(X,Y,Z) .
\end{equation*}
\end{subequations}
%--------------------------------------------------------------------
In the first six terms of the Koszul formula~\eqref{KoszulQ}, we can replace $g^L$ with $-g^F$, just removing the first component everywhere. So we only have to deal with the last three terms of \eqref{KoszulQ}.

When we consider $\hnab^{(1,V)}_{(1,V)}(1,V)$, the last three terms of \eqref{KoszulQ} vanish, and we conclude that
%--------------------------------------------------------------------
\begin{equation*}
g^L_{(1,V)}(\hnab^{(1,V)}_{(1,V)}(1,V),(c_1,X)) = -g^F_V(\nabla^V_V V,X) .
\end{equation*}
%--------------------------------------------------------------------
In particular, the last formula for $X=0$ implies that the first component of $\hnab^{(1,V)}_{(1,V)}(1,V)$ is zero, and for $c_1=0$, it implies that
the second component is $\nabla^V_VV$.  Summing up, we have
%--------------------------------------------------------------------
\begin{equation*}
\hnab^{(1,V)}_{(1,V)}(1,V) = (0,\nabla^V_V V) .
\end{equation*}
%--------------------------------------------------------------------
Using the above identity to handle the last three terms in the Koszul formula, we conclude that 
%--------------------------------------------------------------------
\begin{equation*}
g^L_{(1,V)}(\hnab^{(1,V)}_{(c_1,X)}(1,V),(c_2,Y))=-g^F_V(\nabla^V_X V,Y) .
\end{equation*}
%--------------------------------------------------------------------
By similar reasoning, it follows that
%--------------------------------------------------------------------
\begin{equation*}
\hnab^{(1,V)}_{(c_1,X)}(1,V)=(0,\nabla^V_X V) ,
\end{equation*}
%--------------------------------------------------------------------
and the use of this identity finally allows us to obtain that
%--------------------------------------------------------------------
\begin{equation*}
g^L_{(1,V)}(\hnab^{(1,V)}_{(c_1,X)}(c_2,Y),(c_3,Z))=-g^F_V(\nabla^V_XY,Z) .
\end{equation*}
%--------------------------------------------------------------------
Repeating the previous reasoning leads to
%--------------------------------------------------------------------
\begin{equation*}
\hnab^{(1,V)}_{(c_1,X)}(c_2,Y)=(0,\nabla^V_X Y) ,
\end{equation*}
%--------------------------------------------------------------------
as required.
\end{proof}
%--------------------------------------------------------------------

%--------------------------------------------------------------------
\begin{lemma}
\label{Cherntensors}
The Chern tensors of $\hnab$ and $\nabla$, denoted by $\hat{P}$ and $P$, are related by
%--------------------------------------------------------------------
\begin{equation*}
\hat P_{(1,V)}((c_1,X),(c_2,Y),(c_3,Z)) = (0,P_V(X,Y,Z)) .
\end{equation*}
%--------------------------------------------------------------------
\end{lemma}
%--------------------------------------------------------------------
%--------------------------------------------------------------------
\begin{proof}
It follows directly from Lemma~\ref{connections} and the definitions of $\hat P$, $P$; recall \eqref{ChernTensor}. Indeed,
%--------------------------------------------------------------------
\begin{align*}
\hat P_{(1,V)}((c_1,X),(c_2,Y),(c_3,Z)) &= \left.\frac{\d}{\d t}  \Big( \hnab^{(1,V)+t(c_3,Z)}_{(c_1,X)} (c_2,Y) \Big) \right\vert_{t=0} \\
&= \left( 0 , \left. \frac{\d}{\d t} \big( \nabla^{V+tZ}_X Y \big) \right\vert_{t=0} \right) \\
&= (0,P_V(X,Y,Z)) .
\end{align*}
%--------------------------------------------------------------------
\end{proof}
%--------------------------------------------------------------------

%--------------------------------------------------------------------
\begin{proposition}
\label{curvatures}
Let $\hat R$ and $R$ be the curvature tensors of $L$ and $F$, respectively, where $L$ satisfies \eqref{LFmetric2}. Then
%--------------------------------------------------------------------
\begin{equation*}
\hat R_{(1,V)}((c_1,X),(c_2,Y))((c_3,Z)) = (0,R_V(X,Y)Z) .
\end{equation*}
%--------------------------------------------------------------------
In particular,
%--------------------------------------------------------------------
\begin{equation*}
\hat R_{(1,V)}((0,W),(1,V))(1,V) = (0,R_V(W,V)V) ;
\end{equation*}
%--------------------------------------------------------------------
additionally, if $(1,v)$ is lightlike and $w$ is $g^F_v$-orthogonal to $v$, then the null flag curvature $\K_{(1,v)}(0,w)$ of $L$ is equal to the flag curvature $K_v(w)$ of $F$.
\end{proposition}
%--------------------------------------------------------------------
%--------------------------------------------------------------------
\begin{proof}
It is a straightforward consequence of Lemmas \ref{connections}, \ref{Cherntensors}, and the fact that $g^F_V(V,V)=1$, because $(1,V)$ is lightlike. 
\end{proof}
%--------------------------------------------------------------------

%--------------------------------------------------------------------
\begin{corollary}
\label{KsobreKcor}
Let $\K$ be the null flag curvature of $L$ and $K$ the flag curvature of $F$. For any $v\in TM \setminus \mathbf{0}$, $w\in T_{\pi(v)}M$ not collinear to $v$, and $\tau = \tfrac{1}{F(v)} g_v(v,w)$, the degenerate plane $\spn{(F(v),v),(\tau,w)}$ is a null flag for $L$ and it holds
%--------------------------------------------------------------------
\begin{equation}
\label{KsobreK}
\K_{(F(v),v)}(\tau,w)=F(v)^2 K_v(w).
\end{equation}
%--------------------------------------------------------------------
\end{corollary}
%--------------------------------------------------------------------
%--------------------------------------------------------------------
\begin{proof}
Observe that, by the definition of $\tau$, $(\tau,w)$ is $g^L_{(F(v),v)}$-orthogonal to $(F(v),v)$, and we can apply the definition of null flag curvature.

Using the invariance of the null flag curvature \eqref{NullFlgaCurvWellDef} when adding a multiple of the flagpole to the flag vector $(\tau,w)$, namely
%--------------------------------------------------------------------
\begin{equation*}
(\tau,w)-\tfrac{\tau}{F(v)}(F(v),v) = \left( 0, w - \tau\tfrac{v}{F(v)} \right),
\end{equation*}
%--------------------------------------------------------------------
we obtain that
%--------------------------------------------------------------------
\begin{equation}
\begin{aligned}
\label{Kdeg}
\K_{(F(v),v)}(\tau,w) &= \K_{(F(v),v)}\left( 0, w - \tau\tfrac{v}{F(v)} \right) \\
&= F(v)^2 \K_{\left( 1, \tfrac{v}{F(v)} \right)}\left( 0, w - \tau\tfrac{v}{F(v)} \right) ,
\end{aligned}
\end{equation}
%--------------------------------------------------------------------
where we  have also used that the null flag curvature is positive homogeneous of degree $2$ with respect to the flagpole. From Proposition~\ref{curvatures}, it now follows that
%--------------------------------------------------------------------
\begin{equation}
\label{Kdegflag}
\K_{\left( 1, \tfrac{v}{F(v)} \right)}\left( 0, w - \tau\tfrac{v}{F(v)} \right) = K_{\tfrac{v}{F(v)}} \left( w - \tau\tfrac{v}{F(v)} \right),
\end{equation}
%--------------------------------------------------------------------
and by the properties of the classical flag curvature,
%--------------------------------------------------------------------
\begin{equation}
\label{Kflag}
K_{\tfrac{v}{F(v)}} \left( w - \tau\tfrac{v}{F(v)} \right) = K_v(w).
\end{equation}
%--------------------------------------------------------------------
Putting together \eqref{Kdeg}, \eqref{Kdegflag}, and \eqref{Kflag}, we deduce \eqref{KsobreK}.
\end{proof}
%--------------------------------------------------------------------

%--------------------------------------------------------------------
\begin{corollary}
A Finsler manifold $(M,F)$ has scalar curvature if and only if its associated Finsler spacetime \eqref{LFmetric2} has null Weyl curvature equal to zero in all its lightlike directions.
\end{corollary}
%--------------------------------------------------------------------
%--------------------------------------------------------------------
\begin{proof}
By Corollary~\ref{KsobreKcor}, the flag curvature of $F$ does not depend on $w$ if and only if the null flag curvature of $L$ does not depend on the flag vector $(\tau,w)$. But this is equivalent to the vanishing of the null Weyl curvature map \eqref{Weylmap}, by Proposition~\ref{nullW}.
\end{proof}
%--------------------------------------------------------------------
%--------------------------------------------------------------------
\section{Finsler metrics with the same associated cone structure}
\label{samecone}
%--------------------------------------------------------------------

If $(M,F)$ is a Finsler manifold, and $f:M\rightarrow \R$ is a smooth function such that
%--------------------------------------------------------------------
\begin{equation}
\label{F+df>0}
F(v)+\d f(v)>0 , \; \forall v \in TM\setminus \mathbf{0} ,
\end{equation}
%--------------------------------------------------------------------
then $F+\d f$ is also a Finsler metric; e.g. \cite[Corollary~4.17]{jav14d}. Moreover, $F$ and $F+\d f$ have the same geodesics up to parametrization. Indeed,
%--------------------------------------------------------------------
\begin{equation*}
\mathscr{L}_{F+\d f}(\gamma) = \int_a^b (F(\dot\gamma) + \d f(\dot\gamma)) \d s = \mathscr{L}_F(\gamma) + f(\gamma(b)) - f(\gamma(a)) ,
\end{equation*}
%--------------------------------------------------------------------
i.e. the lengths of a curve between points $p$ and $q$ with respect to $F$ and $F+\d f$ differ by a constant. This means that $F$ and $F+\d f$ are projectively related. By classical results for sprays, the difference of the curvatures of projectively related sprays can be explicitly computed; see \cite[Lemma~12.1.2]{she01d}. In the case of $F$ and $F+\d f$, we will further show that a weighted difference of their flag curvatures depends only on $v$, namely,
%--------------------------------------------------------------------
\begin{equation*}
(F(v)+\d f(v))^2 K^{F+\d f}_v(w) - F(v)^2 K^F_v(w) = \kappa(v) .
\end{equation*}
%--------------------------------------------------------------------
To prove it and compute $\kappa$, we use the null flag curvature of an anisotropically conformal Lorentz--Finsler metric on the associated Finsler spacetime $(\R\times M,L)$, where $L$ is defined by \eqref{LFmetric2}.

First, consider a different splitting of $\R\times M$ using the smooth function $f:M\rightarrow \R$ and the hypersurface $M_f=\{(f(x),x):x\in M\}$. This new splitting can be obtained using the map
%--------------------------------------------------------------------
\begin{equation*}
\phi: \R\times M \to \R\times M
\end{equation*}
%--------------------------------------------------------------------
defined  as $\phi(t,x) = (t - f(x),x)$. In fact, we can induce a new cone structure on $\R\times M$ using the map $\phi$ and the original cone structure $\C$ on $(\R\times M,L)$. Recall that the cone structure associated with $L$ in $\R\times M$ is
%--------------------------------------------------------------------
\begin{equation*}
\C = \{(\tau,v) \in \R\times TM : F(v) = \tau\} .
\end{equation*}
%--------------------------------------------------------------------
Next, observe that the differential of $\phi$ is $\d\phi(\tau,v)=(\tau - \d f(v),v)$. So the new cone structure $\C_f = \phi^{\ast}(\C)$ on $\R\times M$ is given by
%--------------------------------------------------------------------
\begin{equation*}
\C_f = \{ (\tau,v) : F(v) = \tau - \d f(v) \} ,
\end{equation*}
%--------------------------------------------------------------------
which can also be obtained from the Lorentz--Finsler metric
%--------------------------------------------------------------------
\begin{equation}
\label{L_f}
L_f(\tau,v) = \tau^2 - (F(v) + \d f(v))^2 .
\end{equation}
%--------------------------------------------------------------------
Summing up, the map
%--------------------------------------------------------------------
\begin{align*}
\phi: (\R\times M,\C_f) &\to (\R\times M,\C) \\
(t,x) &\mapsto (t - f(x),x)
\end{align*}
%--------------------------------------------------------------------
preserves the cone structures, i.e. $\d\phi(\C_f) = \C$. The cone structure $\C_f = \phi^*(\C)$ is the lightcone of the Lorentz--Finsler metric~\eqref{L_f}, and we can push forward this metric using $\phi$ to define a new one by
%--------------------------------------------------------------------
\begin{align*}
\til{L} (\tau,v) &= L_f (\d\phi^{-1}(\tau,v)) = L_f(\tau + \d f(v),v)\\
&= (\tau + \d f(v))^2 - (F(v) + \d f(v))^2,
\end{align*}
%--------------------------------------------------------------------
which has the same lightcone as $L(\tau,v)=\tau^2-F(v)^2$. The map $\phi: (\R\times M, L_f)\rightarrow (\R\times M, \til{L})$ is now an isometry of Finsler spacetimes. Therefore, the null flag curvature of $L_f$ on the null flag $\Pi = \spn{(F(v) + \d f(v),v),(0,w)}$ coincides with the null flag curvature of $\til{L}$ on
%--------------------------------------------------------------------
\begin{equation*}
\d \phi(\Pi) = \spn{ (F(v),v) , (-\d f(w),w) } ;
\end{equation*}
%--------------------------------------------------------------------
namely,
%--------------------------------------------------------------------
\begin{equation}
\label{dphiPi}
\K^{f}_{(F(v)+\d f(v),v)}(0,w) = \til{\K}_{(F(v),v)}(-\d f(w),w) .
\end{equation}
%--------------------------------------------------------------------
As $\til{L}$ and $L$ have the same lightcone, they are anisotropically conformal, say $\til{L} = \varphi L$, and we can relate both null flag curvatures by Theorem~\ref{AnisotropicConformalFlagCurv}:
%--------------------------------------------------------------------
\begin{equation}
\begin{aligned}
\label{firstRelation}
\til{\K}_{(F(v),v)}(-\d f(w),w) =& \K_{(F(v),v)}(-\d f(w),w)\\
& -\frac{1}{2} \frac{(\hnab^2_{(F(v),v)} \varphi) (F(v),v)}{\varphi(F(v),v)}\\
& +\frac{3}{4} \frac{(\hnab_{(F(v),v)} \varphi)^2(F(v),v)}{\varphi(F(v),v)^2} ,
\end{aligned}
\end{equation}
%--------------------------------------------------------------------
where $\hnab$ is the anisotropic derivative associated with the Chern connection of $L$. Corollary~\ref{KsobreKcor} implies that the null flag curvatures of $L$ and $L_f$ correspond to the flag curvatures of $F$ and $F+\d f$, respectively. As the null flag curvature of $L_f$ equals that of $\til{L}$ (they are isometric), we only need to compute $\varphi$ to relate the flag curvatures of $F$ and $F+\d f$.

%--------------------------------------------------------------------
\begin{lemma}
\label{varphiLema}
Let $L,\til{L}: \mathcal{A} \subset T(\R\times M) \setminus \mathbf{0} \rightarrow \R$ be two Lorentz--Finsler metrics with the same lightcone $\C$, which are then related by $\til{L} = \varphi L$. The function $\varphi$ can be computed at $(\tau,v) \in \C$ as
%--------------------------------------------------------------------
\begin{equation}
\label{varphi}
\varphi(\tau,v)=\frac{\xi(\varphi L)(\tau,v)}{\xi(L)(\tau,v)},
\end{equation}
%--------------------------------------------------------------------
for any $\xi\in T_{\pi(\tau,v)}(\R\times M)\setminus \mathbf{0}$ such that $\xi(L)(\tau,v) \neq 0$.
\end{lemma}
%--------------------------------------------------------------------
%--------------------------------------------------------------------
\begin{proof}
If $(\tau,v)$ lightlike and $\xi \in T_{\pi(\tau,v)}(\R\times M) \setminus \mathbf{0}$, then
%--------------------------------------------------------------------
\begin{equation*}
\xi(\varphi L)(\tau,v) = \xi(\varphi)(\tau,v) L(\tau,v)+\varphi(\tau,v) \xi(L)(\tau,v) = \varphi(\tau,v) \xi(L)(\tau,v) ,
\end{equation*}
%--------------------------------------------------------------------
where the function $L$ is considered only on $T_{\pi(\tau,v)}(\R\times M)$. Hence, we consider only vertical derivatives, and \eqref{varphi} follows straightforwardly.
\end{proof}
%--------------------------------------------------------------------

%--------------------------------------------------------------------
\begin{proposition}\label{F+dfProp}
Let $(M,F)$ be a Finsler manifold and $f:M\rightarrow\R$ a smooth function satisfying \eqref{F+df>0}. Then
%--------------------------------------------------------------------
\begin{equation}
\label{KFdf}
\begin{aligned}
(F(v)+\d f(v))^2 K^{F+\d f}_v(w) = F(v)^2 K^F_v(w) &- \frac{1}{2} \frac{(\nabla^2_{v} \d f) (v)}{(F(v)+\d f(v))} \\
& +\frac{3}{4} \frac{(\nabla_{v} \d f)^2(v)}{(F(v)+\d f(v))^2} ,
\end{aligned}
\end{equation}
where $\nabla$ is the anisotropic derivative associated with the Chern connection of $F$.
%--------------------------------------------------------------------
\end{proposition}
%--------------------------------------------------------------------
%--------------------------------------------------------------------
\begin{proof}
Let us compute $\varphi$ at any point $(\tau,v) \in \C$ using Lemma~\ref{varphiLema} with $\xi=(1,0)$. It turns out that
%--------------------------------------------------------------------
\begin{equation*}
(1,0)(L)(\tau,v)=2\tau=2F(v) ,
\end{equation*}
%--------------------------------------------------------------------
because $\tau = F(v)$ (once $(\tau,v)$ is lightlike). Also,
%--------------------------------------------------------------------
\begin{equation*}
(1,0)(\til{L})(\tau,v) = 2(\tau + \d f(v)) = 2(F(v) + \d f(v)) .
\end{equation*}
%--------------------------------------------------------------------
From~\eqref{varphi}, we obtain that
%--------------------------------------------------------------------
\begin{equation*}
\varphi(\tau,v) = \frac{(1,0)(\til{L})(\tau,v)}{(1,0)(L)(\tau,v)} = \frac{2(F(v) + \d f(v))}{2F(v)} = 1 + \d f\left( \frac{v}{F(v)} \right).
\end{equation*}
%--------------------------------------------------------------------
 
Keeping in mind we denote by $\hnab,\nabla$ the anisotropic derivatives associated with the Chern connections of $L,F$, respectively, it follows from \eqref{DerivAFunct} and Lemma~\ref{connections} that, for any anisotropic function $\psi$ on $M$ regarded as a function of $(\tau,v)$, $(\hnab_{(c_1,X)} \psi)(\tau,v) = (\nabla_X \psi)(v)$. In particular, $(\hnab_{(F(v),v)} F)(F(v),v) = (\nabla_{v} F)(v) = 0$, where $F$ is regarded as a function of $(\tau,v)$. Then, by the above formula, $(\hnab_{(F(v),v)}\varphi) (F(v),v) = \tfrac{1}{F(v)} (\nabla_v \d f) (v)$, and analogously for the second derivative. Therefore, from~\eqref{firstRelation},
%--------------------------------------------------------------------
\begin{equation}
\label{intermediate}
\begin{aligned}
\til{\K}_{(F(v),v)}(-\d f(w),w) = \K_{(F(v),v)}(-\d f(w),w) &- \frac{1}{2} \frac{(\nabla^2_{v} \d f) (v)}{(F(v)+\d f(v))}\\
&+ \frac{3}{4} \frac{(\nabla_{v} \d f)^2(v)}{(F(v)+\d f(v))^2} .
\end{aligned}
\end{equation}
%--------------------------------------------------------------------

On one hand, using \eqref{dphiPi} and \eqref{KsobreK} for the pair $L^f$ and $F+\d f$, we have that
%--------------------------------------------------------------------
\begin{align*}
\til{\K}_{(F(v),v)}(-\d f(w),w) &= \K^f_{(F(v)+\d f(v),v)}(0,w) \\
&= (F(v)+\d f(v))^2 K^{F+\d f}_v(w) . 
\end{align*}
%--------------------------------------------------------------------
On the other hand, applying equation~\eqref{KsobreK} to $L$ and $F$, we get that
%--------------------------------------------------------------------
\begin{equation*}
\K_{(F(v),v)}(-\d f(w),w) = F(v)^2 K^F_v(w) .
\end{equation*}
%--------------------------------------------------------------------
Replacing the last two identities in \eqref{intermediate}, we conclude \eqref{KFdf}.
\end{proof}
%--------------------------------------------------------------------

To the best of our knowledge, this formula is only known when $F$ is the Euclidean metric; see \cite[p. 71]{bao12i}. In practice, there are also other related formulas, as the comparison of the curvature operator for projectively related Finsler metrics (see \cite[Lemma 12.1.2]{she01d}) and the particular cases of Randers metrics considered in \cite{xia13f}.
%--------------------------------------------------------------------
\section{Flag curvature of Randers metrics from \texorpdfstring{$(\alpha,\beta)$}{(α,β)}-data}
\label{Randersalphabeta}
%--------------------------------------------------------------------

There is a similar application to compute the flag curvature of a Randers metric. As there are multiple forms of expressing the data of a Randers metric, we will see how to compute the curvature in each case. To do so, we will need to consider a different expression for the metric of a stationary spacetime. In this section, we will consider a Randers metric given by its $(\alpha,\beta)$-data, that is,
%--------------------------------------------------------------------
\begin{equation}
\label{Randers}
F(v) = \sqrt{g(v,v)} + \omega(v) ,
\end{equation}
%--------------------------------------------------------------------
for a Riemannian metric $g$ on $M$ and a one-form $\omega$ with $g$-norm less than one at every point.

Let us consider on $\R \times M$ the stationary Lorentz metric
%--------------------------------------------------------------------
\begin{equation}
\label{Stationary}
L(\tau,v) = (\tau - \omega(v))^2 - g(v,v) ,
\end{equation}
%--------------------------------------------------------------------
and the Lorentz--Finsler one
%--------------------------------------------------------------------
\begin{equation}
\til{L}(\tau,v) = \tau^2 - (\sqrt{g(v,v)} + \omega(v))^2 .
\end{equation}
%--------------------------------------------------------------------
Both of them have the same lightcone, $\C=\{(F(v),v): v\in TM\}$, where $F$ is given by \eqref{Randers}. Then there exists $\varphi$ such that $\til{L} = \varphi L$.

%--------------------------------------------------------------------
\begin{lemma}
\label{conformalStat}
The conformal factor relating $\til{L}$ and $L$ by $\til{L} = \varphi L$, when restricted to the lightcone, is
%--------------------------------------------------------------------
\begin{equation}
\label{varphiRanders}
\varphi(F(v),v) = \frac{F(v)}{\sqrt{g(v,v)}} = 1 + \frac{\omega(v)}{\sqrt{g(v,v)}} .
\end{equation}
%--------------------------------------------------------------------
\end{lemma}
%--------------------------------------------------------------------
%--------------------------------------------------------------------
\begin{proof}
Let us apply Lemma~\ref{varphiLema} with $\xi=(1,0)$. We compute
%--------------------------------------------------------------------
\begin{equation*}
(1,0)(\til{L})(\tau,v) = 2\tau \text{ and } (1,0)(L)(\tau,v) = 2(\tau-\omega(v))  .
\end{equation*}
%--------------------------------------------------------------------
For lightlike vectors, we conclude
%--------------------------------------------------------------------
\begin{equation*}
\varphi(F(v),v) = \frac{2F(v)}{2(F(v) - \omega(v))} = \frac{F(v)}{\sqrt{g(v,v)}} .
\end{equation*}
%--------------------------------------------------------------------
\end{proof}
%--------------------------------------------------------------------

We can compute the flag curvature of the Randers metric \eqref{Randers} using the function \eqref{varphiRanders} and the null sectional curvature of the stationary manifold \eqref{Stationary}.

%--------------------------------------------------------------------
\begin{proposition}
\label{computingKF}
Let $g$ be a Riemannian metric and $\omega$ a one-form with $g$-norm less than one at every point. The null sectional curvature $\K$ of $L$ in \eqref{Stationary} is related to the flag curvature $K^F$ of the Randers metric $F$ in \eqref{Randers} by
%--------------------------------------------------------------------
\begin{multline}
\label{firstRelationa}
F(v)^2K^{F}_{v}(w) = \K_{(F(v),v)}(\tau,w) - \frac{\sqrt{g(v,v)}}{2F(v)} (\hnab^2_{(F(v),v)}\varphi) (F(v),v)\\
+ \frac{3g(v,v)}{4 F(v)^2}  (\hnab_{(F(v),v)} \varphi)^2(F(v),v) , 
\end{multline}
%--------------------------------------------------------------------
where $\hat \nabla$ denotes the Levi-Civita connection of $L$, $\varphi$ is defined by \eqref{varphiRanders}, and $(\tau,w)$ must be $L$-orthogonal to $(F(v),v)$.
\end{proposition}
%--------------------------------------------------------------------
%--------------------------------------------------------------------
\begin{proof}
$\til{L}$ and $L$ are anisotropically conformal and, by Lemma~\ref{conformalStat}, the conformal factor on lightlike vectors is given by $\varphi(F(v),v) = \frac{F(v)}{\sqrt{g(v,v)}}$. Then, applying \eqref{ConformalFlagCurvEq}, we get
%--------------------------------------------------------------------
\begin{align*}
\til{\K}_{(F(v),v)} (\tau,w) = \K_{(F(v),v)} (\tau,w) & - \frac{\sqrt{g(v,v)}}{2F(v)} (\hnab^2_{(F(v),v)} \varphi) (F(v),v) \\
& + \frac{3 g(v,v)}{4 F(v)^2} (\hnab_{(F(v),v)} \varphi)^2 (F(v),v) ,
\end{align*}
%--------------------------------------------------------------------
where $\til{\K}$ denotes the null flag curvature of $\til{L}$. Taking into account that $\til{\K}_{(F(v),v)} (\tau,w) = F(v)^2 K^F_v(w)$ by \eqref{KsobreK}, we deduce \eqref{firstRelationa}.
\end{proof}
%--------------------------------------------------------------------

%--------------------------------------------------------------------
\begin{proposition}
\label{nablas}
Let $\hnab$ be the Levi-Civita connection of $L$, and $\nabla$ the Levi-Civita connection of $g$. If $\varphi$ is the conformal factor relating $\til{L}$ and $L$ by $\til{L} = \varphi L$, then
%--------------------------------------------------------------------
\begin{subequations}
\begin{equation}
\label{nablaomega}
(\hnab_{(F(v),v)} \varphi) (F(v),v) = \frac{1}{\sqrt{g(v,v)}} \nabla_v \omega (v) - \omega(\d\omega^{\sharp} (v)) ,
\end{equation}
\begin{multline}
\label{nabla2omega}
(\hnab^2_{(F(v),v)} \varphi) (F(v),v) = \\
= \frac{1}{\sqrt{g(v,v)}} \nabla^2_v \omega (v) - \nabla_v (\omega \circ \d\omega^{\sharp}) (v) - \nabla_{\d\omega^{\sharp} (v)} \omega (v) \\
- \nabla_v \omega(\d\omega^{\sharp} (v)) + \sqrt{g(v,v)} \omega(\d\omega^{\sharp}(\d\omega^{\sharp} (v))) ,
\end{multline}
\end{subequations}
%--------------------------------------------------------------------
where $\sharp$ denotes the musical isomorphism of raising one index, i.e. $\d\omega^{\sharp}$ is determined by $g(\d\omega^{\sharp}(X),Y) = \d\omega (X,Y)$, $\forall X,Y \in \X(M)$.
\end{proposition}
%--------------------------------------------------------------------
%--------------------------------------------------------------------
\begin{proof}
First observe that $\hnab$ is precisely the Levi-Civita connection of $\hat{g}$ defined in \eqref{StationaryMetric}, regardless of the sign convention. So all results in Appendix~\ref{stationarysection} are valid in the present context.

Consider an extension $V$ of $v$ such that $\nabla_V V\vert_{\pi(v)} = 0$. Let $\partial_t = (1,0)$ and, for any $X \in \X(M)$, let $\hor{X} = (\omega(X),X)$; see \eqref{lift}. Then $(F(V),V) = \mu\partial_t + \hor{V}$ with $\mu = \sqrt{g(V,V)}$ and
%--------------------------------------------------------------------
\begin{equation*}
\hnab_{(F(V),V)} (F(V),V)\vert_{(0,\pi(v))} = (\mu \hnab_{\partial_t} \hor{V} + \hor{V}(\mu) \partial_t + \mu\hnab_{\hor{V}} \partial_t + \hnab_{\hor{V}} \hor{V})\vert_{(0,\pi(v))} ,
\end{equation*}
%--------------------------------------------------------------------
because $\mu$ is independent of the time coordinate $t$ and $\hnab_{\partial_t} \partial_t = 0$; see proof of Lemma~\ref{AXY}. Once picked $V$ so that $\nabla_V V\vert_{\pi(v)} = 0$, we have $\hor{V}(\mu)\vert_{(0,\pi(v))} = V(\mu)\vert_{\pi(v)} = 0$. Moreover, $\hnab_{\hor{V}} \hor{V}\vert_{(0,\pi(v))} = (\frac{1}{2} \d\omega(V,V) \partial_t + \hor{\nabla_V V})\vert_{(0,\pi(v))} = 0$; refer to \eqref{nabHorVectors}. At the same time, $\hnab_{\hor{V}} \partial_t = \hnab_{\partial_t} \hor{V}$ as $[\partial_t, \hor{V}] = 0$; see \eqref{br1}. It follows that
%--------------------------------------------------------------------
\begin{equation}
\label{nablaFV}
\begin{aligned}
\hnab_{(F(V),V)} (F(V),V)\vert_{(0,\pi(v))} &= 2\mu \hnab_{\partial_t} \hor{V}\vert_{(0,\pi(v))} \\
&= \sqrt{g(v,v)} (\omega(\d\omega^{\sharp}(v)), \d\omega^{\sharp}(v)) ,
\end{aligned}
\end{equation}
%--------------------------------------------------------------------
where we used $\hnab_{\partial_t} \hor{X} = \frac{1}{2} \hor{\d\omega^{\sharp}(X)}$ from Lemma~\ref{lemmaAxi}.

Next, observe that to compute the desired anisotropic derivatives, we can consider any extension of $\varphi$ to the whole tangent bundle. In particular, considering the extension $\varphi(\tau,v)=\frac{\tau}{\tau-\omega(v)}$, we obtain
%--------------------------------------------------------------------
\begin{equation*}
\vd\varphi_{(\tau,v)} (\eta,X) = \frac{\eta(\tau-\omega(v)) - \tau(\eta-\omega(X))}{(\tau-\omega(v))^2} = \frac{-\eta\omega(v) + \tau\omega(X)}{(\tau-\omega(v))^2} .
\end{equation*}
%--------------------------------------------------------------------
So
%--------------------------------------------------------------------
\begin{equation}
\label{partialomega}
\vd\varphi_{(F(v),v)} (\eta,X) = \frac{-\eta\omega(v) + F(v)\omega(X)}{g(v,v)} .
\end{equation}
%--------------------------------------------------------------------
Putting together \eqref{nablaFV} and \eqref{partialomega},
%--------------------------------------------------------------------
\begin{align*}	
&(\hnab_{(F(v),v)} \varphi) (F(v),v) = \\
&= (F(v),v) (\varphi(F(V),V)) - \vd\varphi_{(F(v),v)} (\hnab_{(F(V),V)} (F(V),V)) \\
&= \frac{1}{\sqrt{g(v,v)}} \nabla_v \omega (v) - \frac{\sqrt{g(v,v)} \big( - \omega(\d\omega^{\sharp}(v)) \omega(v) + F(v) \omega(\d\omega^{\sharp}(v)) \big)}{g(v,v)} \\
&= \frac{1}{\sqrt{g(v,v)}} \nabla_v \omega (v) - \omega(\d\omega^{\sharp}(v)) ;
\end{align*}
%--------------------------------------------------------------------
here we have used that $\varphi (F(V),V) = 1 + \frac{\omega(V)}{\sqrt{g(V,V)}}$, which has no dependence on the time coordinate, and that $V(g(V,V))\vert_{\pi(v)} = 0$ and $V(\omega(V))\vert_{\pi(v)} = \nabla_v \omega (v)$, since $\nabla_V V\vert_{\pi(v)}=0$. This concludes \eqref{nablaomega}.

Let us now consider the extension of $\hnab_{(F(v),v)} \varphi (F(v),v)$ to the whole tangent bundle $\phi(\tau,v) = \frac{1}{\sqrt{g(v,v)}} \nabla_v \omega(v) - \omega(\d\omega^{\sharp}(v))$, which does not depend on $\tau$. Then
%--------------------------------------------------------------------
\begin{multline*}
\vd_{(F(v),v)} \phi (\eta,X) = \frac{1}{\sqrt{g(v,v)}} \big( \nabla_X \omega(v) + \nabla_v \omega(X) \big) \\
- \frac{g(v,X)}{g(v,v)\sqrt{g(v,v)}} \nabla_v \omega(v) -\omega(\d\omega^{\sharp}(X)) .
\end{multline*}
%--------------------------------------------------------------------
Using the last identity and \eqref{nablaFV}, we get to \eqref{nabla2omega}:
%--------------------------------------------------------------------
\begin{align*}
&(\hnab^2_{(F(v),v)} \varphi) (F(v),v) = \\
&= (F(v),v)(\phi(F(V),V)) - \vd_{(F(v),v)} \phi (\hnab_{(F(V),V)} (F(V),V)) \\
&=\frac{1}{\sqrt{g(v,v)}} \nabla^2_v \omega (v) - \nabla_v (\omega \circ \d\omega^{\sharp}) (v) - \nabla_{\d\omega^{\sharp}(v)} \omega (v) -\nabla_v \omega(\d\omega^{\sharp}(v)) \\
&\pushright{+ \sqrt{g(v,v)} \omega(\d\omega^{\sharp}(\d\omega^{\sharp}(v))) .}
\end{align*}
%--------------------------------------------------------------------
\end{proof}
%--------------------------------------------------------------------

Observe that in \cite[Lemma 2.1]{lis09r} there is a formula for the flag curvature of a Randers metric in terms of $(\alpha,\beta)$-data. In the next theorem, we obtain a more compact one with an explicit appearance of $\omega$, $\d \omega$ and tensorial covariant derivatives. Moreover, we will use $\nabla^2\omega$ for the covariant derivative $\nabla(\nabla\omega)$.

%--------------------------------------------------------------------
\begin{theorem}
\label{Randerssub}
Let $F$ be a Randers metric on a manifold $M$ given by $F(v)=\sqrt{g(v,v)} + \omega(v)$ with $g$ a Riemannian metric on $M$  and $\omega$ a one-form with $g$-norm less than one at every point. Then the flag curvature of $F$ is given by
%--------------------------------------------------------------------
\begin{multline}
\label{A}
F(v)^2K^{F}_{v}(w) = g(v,v) K(v,w) \\
+ \frac{1}{g(w,w)} \left( \frac{3}{4}\d \omega(v,w)^2 + \sqrt{g(v,v)} (\nabla_{w} \d\omega)(v,w) + \frac{g(v,v)}{4} \Vert \d\omega^\sharp(w) \Vert^2 \right) \\
\pushleft{- \frac{1}{2F(v)} \bigg( \nabla^2_v \omega(v) + g(v,v) \omega(\d\omega^{\sharp}(\d\omega^{\sharp}(v)))} \\
- \sqrt{g(v,v)} \Big( \nabla_v (\omega \circ \d\omega^{\sharp})(v) + \nabla_{\d\omega^{\sharp}(v)} \omega(v) + \nabla_v \omega(\d\omega^{\sharp}(v)) \Big) \bigg) \\
+ \frac{3}{4 F(v)^2} \Big( \nabla_v \omega(v) - \sqrt{g(v,v)} \omega(\d\omega^{\sharp}(v)) \Big)^2 ,
\end{multline}
%--------------------------------------------------------------------
for every $v\in TM\setminus \mathbf{0}$ and $w\in T_{\pi(v)}M\setminus \mathbf{0}$ $g$-orthogonal to $v$, where $\nabla$ is the Levi-Civita connection of $g$, $\Vert \cdot \Vert$ is the norm defined by $g$, $\d \omega$ is the exterior differential of $\omega$, $\d\omega^\sharp$ is the operator metrically equivalent to $\d\omega$ (raising one index), and $K$ is the sectional curvature of $g$. 
\end{theorem}
%--------------------------------------------------------------------
%--------------------------------------------------------------------
\begin{proof}
Let $g^L$ denote the fundamental tensor of the metric $L$ in \eqref{Stationary}. If $\hat{\nabla}$ and $\hat{R}$ are, respectively, the Levi-Civita connection and curvature tensor of $L$, then the metric $\hat{g} = - g^L$ considered in \eqref{StationaryMetric} has $\hat{\nabla}$ and $\hat{R}$ as its Levi-Civita connection and curvature tensor as well.

Consider $\hor{X} = (\omega(X),X)$, for all $X \in \X(M)$, as defined in \eqref{lift}. For $\partial_t = (1,0)$ and $\mu = \sqrt{g(v,v)}$, we have that $\mu\partial_t + \hor{v} = (\sqrt{g(v,v)} + \omega(v), v) = (F(v),v)$ is a lightlike vector for $L$ and, equivalently, for $\hat{g}$. Then the null sectional curvature for $L$ can be computed as
%--------------------------------------------------------------------
\begin{equation*}
\K_{(F(v),v)} (\omega(w),w) = \frac{\hat{g}(\hat{R}((F(v),v), (\omega(w),w)) (\omega(w),w), (F(v),v))}{\hat{g}((\omega(w),w),(\omega(w),w))} ,
\end{equation*}
%--------------------------------------------------------------------
using the symmetries of the curvature tensor $\hat{R}$, regardless of the signature convention for $\hat{g}$. By Proposition~\ref{CurvTen},
%--------------------------------------------------------------------
\begin{multline*}
\hat{g}(\hat{R}((F(v),v), (\omega(w),w)) (\omega(w),w), (F(v),v)) = \\
= g(R(v,w)w,v) + \frac{3}{4}\d \omega(v,w)^2  \\
+ \sqrt{g(v,v)} (\nabla_w \d\omega)(v,w) + \frac{g(v,v)}{4} \Vert \d\omega^{\sharp}(w) \Vert^2 ,
\end{multline*}
%--------------------------------------------------------------------
where $R$ is the curvature tensor of $g$. Taking into account that $g(v,w) = \hat{g}((\omega(w),w), (F(v),v)) = 0$ and $\hat{g}((\omega(w),w), (\omega(w),w)) = g(w,w)$, the above two identities imply 
%--------------------------------------------------------------------
\begin{equation}
\label{precedent}
\begin{aligned}
\K_{(F(v),v)} (\omega(w),w) =& g(v,v) K(v,w) + \frac{1}{g(w,w)} \left( \frac{3}{4}\d \omega(v,w)^2 \right. \\
&\left. + \sqrt{g(v,v)} (\nabla_{w} \d\omega)(v,w) + \frac{g(v,v)}{4} \Vert \d\omega^\sharp(w) \Vert^2 \right) .
\end{aligned}
\end{equation}
%--------------------------------------------------------------------
Substituting \eqref{precedent} and Proposition~\ref{nablas} into \eqref{firstRelationa}, we deduce the expression for the flag curvature $K^F_v(w)$.
\end{proof}
%--------------------------------------------------------------------

%--------------------------------------------------------------------
\begin{remark}
We can deduce an identity for the flag curvature of $F$ even when $w$ is arbitrary. In that case,
%--------------------------------------------------------------------
\begin{equation*}
\til{w} = w - \frac{g(v,w)}{g(v,v)} v
\end{equation*}
%--------------------------------------------------------------------
is $g$-orthogonal to $v$, and we can apply Theorem~\ref{Randerssub} to the pair $v,\til{w}$. The only term that changes formally is
%--------------------------------------------------------------------
\begin{equation*}
\frac{1}{g(w,w)} \left( \frac{3}{4}\d \omega(v,w)^2 + \sqrt{g(v,v)} (\nabla_{w} \d\omega)(v,w) + \frac{g(v,v)}{4} \Vert \d\omega^\sharp(w) \Vert^2 \right) ,
\end{equation*}
%--------------------------------------------------------------------
that becomes
%--------------------------------------------------------------------
\begin{multline*}
\frac{g(v,v)}{g(v,v)g(w,w) - g(v,w)^2} \left( \frac{3}{4}\d \omega(v,w)^2 \right. \\
\left. + \sqrt{g(v,v)} (\nabla_{w - \frac{g(v,w)}{g(v,v)}v} \d\omega)(v,w) + \frac{g(v,v)}{4} \left\Vert \d\omega^\sharp \left( w - \frac{g(v,w)}{g(v,v)}v \right) \right\Vert^2 \right) .
\end{multline*}
%--------------------------------------------------------------------
\end{remark}
%--------------------------------------------------------------------

The next result is closely related to the examples in \cite{che12n}, which have been a source of inspiration.

%--------------------------------------------------------------------
\begin{corollary}
\label{Kahler}
Let $(M,J,g)$ be a K\"ahler manifold with $\hat\omega(u,v) = g(Ju,v)$. Let $\Phi$ be a local function such that $\hat\omega = \frac{i}{2} \partial \bar\partial \Phi$ and consider $\omega(v) = c \d\Phi (v - Jv)$ such that $\|\omega\|_g<1$. Then $\d\omega = 4c \hat\omega$, $\d\omega^{\sharp} = 4cJ$, and
%--------------------------------------------------------------------
\begin{multline*}
F(v)^2 K^{F}_{v}(w) = g(v,v) K(v,w) + 12 c^2\frac{g(Jv,w)^2}{g(w,w)} + 4c^2 g(v,v) \\
- \frac{1}{2F(v)} \Big( c g(\nabla^2_v X, v-J(v)) - 16c^3 g(v,v) g(X, v-Jv) \\
\qquad - 4c^2 \sqrt{g(v,v)} \big( 2 g(\nabla_v X, v+Jv) + g(\nabla_{Jv} X, v - Jv) \big) \Big) \\
+ \frac{3}{4 F(v)^2} \big( c g(\nabla_v X, v-Jv) - 4c^2 \sqrt{g(v,v)} g(X, v+Jv) \big)^2 ,
\end{multline*}
%--------------------------------------------------------------------
where $X = \nabla\Phi$ is the Riemannian gradient of $\Phi$, namely, $\d\Phi(u) = g(\nabla\Phi,u)$, for all $u\in TM$. Moreover, if $g$ has constant negative holomorphic sectional curvature $C$, then taking $c = \frac{\sqrt{-C}}{4}$ in the above setting, we obtain a Randers metric $F = \sqrt{g(v,v)} + \omega(v)$ with scalar flag curvature.
\end{corollary}
%--------------------------------------------------------------------
%--------------------------------------------------------------------
\begin{proof}
In order to check that $\d\omega = 4c\hat\omega$, just observe that $\omega(v) = c (\d\Phi(v) + 2\d^c\Phi(v))$, where $\d^c\Phi(v) = -\frac{1}{2} \d\Phi(Jv)$. It is well-known that $\d = \partial + \bar\partial$ and $\d^c = - \frac{i}{2}(\partial - \bar\partial)$, where $\partial$ and $\bar\partial$ are the Dolbeault operators. Taking into account that $\partial^2 = \bar\partial^2 = \partial\bar\partial + \bar\partial\partial = 0$, we conclude that
%--------------------------------------------------------------------
\begin{equation*}
\d\omega = 2c \d\d^c \Phi = 2ic \partial \bar\partial \Phi = 4c \hat\omega ,
\end{equation*}
%--------------------------------------------------------------------
as required. Moreover, as
%--------------------------------------------------------------------
\begin{equation*}
g(\d\omega^\sharp(u),z) = \d\omega(u,z) = 4c\hat\omega(u,z) = 4cg(Ju,z) ,
\end{equation*}
%--------------------------------------------------------------------
for all $u,z\in TM$, it follows that $\d\omega^{\sharp} = 4cJ$.

To obtain the formula of the flag curvature, we have to use \eqref{A}. Let us compute term by term. Begin by observing that $g$ and $J$ are $\nabla$-parallel, and consequently $\nabla_v \d\omega = 0$. Also taking into account that $\omega(v) = c g(X,v-Jv)$, we compute
%--------------------------------------------------------------------
\begin{align*}
\d\omega(v,w)^2 &= (4c g(Jv,w))^2 = 16c^2 g(Jv,w)^2 ,\\
\Vert \d\omega^{\sharp}(w) \Vert^2 &= g(4cJw, 4cJw) = 16 c^2 g(Jw, Jw) = 16 c^2 g(w,w) ,\\
\omega(\d\omega^{\sharp}(\d\omega^{\sharp}(v))) &= 16c^2 \omega(J^2 v) = -16c^3 g(X, v-Jv) ,\\
\nabla_v \omega(v) &= c g(\nabla_v X, v-Jv) ,\\
\nabla^2_v \omega(v) &= c g(\nabla^2_v X, v-Jv) ,\\
\nabla_v (\omega\circ\d\omega^{\sharp})(v) &= 4c \nabla_v \omega(Jv) = 4c^2 g(\nabla_v X, Jv+v) ,\\
\nabla_{\d\omega^{\sharp}(v)} \omega(v) &= 4c \nabla_{Jv} \omega(v) = 4c^2 g(\nabla_{Jv} X, v-Jv) ,\\
\nabla_v \omega(\d\omega^{\sharp}(v)) &= 4c \nabla_v \omega (Jv) = 4c^2 g(\nabla_v X, Jv+v) ,\\
\omega(\d\omega^{\sharp}(v)) &= 4c \omega(Jv) = 4c^2 g(X, Jv+v) .
\end{align*}
%--------------------------------------------------------------------
We obtain the required expression for $K^F_v (w)$ by replacing all these terms in \eqref{A}. Notice $2 g(\nabla_v X, v+Jv) + g(\nabla_{Jv} X, v - Jv)$ could be replaced by $g(\nabla_v X, 2v+3Jv) - g(\nabla_{Jv} X, Jv)$, because $g(\nabla_uX,z)$ is symmetric in $u$ and $z$ (it is the {\em Riemannian Hessian} of $\Phi$).

For the last statement, observe that if $(M,g,J)$ has constant holomorphic sectional curvature $C$, then
%--------------------------------------------------------------------
\begin{equation*}
K(v,w)=\frac{C}{4}\big(1+3g(v,Jw)^2\big) ,
\end{equation*}
%--------------------------------------------------------------------
for an orthonormal basis $\{v,w\}$; see \cite[Proposition 7.3]{kob69f}. The second term in the sum cancels out with  $12 c^2\frac{g(Jv,w)^2}{g(w,w)}$, which removes all dependence of $K^F_v(w)$ on $w$.
\end{proof}
%--------------------------------------------------------------------
%--------------------------------------------------------------------
\section{Flag curvature of Randers metrics from Zermelo data}
\label{RandersZermelo}
%--------------------------------------------------------------------

In order to compute the flag curvature of a Randers metric with respect to its Zermelo data $(h,W)$, where $h$ is a Riemannian metric and $W$ a vector field with $h$-norm less than one at every point, we will first look at the metric
%--------------------------------------------------------------------
\begin{equation}
\label{RandersStat}
F(v) = \sqrt{\frac{1}{\beta} g_0(v,v) + \frac{1}{\beta^2} g_0(v,\delta)^2} + \frac{1}{\beta} g_0(v,\delta) ,
\end{equation}
%--------------------------------------------------------------------
where $g_0$, $\delta$ and $\beta$ are, respectively, a Riemannian metric, a vector field and a positive smooth function on $M$. This is a well-defined (positive-definite) Randers metric, because the one-form $\omega(v) = \frac{1}{\beta} g_0(v,\delta)$ has norm less than one at every point with respect to the metric $g(v,v) = \frac{1}{\beta} g_0(v,v) + \frac{1}{\beta^2} g_0(v,\delta)^2$. Indeed, as $\beta > 0$, the condition $\frac{g_0(\delta,\delta)}{\beta + g_0(\delta,\delta)} < 1$ is always satisfied. Then we will consider the product manifold $\R\times M$ endowed with the standard stationary Lorentz metric
%--------------------------------------------------------------------
\begin{equation}
\label{stdstationary}
\til{g}((\tau,v),(\tau,v)) = g_0(v,v) + 2 g_0(v,\delta)\tau - \beta\tau^2 .
\end{equation}
%--------------------------------------------------------------------
It is well-known that the Randers metric \eqref{RandersStat} controls the lightlike geodesics and the causality of $\til{g}$; see \cite{cap11i,cap11e}. Moreover, the Finsler spacetime on $\R\times M$ given by
%--------------------------------------------------------------------
\begin{equation*}
L(\tau,v) = \tau^2 - F(v)^2
\end{equation*}
%--------------------------------------------------------------------
has the same lightcone as $\til{L}(\tau,v) = -\til{g}((\tau,v),(\tau,v))$. As a consequence, $L$ and $\til{L}$ are anisotropically conformal Lorentz--Finsler metrics. Observe that in the following,  the roles of $L$ and $\tilde L$ are swapped with respect to the previous sections, namely, $L=\varphi \tilde L$; this requires some care when applying the results on anisotropic conformal changes.

%--------------------------------------------------------------------
\begin{lemma}
\label{conformalStat2}
The conformal factor relating $L$ and $\til{L}$ by $L = \varphi \til{L}$, when restricted to the lightcone, is
%--------------------------------------------------------------------
\begin{equation}
\label{varphiRanders2}
\begin{aligned}
\varphi(F(v),v) &= \frac{F(v)}{\sqrt{\beta g_0(v,v) + g_0(v,\delta)^2}} \\
&= \frac{1}{\beta} \left( 1 + \frac{g_0(v,\delta)}{\sqrt{\beta g_0(v,v) + g_0(v,\delta)^2}} \right) .
\end{aligned}
\end{equation}
%--------------------------------------------------------------------
\end{lemma}
%--------------------------------------------------------------------
%--------------------------------------------------------------------
\begin{proof}
Let us apply Lemma~\ref{varphiLema} with $\xi=(1,0)$. We first compute
%--------------------------------------------------------------------
\begin{equation*}
(1,0)(L)(\tau,v) = 2\tau \text{ and } (1,0)(\til{L})(\tau,v) = 2(\beta\tau - g_0(v,\delta)) .
\end{equation*}
%--------------------------------------------------------------------
As $(\tau,v)$ is lightlike iff $\tau=F(v)$, we conclude
%--------------------------------------------------------------------
\begin{align*}
\varphi(F(v),v) =& \frac{2F(v)}{2(\beta F(v) - g_0(v,\delta))} = \frac{F(v)}{\beta F(v) - g_0(v,\delta)} \\
=& \frac{F(v)}{\sqrt{\beta g_0(v,v) + g_0(v,\delta)^2}} = \frac{1}{\beta} + \frac{g_0(v,\delta)}{\beta \sqrt{\beta g_0(v,v) + g_0(v,\delta)^2}} .
\end{align*}
%--------------------------------------------------------------------
\end{proof}
%--------------------------------------------------------------------

As in the previous section, we will be able to compute the flag curvature of a Randers metric as in \eqref{RandersStat} using the function $\varphi$ above and the null sectional curvature of a stationary manifold as in \eqref{stdstationary} in terms of $g_0$, $\delta$ and $\beta$. To proceed, let us define the auxiliary function
%--------------------------------------------------------------------
\begin{equation*}
E(v) = \sqrt{\beta g_0(v,v)+g_0(v,\delta)^2} .
\end{equation*}
%--------------------------------------------------------------------

%--------------------------------------------------------------------
\begin{proposition}
Given a Riemannian metric $g_0$, a vector field $\delta$ and a positive function $\beta$, the null sectional curvature $\til{\K}$ of $\til{g}$ in \eqref{stdstationary} is related to the flag curvature $K^F$ of the Randers metric $F$ in \eqref{RandersStat} by
%--------------------------------------------------------------------
\begin{multline}
\label{firstRelation2}
F(v)^2 K^{F}_{v}(w) = \til{\K}_{(F(v),v)}(0,w) - \frac{E(v)}{2 F(v)} (\tnab^2_{(F(v),v)} \varphi) (F(v),v) \\
+ \frac{3 E(v)^2}{4 F(v)^2} (\tnab_{(F(v),v)} \varphi)^2 (F(v),v) ,
\end{multline}
%--------------------------------------------------------------------
where $\tnab$ is the Levi-Civita connection of $\til{g}$, $\varphi = \frac{F}{E}$ by \eqref{varphiRanders2}, and $(0,w)$ must be $\til{g}$-orthogonal to $(F(v),v)$.
\end{proposition}
%--------------------------------------------------------------------
%--------------------------------------------------------------------
\begin{proof}
Analogous to the proof of Proposition~\ref{computingKF}, taking into account $\tnab, \til{\K}$ are also the Levi-Civita connection and the null sectional curvatu\-re of $\til{L}$, respectively.
\end{proof}
%--------------------------------------------------------------------

To compute the derivatives, we need to define the function
%--------------------------------------------------------------------
\begin{equation*}
\Lambda = - \frac{1}{\beta + g_0(\delta,\delta)} ,
\end{equation*}
%--------------------------------------------------------------------
and, for $X,Y \in \X(M)$, the operators
%--------------------------------------------------------------------
\begin{equation*}
\begin{aligned}
\Omega_{+} (X,Y) &= g_0(\nabla_X \delta, Y) + g_0(\nabla_Y \delta, X) = g_0(\Omega_{+}^{\sharp}(X), Y) , \\
\Omega_{-} (X,Y) &= g_0(\nabla_X \delta, Y) - g_0(\nabla_Y \delta, X) = g_0(\Omega_{-}^{\sharp}(X), Y) ;
\end{aligned}
\end{equation*}
%--------------------------------------------------------------------
same as \eqref{Lambda} and \eqref{Omega}.

%--------------------------------------------------------------------
\begin{proposition}
If $\tnab$ is the Levi-Civita connection of $\til{g}$ in \eqref{stdstationary}, $\nabla$ is the Levi-Civita connection of $g_0$, and $\nabla\beta$ denotes the gradient of $\beta$ with respect to $g_0$, then
%--------------------------------------------------------------------
\begin{subequations}
\begin{multline}
\label{nablaomegab}
\frac{E(v)}{\Lambda} (\tnab_{(F(v),v)} \varphi) (F(v),v) = \frac{F(v)^2}{2} g_0(\nabla\beta, \delta) - \frac{1}{2}\Omega_{+}(v,v) \\
+ F(v) \big( g_0(\nabla\beta, v) + \Omega_{-}(v,\delta) \big) ,
\end{multline}
\begin{multline}
\label{nabla2omegab}
\frac{E(v)}{\Lambda} (\tnab^2_{(F(v),v)} \varphi) (F(v),v) = \\
= 2 \frac{v(\Lambda)}{\Lambda} H(v) - \frac{F(v)^3}{2} \big( g_0(\nabla\beta,\nabla\beta) + \Omega_{-}(\nabla\beta,\delta) \big) \\
+ \frac{F(v)^2}{2} \big( g_0(\nabla_{v} \nabla\beta, \delta) + g_0(\nabla\beta, \nabla_{v} \delta) + \Omega_{+}(\nabla\beta,v) \big) \\
- F(v)^2 \big( \Omega_{-}(v,\nabla\beta) + \Omega_{-}(\Omega_{-}^{\sharp}(v),\delta) \big) \\
+ F(v) \big( g_0(\nabla_{v} \nabla\beta,v) + (\nabla_{v}\Omega_{-})(v,\delta) + \Omega_{-}(v,\nabla_{v} \delta) \big)\\
+ F(v) \Omega_{+}(\Omega_{-}^{\sharp}(v),v) - \frac{1}{2} (\nabla_{v} \Omega_{+})(v,v) ,
\end{multline}
\end{subequations}
%--------------------------------------------------------------------
for $E, \Lambda, \Omega_{+}, \Omega_{-}$ defined above, and $H(v) = \frac{E(v)}{\Lambda} (\tnab_{(F(v),v)} \varphi) (F(v),v)$ expressed by \eqref{nablaomegab}.
\end{proposition}
%--------------------------------------------------------------------
%--------------------------------------------------------------------
\begin{proof}
Let $\alpha(s) = (t(s),x(s))$ be a lightlike geodesic for $\til{g}$ such that $\dot\alpha(0) = (F(v),v)$. As a first observation, we can express $E(v) = \beta F(v) - g_0(v, \delta)$. Second, as $\dot\alpha(s)$ is lightlike iff $\dot t(s) = F(\dot x(s))$, we notice that
%--------------------------------------------------------------------
\begin{equation*}
E(\dot x(s)) = \beta \dot t(s) - g_0(\dot x(s), \delta) = - \til{g}(\dot\alpha(s), \partial_t)
\end{equation*}
%--------------------------------------------------------------------
must be constant, because $\partial_t = (1,0)$ is a Killing vector field. Now, the conformal function $\varphi$ along $\alpha(s)$ reads
%--------------------------------------------------------------------
\begin{equation*}
\varphi(\dot t(s),\dot x(s)) = \frac{1}{\beta} \left( 1 + \frac{g_0(\dot x(s), \delta)}{E} \right) = \frac{F(\dot x(s))}{E} = \frac{\dot t(s)}{E} ,
\end{equation*}
%--------------------------------------------------------------------
for $E$ constant. Then,
%--------------------------------------------------------------------
\begin{align*}
&(\tnab_{(F(v),v)} \varphi)(F(v),v) = \left.\frac{\d}{\d s} \varphi(\dot t(s), \dot x(s))\right\vert_{s=0} = \frac{\ddot t(0)}{E} ,\\
&(\tnab^2_{(F(v),v)} \varphi)(F(v),v) = \left.\frac{\d^2}{\d s^2} \varphi(\dot t(s), \dot x(s))\right\vert_{s=0 }= \frac{t^{(3)}(0)}{E} ,
\end{align*}
%--------------------------------------------------------------------
and we only need to compute $\ddot t$ and $t^{(3)}$.

A lightlike geodesic $\alpha(s)=(t(s),x(s))$ is determined by
%--------------------------------------------------------------------
\begin{equation}
\label{nablaxx}
\begin{cases}
\nabla_{\dot x} \dot x = - \ddot t \delta - \frac{1}{2} \dot t^2 \nabla\beta - \dot t \Omega_{-}^{\sharp}(\dot x) \\
\beta \dot t - g_0(\dot x,\delta) = E \\
\dot t = F(\dot x)
\end{cases};
\end{equation}
%--------------------------------------------------------------------
cf. \cite{flo03g}. Differentiating the second equation with respect to $s$, having in mind that
%--------------------------------------------------------------------
\begin{equation*}
\frac{\d}{\d s} \beta(x(s)) = g_0(\nabla\beta, \dot x(s)) ,
\end{equation*}
%--------------------------------------------------------------------
we obtain
%--------------------------------------------------------------------
\begin{equation*}
g_0(\nabla\beta, \dot x) \dot t + \beta \ddot t - g_0(\nabla_{\dot x} \dot x, \delta) - g_0(\dot x, \nabla_{\dot x} \delta) = 0 .
\end{equation*}
%--------------------------------------------------------------------
Using the first equation in \eqref{nablaxx} to replace $\nabla_{\dot x} \dot x$, it follows that
%--------------------------------------------------------------------
\begin{equation*}
\frac{1}{\Lambda} \ddot t = \frac{1}{2} \dot t^2 g_0(\nabla\beta, \delta) + \dot t \big( g_0(\nabla\beta, \dot x) + \Omega_{-}(\dot x, \delta) \big) - \frac{1}{2} \Omega_{+}(\dot x,\dot x) .
\end{equation*}
%--------------------------------------------------------------------
Evaluating at $s=0$, we get \eqref{nablaomegab}.

Differentiating the above identity with respect to $s$ one more time, we have
%--------------------------------------------------------------------
\begin{align*}
- \frac{\dot x(\Lambda)}{\Lambda^2}\ddot t + \frac{1}{\Lambda} t^{(3)} =& \dot t \ddot t g_0(\nabla\beta, \delta) + \frac{1}{2} \dot t^2 g_0(\nabla_{\dot x} \nabla\beta, \delta) + \frac{1}{2} \dot t^2 g_0(\nabla\beta, \nabla_{\dot x} \delta) \\
& + \ddot t \big( g_0(\nabla\beta, \dot x) + \Omega_{-}(\dot x,\delta) \big) \\
& + \dot t \big( g_0(\nabla_{\dot x} \nabla\beta, \dot x) + g_0(\nabla\beta, \nabla_{\dot x} \dot x) \big) \\
& + \dot t \big( (\nabla_{\dot x} \Omega_{-})(\dot x,\delta) + \Omega_{-}(\nabla_{\dot x} \dot x,\delta) + \Omega_{-}(\dot x,\nabla_{\dot x} \delta) \big) \\
& - \frac{1}{2} (\nabla_{\dot x} \Omega_{+})(\dot x,\dot x) - \Omega_{+}(\nabla_{\dot x} \dot x,\dot x) ,
\end{align*}
%--------------------------------------------------------------------
because $\Omega_{+}$ is symmetric. Reorganizing the terms,
%--------------------------------------------------------------------
\begin{align*}
\frac{1}{\Lambda} t^{(3)} =& \ddot t \left( \frac{\dot x(\Lambda)}{\Lambda^2} + g_0(\nabla\beta, \dot x + \dot t \delta) + \Omega_{-}(\dot x,\delta) \right) \\
& + \frac{1}{2} \dot t^2 \big( g_0(\nabla_{\dot x} \nabla\beta,\delta) + g_0(\nabla\beta, \nabla_{\dot x} \delta) \big) - \frac{1}{2} (\nabla_{\dot x} \Omega_{+})(\dot x,\dot x) \\
& + \dot t \big( g_0(\nabla_{\dot x} \nabla\beta,\dot x) + (\nabla_{\dot x} \Omega_{-})(\dot x,\delta) + \Omega_{-}(\dot x,\nabla_{\dot x} \delta) \big) \\
& + \dot t \big( g_0(\nabla\beta, \nabla_{\dot x} \dot x) + \Omega_{-}(\nabla_{\dot x} \dot x, \delta) \big) - \Omega_{+}(\nabla_{\dot x} \dot x, \dot x) .
\end{align*}
%--------------------------------------------------------------------
Replacing $\nabla_{\dot x}\dot x$ from \eqref{nablaxx} in the last line, we get
%--------------------------------------------------------------------
\begin{align*}
\dot t &\big( g_0(\nabla\beta, \nabla_{\dot x} \dot x) + \Omega_{-}(\nabla_{\dot x} \dot x, \delta) \big) - \Omega_{+}(\nabla_{\dot x} \dot x, \dot x) = \\
=& \ddot t \big( \Omega_{+}(\delta,\dot x) - \dot t g_0(\nabla\beta,\delta) \big) - \frac{1}{2} \dot t^3 \big( g_0(\nabla\beta, \nabla\beta) + \Omega_{-}(\nabla\beta,\delta) \big) \\
&- \dot t^2 \big( \Omega_{-}(\dot x,\nabla\beta) + \Omega_{-}(\Omega_{-}^{\sharp}(\dot x),\delta) \big) + \frac{1}{2} \dot t^2 \Omega_{+}(\nabla\beta,\dot x) + \dot t \Omega_{+}(\Omega_{-}^{\sharp}(\dot x),\dot x)
\end{align*}
%--------------------------------------------------------------------
as $\Omega_{-}(\delta,\delta) = 0$ (for $\Omega_{-}$ is antisymmetric). It follows that
%--------------------------------------------------------------------
\begin{align*}
\frac{1}{\Lambda} t^{(3)} =& \ddot t \left( \frac{\dot x(\Lambda)}{\Lambda^2} + g_0(\nabla\beta, \dot x) + 2g_0(\nabla_{\dot x} \delta, \delta) \right) \\
&- \frac{1}{2} \dot t^3 \big( g_0(\nabla\beta, \nabla\beta) + \Omega_{-}(\nabla\beta,\delta) \big) \\
&+ \frac{1}{2} \dot t^2 \big( g_0(\nabla_{\dot x} \nabla\beta,\delta) + g_0(\nabla\beta, \nabla_{\dot x} \delta) + \Omega_{+}(\nabla\beta,\dot x) \big) \\
&- \dot t^2 \big( \Omega_{-}(\dot x,\nabla\beta) + \Omega_{-}(\Omega_{-}^{\sharp}(\dot x),\delta) \big) \\
&+ \dot t \big( g_0(\nabla_{\dot x} \nabla\beta,\dot x) + (\nabla_{\dot x} \Omega_{-})(\dot x,\delta) + \Omega_{-}(\dot x,\nabla_{\dot x} \delta) \big) \\
&+ \dot t \Omega_{+}(\Omega_{-}^{\sharp}(\dot x),\dot x) - \frac{1}{2} (\nabla_{\dot x} \Omega_{+})(\dot x,\dot x) ,
\end{align*}
%--------------------------------------------------------------------
using the identity $\Omega_{-}(\dot x,\delta) + \Omega_{+}(\delta,\dot x) = 2g_0(\nabla_{\dot x} \delta, \delta)$. Evaluating at $s=0$, gives \eqref{nabla2omegab}. Indeed, the coefficient of $\ddot t$ can be simplified as
%--------------------------------------------------------------------
\begin{align*}
\frac{\dot x(\Lambda)}{\Lambda^2} + g_0(\nabla\beta, \dot x) + 2g_0(\nabla_{\dot x} \delta, \delta) &= \frac{\dot x(\Lambda)}{\Lambda^2}+\dot x(\beta+g_0(\delta,\delta)) \\
&= \frac{\dot x(\Lambda)}{\Lambda^2} - \dot x\left( \frac{1}{\Lambda} \right) = 2\frac{\dot x(\Lambda)}{\Lambda^2}.
\end{align*}
%--------------------------------------------------------------------
\end{proof}
%--------------------------------------------------------------------

We compute the null sectional curvature $\til{\K}$ of $\til{g}$ in Appendix~\ref{stdstationarysection}. Namely, taking $\mu = F(v)$, $X=v$, and $Y=w$ in \eqref{finalresult} yields
%--------------------------------------------------------------------
\begin{multline*}
g_0(w,w)\til{\K}_{(F(v),v)}(0,w) = \til{g}(\til{R}((F(v),v), (0,w)) (0,w), (F(v),v)) \\
= \frac{\Lambda}{4} \Omega_{+}(w,w) \big( F(v) ((v + F(v)\delta)(\beta) + \Omega_{-}(v,\delta)) - \Omega_{+}(v,v + F(v)\delta) \big) \\
- \frac{\Lambda}{4} \Omega_{+}(w,v + F(v)\delta) \big( F(v) (w(\beta) + \Omega_{-}(w,\delta)) - \Omega_{+}(v,w) \big) \\
+ \frac{F(v)}{4\Lambda} \Big( v(\Lambda)\Omega_{+}(w,w) + w(\Lambda) \big( F(v) (w(\beta) + \Omega_{-}(w,\delta)) - \Omega_{+}(v,w) \big) \Big) \\
+ \frac{F(v)}{2} \big( (\nabla_v \Omega_{+})(w,w) - (\nabla_w \Omega_{+})(v,w) \big)\\
+ \frac{F(v)}{2} \big( (\nabla_{v + F(v)\delta} \Omega_{+})(w,w) - (\nabla_w \Omega_{+})(v + F(v)\delta,w) \big) \\
+ \frac{F(v)^2}{2} \big( (\nabla_w \d\beta)(w) + (\nabla_w \Omega_{-})(w,\delta) + \Omega_{-}(w,\nabla_w \delta) \big) \\
- \frac{F(v)^2}{4} \Omega_{+}(w,\Omega_{-}^{\sharp}(w)) + g_0(R(v + F(v)\delta, w) w, v + F(v)\delta) ,
\end{multline*}
%--------------------------------------------------------------------
for any $v \in TM \setminus \mathbf{0}$ and $w \in T_{\pi(v)} M \setminus \mathbf{0}$ $g_0$-orthogonal to $v + F(v)\delta$, where $R$ is the curvature tensor of $g_0$. Substituting the above equality, \eqref{nablaomegab} and \eqref{nabla2omegab} into \eqref{firstRelation2} produces an expression for the flag curvature of the Randers metric as in \eqref{RandersStat}, but it is unnecessarily long. For the sake of simplicity, we write any Randers metric with respect to its Zermelo data $(h,W)$ as
%--------------------------------------------------------------------
\begin{equation}
\label{RanderZermeloData}
F(v) = \sqrt{\frac{1}{\lambda} h(v,v) + \frac{1}{\lambda^2} h(v,W)^2} - \frac{1}{\lambda} h(v,W) ,
\end{equation}
%--------------------------------------------------------------------
where $\lambda = 1 - h(W,W)$ and $h(W,W) < 1$. This can be obtained from \eqref{RandersStat} by taking $g_0 = h$, $\delta = -W$, and $\beta = \lambda$. In this case, $\beta + g_0(\delta, \delta) = 1$, $\Lambda = -1$ and, for all $X,Y \in \X(M)$,
%--------------------------------------------------------------------
\begin{equation*}
\Omega(X,Y) = h(\nabla_X W, Y) + h(\nabla_Y W, X) = - \Omega_{+}(X,Y) .
\end{equation*}
%--------------------------------------------------------------------

%--------------------------------------------------------------------
\begin{theorem}
\label{ZermeloFlagCurvature}
Let $F$ be a Randers metric on a manifold $M$ with Zermelo navigation data $(h, W)$. Then the flag curvature of $F$ is given by
%--------------------------------------------------------------------
\begin{multline*}
F(v)^2 K^F_v(w) = \Vert v - F(v)W \Vert^2_h K(v - F(v)W, w) \\
+ \frac{1}{4 \Vert w \Vert_h^2} \big( \Omega(w,w) \Omega(v - F(v)W, v - F(v)W)  - \Omega(w, v - F(v)W)^2 \big) \\
+ \frac{F(v)}{\Vert w \Vert_h^2} \big( (\nabla_w \Omega)(v - F(v)W, w) - (\nabla_{v - F(v)W} \Omega)(w,w) \big) \\
+ \frac{F(v)^2}{2 \Vert w \Vert_h^2} \left( \frac{1}{2} \Vert \Omega^{\sharp}(w) \Vert^2_h - 2\Omega(w,\nabla_w W) - (\nabla_W \Omega)(w,w) \right) \\
- \frac{1}{2} \Omega(v - F(v)W, \Omega^\sharp (v - F(v)W) - \nabla_{v - F(v)W} W) \\
+ \frac{1}{4 F(v)} (\nabla_v \Omega)(v - F(v)W, v - F(v)W) \\
+ \frac{3}{16 F(v)^2} \Omega(v - F(v)W, v - F(v)W)^2 ,
\end{multline*}
%--------------------------------------------------------------------
for all $v \in TM \setminus \mathbf{0}$ and $w \in T_{\pi(v)}M \setminus \mathbf{0}$ $h$-orthogonal to $v-F(v)W$, where $\nabla$ is the Levi-Civita connection of $h$, $K$ is the sectional curvature of $h$, $\Vert \cdot \Vert_h$ the norm associated with $h$, $\Omega$ is the Lie derivative of $h$ with respect to $W$ (explicitly defined above), and $\Omega^\sharp$ is the operator $h$-metrically equivalent to $\Omega$ (raising one index).
\end{theorem}
%--------------------------------------------------------------------
%--------------------------------------------------------------------
\begin{proof}
Suppose the metric $F$ in \eqref{RandersStat} satisfies $\beta + g_0(\delta, \delta) = 1$ and $g_0(\delta,\delta) < 1$. Then $\Lambda = -1$, and $\d\beta, \nabla\beta$ can be expressed as
%--------------------------------------------------------------------
\begin{subequations}
\begin{equation*}
\d\beta(X) = -2 g_0(\nabla_X \delta, \delta) = -2(\Omega_{+}(\delta,X) - g_0(\nabla_{\delta} \delta, X)) , \; \forall X \in \X(M) ,
\end{equation*}
\begin{equation*}
\nabla\beta = - 2 (\Omega_{+}^{\sharp}(\delta) - \nabla_{\delta} \delta) = 2 (\nabla_{\delta} \delta - \Omega_{+}^{\sharp}(\delta)).
\end{equation*}
\end{subequations}
%--------------------------------------------------------------------
Also take into consideration the general identity 
%--------------------------------------------------------------------
\begin{equation}
\label{omegassum}
\Omega_{+}(X,Y) + \Omega_{-}(X,Y) = 2g_0(\nabla_X \delta, Y) , \; \forall X,Y \in \X(M) .
\end{equation}
%--------------------------------------------------------------------
So \eqref{finalresult} with $\mu = F(v)$, $X=v$ and $Y=w$ simplifies to:
%--------------------------------------------------------------------
\begin{align*}
\til{g} (\til{R}((&F(v),v), (0,w)) (0,w), (F(v),v)) = \\
=& \frac{1}{4} \big( \Omega_{+}(w,w) \Omega_{+}(v+F(v)\delta, v + F(v)\delta) - \Omega_{+}(w, v + F(v)\delta)^2 \big) \\
&+ \frac{F(v)}{2} \big( (\nabla_v \Omega_{+})(w,w) - (\nabla_w \Omega_{+})(v,w) \big)\\
&+ \frac{F(v)}{2} \big( (\nabla_{v + F(v)\delta} \Omega_{+})(w,w) - (\nabla_w \Omega_{+})(v + F(v)\delta,w) \big) \\
&+ \frac{F(v)^2}{2} \big( (\nabla_w \d\beta)(w) + (\nabla_w \Omega_{-})(w,\delta) + \Omega_{-}(w,\nabla_w \delta) \big) \\
&- \frac{F(v)^2}{4} \Omega_{+}(w,\Omega_{-}^{\sharp}(w)) + g_0(R(v + F(v)\delta, w) w, v + F(v)\delta) ,
\end{align*}
%--------------------------------------------------------------------
where $\til{R}, R$ are the curvature tensors of $\til{g},g_0$, respectively. To simplify the above formula, we can choose $Y \in \X(M)$ an extension of $w \in T_p M$ satisfying $\nabla_Y Y\vert_p = 0$. Now,
%--------------------------------------------------------------------
\begin{align*}
(\nabla_w \d\beta)(w) + &(\nabla_w \Omega_{-})(w,\delta) + \Omega_{-}(w,\nabla_w \delta) = \\
&= \big( (\nabla_Y \d\beta)(Y) + (\nabla_Y \Omega_{-})(Y,\delta) + \Omega_{-}(Y,\nabla_Y \delta) \big)\vert_p \\
&= Y ( \d\beta(Y) + \Omega_{-}(Y,\delta) )\vert_p \\
&= Y ( -\Omega_{+}(Y,\delta) )\vert_p \\
&= \big( -(\nabla_Y \Omega_{+})(Y,\delta) -\Omega_{+}(Y,\nabla_Y \delta) \big)\vert_p \\
&= -(\nabla_w \Omega_{+})(w,\delta) -\Omega_{+}(w,\nabla_w \delta) .
\end{align*}
%--------------------------------------------------------------------
Moreover, by \eqref{omegassum}, we have
%--------------------------------------------------------------------
\begin{equation*}
\Omega_{+}^{\sharp}(X) + \Omega_{-}^{\sharp}(X) = 2\nabla_X \delta , \; \forall X \in \X(M) ,
\end{equation*}
%--------------------------------------------------------------------
and it follows that
%--------------------------------------------------------------------
\begin{align*}
\Omega_{+}(w,\Omega_{-}^{\sharp}(w)) &= \Omega_{+}(w, 2\nabla_w \delta - \Omega_{+}^{\sharp}(w)) \\
&= 2 \Omega_{+}(w,\nabla_w \delta) - g_0(\Omega_{+}^{\sharp}(w),\Omega_{+}^{\sharp}(w)) .
\end{align*}
%--------------------------------------------------------------------
Substituting these new expressions, adding and subtracting 
%--------------------------------------------------------------------
\begin{multline*}
\frac{F(v)^2}{2} \big( (\nabla_{\delta} \Omega_{+})(w,w) - (\nabla_w \Omega_{+})(\delta,w) \big) = \\
= \frac{F(v)}{2} \big( (\nabla_{F(v)\delta} \Omega_{+})(w,w) - (\nabla_w \Omega_{+})(F(v)\delta,w) \big) 
\end{multline*}
%--------------------------------------------------------------------
into the curvature equation, and observing $\nabla_w \Omega_{+}$ is symmetric, we conclude
%--------------------------------------------------------------------
\begin{equation}
\label{final1}
\begin{aligned}
\til{g} (\til{R}((&F(v),v), (0,w)) (0,w), (F(v),v)) = \\
=& \frac{1}{4} \big( \Omega_{+}(w,w) \Omega_{+}(v+F(v)\delta, v + F(v)\delta) - \Omega_{+}(w, v + F(v)\delta)^2 \big) \\
& + F(v) \big( (\nabla_{v + F(v)\delta}\Omega_{+})(w,w) - (\nabla_w \Omega_{+})(v + F(v)\delta, w) \big) \\
& + \frac{F(v)^2}{2} \left( \frac{1}{2} \Vert \Omega_{+}^{\sharp}(w) \Vert_0^2 - 2 \Omega_{+}(w, \nabla_w \delta) - (\nabla_{\delta} \Omega_{+})(w,w) \right) \\
& + g_0(R(v + F(v)\delta, w) w, v + F(v)\delta) ,
\end{aligned}
\end{equation}
%--------------------------------------------------------------------
where $\Vert \cdot \Vert_{0}$ is the norm determined by $g_0$.

Meanwhile, after some long but elementary computations, \eqref{nablaomegab} and \eqref{nabla2omegab} become:
%--------------------------------------------------------------------
\begin{subequations}
\begin{align*}
E(v)(\tnab_{(F(v),v)} \varphi)(F(v),v) &= \tfrac{F(v)^2}{2} \Omega_{+}(\delta,\delta) + \tfrac{1}{2} \Omega_{+}(v,v) + F(v) \Omega_{+}(v, \delta) \\
&= \frac{1}{2}\Omega_{+}(v + F(v)\delta,v + F(v)\delta) ,
\end{align*}
\begin{align*}
E(v)&(\tnab^2_{(F(v),v)} \varphi)(F(v),v) = \\
= & F(v)^3 \Omega_{+}(\delta, \Omega_{+}^{\sharp}(\delta) - \nabla_{\delta} \delta) \\
& + F(v)^2 \bigg( \Omega_{+}(v, \Omega_{+}^{\sharp}(\delta) - \nabla_{\delta} \delta) + \Omega_{+}(\delta, \Omega_{+}^{\sharp}(v) - \nabla_v \delta) \\
&\pushright{+ \frac{1}{2}(\nabla_v \Omega_{+})(\delta,\delta) \bigg)} \\
& + F(v) \left( (\nabla_v \Omega_{+})(\delta,v) + \Omega_{+}(v, \Omega_{+}^{\sharp}(v) - \nabla_v \delta) \right) + \frac{1}{2} (\nabla_{v}\Omega_{+})(v,v) \\
= & \frac{1}{2} (\nabla_v \Omega_{+})(v + F(v)\delta, v + F(v)\delta) \\
& + F(v) \Omega_{+}(v + F(v)\delta, \Omega_{+}^\sharp (v + F(v)\delta) - \nabla_{v + F(v)\delta} \delta) .
\end{align*}
\end{subequations}
%--------------------------------------------------------------------

Putting \eqref{final1} and the previous two equations together in \eqref{firstRelation2} for $v,w$ such that $g_0(v+F(v)\delta,w) = \til{g}( (F(v),v), (0,w)) = 0$ yields the flag curvature of the Randers metric $F$ in \eqref{RandersStat}, under the conditions $\beta + g_0(\delta,\delta) = 1$ and $g_0(\delta,\delta) <1$; namely,
%--------------------------------------------------------------------
\begin{multline*}
F(v)^2 K^F_v(w) = \Vert v + F(v)\delta \Vert^2_0 K(v + F(v)\delta, w) \\
+ \frac{1}{4 \Vert w \Vert_0^2} \big( \Omega_{+}(w,w) \Omega_{+}(v + F(v)\delta, v + F(v)\delta)  - \Omega_{+}(w, v + F(v)\delta)^2 \big) \\
+ \frac{F(v)}{\Vert w \Vert_0^2} \big( (\nabla_{v + F(v)\delta} \Omega_{+})(w,w) - (\nabla_w \Omega_{+})(v + F(v)\delta, w) \big) \\
+ \frac{F(v)^2}{2 \Vert w \Vert_0^2} \left( \frac{1}{2} \Vert \Omega_{+}^{\sharp}(w) \Vert^2_0 - 2\Omega_{+}(w,\nabla_w \delta) - (\nabla_\delta \Omega_{+})(w,w) \right) \\
- \frac{1}{2} \Omega_{+}(v + F(v)\delta, \Omega_{+}^\sharp (v + F(v)\delta) - \nabla_{v + F(v)\delta} \delta) \\
-\frac{1}{4 F(v)} (\nabla_v \Omega_{+})(v + F(v)\delta, v + F(v)\delta) \\
+ \frac{3}{16 F(v)^2} \Omega_{+}(v + F(v)\delta, v + F(v)\delta)^2.
\end{multline*}
%--------------------------------------------------------------------
The result for a Randers metric $F$ in terms of Zermelo navigation data $(h,W)$ as in \eqref{RanderZermeloData} follows from direct substitution of $g_0 = h$, $\delta = -W$ and $\Omega_{+} = -\Omega$ in the previous equation.
\end{proof}
%--------------------------------------------------------------------

To shorten the expression of the flag curvature, since we have positive homogeneity, it is worthwhile to obtain a version of this equation with $h$-unit vectors.

%--------------------------------------------------------------------
\begin{remark}
Let $v \in TM \setminus \mathbf{0}$ and $w \in T_{\pi(v)}M \setminus \mathbf{0}$ be such that $w$ is $h$-orthogonal to $v-F(v)W$. Denote by $\hat u$ and $\hat w$ the $h$-unit vectors in the directions $u = v-F(v)W$ and $w$, respectively. Then
%--------------------------------------------------------------------
\begin{equation}
\label{ZermeloFlagCurvatureUnit}
\begin{aligned}
K^F_{v} (w) &= K(v-F(v)W, w) \\
& +\frac{1}{4} \left( \Omega(\hat u, \hat u) \Omega(\hat w, \hat w) - \Omega(\hat u, \hat w)^2 + \frac{3}{4} \Omega(\hat u,\hat u)^2 \right) \\
& -\frac{1}{2} \Vert \Omega^{\sharp}(\hat u) \Vert_{h}^2 + \frac{1}{4} \Vert \Omega^{\sharp}(\hat w) \Vert_{h}^2\\
& +\frac{1}{2} \Omega(\hat u,\nabla_{\hat u}W) - \Omega(\hat w,\nabla_{\hat w}W) \\
& + (\nabla_{\hat w}\Omega)(\hat u, \hat w) - (\nabla_{\hat u}\Omega)(\hat w,\hat w) + \frac{1}{4}(\nabla_{\hat u}\Omega)(\hat u, \hat u) \\
&+ \frac{1}{4}(\nabla_{W}\Omega)(\hat u, \hat u) - \frac{1}{2}(\nabla_{W}\Omega)(\hat w,\hat w) ,
\end{aligned}
\end{equation}
%--------------------------------------------------------------------
as a direct consequence of Theorem~\ref{ZermeloFlagCurvature}, taking into account that $F(v) = \Vert v-F(v)W \Vert_{h}$; recall $(F(v),v)$ is lightlike for $\til{g}$, which amounts to $h(v - F(v)W,v - F(v)W) - F(v)^2 = \til{g} ((F(v),v), (F(v),v)) = 0$.
\end{remark}
%--------------------------------------------------------------------

The flag curvature of the Randers metric $F$ in \eqref{RanderZermeloData} with Zermelo data $(h,W)$ is further simplified when $W$ is conformal on $(M,h)$ (cf. \cite{hua15f, she13c, xia13f}) or when it is projective on $(M,h)$.

%--------------------------------------------------------------------
\begin{corollary}
\label{conformalW}
If $W$ is conformal on $(M,h)$, i.e. $\mathcal{L}_W h = \sigma h$, for some $\sigma \in C^{\infty}(M)$, then
%--------------------------------------------------------------------
\begin{equation*}
K^F_v(w) = K(v - F(v)W, w) - \frac{\sigma^2}{16} - \frac{3}{4} \frac{\d\sigma(v)}{F(v)} + \frac{1}{2} \d\sigma(W) ,
\end{equation*}
%--------------------------------------------------------------------
where $v \in TM \setminus \mathbf{0}$ and $w \in T_{\pi(v)}M \setminus \mathbf{0}$ is $h$-orthogonal to $v-F(v)W$.
\end{corollary}
%--------------------------------------------------------------------
%--------------------------------------------------------------------
\begin{proof}
Since $\Omega = \mathcal{L}_W h = \sigma h$, we have
%--------------------------------------------------------------------
\begin{align*}
&\Omega^{\sharp}(X) = \sigma X , \\
&\Omega(X, \nabla_Y W) + \Omega(Y, \nabla_X W) = \sigma^2 h(X,Y), \\
&\nabla_X \Omega (Y,Z) = X(\sigma) h(Y,Z) ,
\end{align*}
%--------------------------------------------------------------------
for all $X,Y,Z \in \X(M)$. The result is straightforward from \eqref{ZermeloFlagCurvatureUnit}.
\end{proof}
%--------------------------------------------------------------------

%--------------------------------------------------------------------
\begin{remark}
We stress two particular cases of conformal vector fields.

When $W$ is homothetic (i.e. conformal with $\sigma$ constant),
%--------------------------------------------------------------------
\begin{equation*}
K^F_v(w) = K(v - F(v)W, w) - \frac{\sigma^2}{16} .
\end{equation*}
%--------------------------------------------------------------------

When $W$ is Killing (i.e. conformal with $\sigma = 0$),
%--------------------------------------------------------------------
\begin{equation*}
K^F_v(w) = K(v - F(v)W,w) .
\end{equation*}
%--------------------------------------------------------------------
These formulas are consistent with the literature; see \cite{bao04z,mo07o}.
\end{remark}
%--------------------------------------------------------------------

%--------------------------------------------------------------------
\begin{corollary}
\label{projectiveW}
Suppose $W$ is projective on $(M,h)$, i.e. there exists $\mu \in \X^\ast(M)$ such that $\mathcal{L}_W \nabla (X,Y) = \mu(X)Y + \mu(Y)X$, for all $X,Y\in\mathfrak{X}(M)$. Then
%--------------------------------------------------------------------
\begin{align*}
K^F_v(w) &= K(u,w) - \frac{1}{2}\mu(W) + \frac{1}{2}\mu(\hat u) h(W,\hat u) - \mu(\hat w) h(W,\hat w) \\
& + \frac{1}{4} \left( \Omega(\hat u, \hat u) \Omega(\hat w, \hat w) - \Omega(\hat u,\hat w)^2 + \frac{3}{4} \Omega(\hat u, \hat u)^2 \right) \\
& - \frac{1}{2} \Vert \Omega^{\sharp}(\hat u) \Vert_{h}^2 + \frac{1}{4} \Vert \Omega^{\sharp}(\hat w) \Vert_{h}^2 + \frac{1}{2} \Omega(\hat u, \nabla_{\hat u} W) - \Omega(\hat w, \nabla_{\hat w} W) ,
\end{align*}
%--------------------------------------------------------------------
where $v \in TM \setminus \mathbf{0}$, $w \in T_{\pi(v)}M \setminus \mathbf{0}$ is $h$-orthogonal to $u=v-F(v)W$, and $\hat{z}$ denotes the $h$-unit vector in the direction of $z$.
\end{corollary}
%--------------------------------------------------------------------
%--------------------------------------------------------------------
\begin{proof}
As the Levi-Civita connection $\nabla$ is torsion-free,
%--------------------------------------------------------------------
\begin{align*}
\mathcal{L}_Z \nabla (X,Y) &= [Z,\nabla_X Y] - \nabla_{[Z,X]} Y - \nabla_X [Z,Y] \\
&= \nabla_X \nabla_Y Z - \nabla_{\nabla_X Y} Z + R(Z,X)Y ,
\end{align*}
%--------------------------------------------------------------------
for all $X,Y,Z \in \X(M)$. Using the metric compatibility of $\nabla$ and the antisymmetry of the curvature tensor $h(R(X,Y)U,V)$ in $U$ and $V$, it follows that
%--------------------------------------------------------------------
\begin{align*}
\nabla_Z \Omega(X,Y) &= Z(\Omega(X,Y)) - \Omega(\nabla_Z X, Y) - \Omega(X, \nabla_Z Y) \\
&= h(\nabla_Z \nabla_X W - \nabla_{\nabla_Z X} W, Y) + h(\nabla_Z \nabla_Y W - \nabla_{\nabla_Z Y} W, X) \\
&= h(\mathcal{L}_W \nabla (Z,X), Y) + h(\mathcal{L}_W \nabla (Z,Y), X) ,
\end{align*}
%--------------------------------------------------------------------
for all $X,Y,Z \in \X(M)$. So the projective condition implies
%--------------------------------------------------------------------
\begin{equation*}
\nabla_Z \Omega(X,Y) = 2 \mu(Z) h(X,Y) + \mu(X) h(Y,Z) + \mu(Y) h(X,Z) ,
\end{equation*}
%--------------------------------------------------------------------
for all $X,Y,Z \in \X(M)$, and the result easily follows from \eqref{ZermeloFlagCurvatureUnit}.
\end{proof}
%--------------------------------------------------------------------

%--------------------------------------------------------------------
\begin{remark}
A particular case of projective vector field is when $W$ is affine on $(M,h)$ (i.e. projective with $\mu = 0$). Then,
%--------------------------------------------------------------------
\begin{align*}
K^F_v(w) =& K(u,w) + \frac{1}{4}\left( \Omega(\hat u, \hat u) \Omega(\hat w, \hat w) - \Omega(\hat u, \hat w)^2 + \frac{3}{4} \Omega(\hat u, \hat u)^2 \right)\\
&- \frac{1}{2}\Vert \Omega^{\sharp}(\hat u) \Vert_{h}^2 + \frac{1}{4}\Vert \Omega^{\sharp}(\hat w) \Vert_{h}^2 + \frac{1}{2}\Omega(\hat u, \nabla_{\hat u} W) - \Omega(\hat w, \nabla_{\hat w} W).
\end{align*}
%--------------------------------------------------------------------
\end{remark}
%--------------------------------------------------------------------

\appendix
%--------------------------------------------------------------------
\section{Curvature of stationary spacetimes as Lorentzian submersions}
\label{stationarysection}
%--------------------------------------------------------------------

Let $(M,g)$ be a Riemannian manifold and let $\omega$ be an arbitrary one-form on $M$. Consider on $\R\times M$ the Lorentzian metric
%--------------------------------------------------------------------
\begin{equation}
\label{StationaryMetric}
\hat{g}((\tau_1,v_1),(\tau_2,v_2)) = -(\tau_1-\omega(v_1))(\tau_2-\omega(v_2)) + g(v_1,v_2) ,
\end{equation}
%--------------------------------------------------------------------
which is the one determined by \eqref{Stationary} with the convention about sign changed. The pair $(\R\times M,\hat{g})$ is called a {\em stationary spacetime}, and the projection map
%--------------------------------------------------------------------
\begin{equation}
\label{submMap}
\hat{\pi}: (\R\times M, \hat{g}) \to (M,g)
\end{equation}
%--------------------------------------------------------------------
is a Lorentzian submersion. Indeed, the vertical subspace $\mathcal{V}$ (the tangent space to the fibers) is generated by $\partial_t = (1,0)$ and the horizontal subspace is
%--------------------------------------------------------------------
\begin{equation*}
\mathcal{H} = \{ (\omega(v),v) : v \in TM \} .
\end{equation*}
%--------------------------------------------------------------------
To check that it is a submersion, just observe that
%--------------------------------------------------------------------
\begin{equation*}
\hat{g}( (\omega(v_1),v_1) , (\omega(v_2),v_2) ) = g(v_1,v_2) = g(\d\hat{\pi}(\omega(v_1),v_1),\d\hat{\pi}(\omega(v_2),v_2)) ;
\end{equation*}
%--------------------------------------------------------------------
i.e. the restriction of $\d\hat{\pi}$ to horizontal vectors is an isometry. For each $X \in \X(M)$, define its \emph{horizontal lift} by
%--------------------------------------------------------------------
\begin{equation}
\label{lift}
\hor{X} = \left( \omega(X),\,X \right) \in \X(\R \times M) .
\end{equation}
%--------------------------------------------------------------------

%--------------------------------------------------------------------
\begin{lemma}
\label{lem:br}
For all $X,Y\in\mathfrak{X}(M)$,
%--------------------------------------------------------------------
\begin{subequations}
\begin{equation}
\label{br1}
[\partial_t,\hor{X}] = 0 ;
\end{equation}
\begin{equation}
\label{br2}
[\hor{X},\hor{Y}] = \d\omega(X,Y) \partial_t + \hor{[X,Y]} .
\end{equation}
\end{subequations}
%--------------------------------------------------------------------
\end{lemma}
%--------------------------------------------------------------------
%--------------------------------------------------------------------
\begin{proof}
\eqref{br1} follows from the fact that $\omega,X,Y$ are independent of the real coordinate $t$ for $\R\times M$. For the same reason, the $\partial_t$-component of $[\hor{X},\hor{Y}]$ is $X(\omega(Y)) - Y(\omega(X))$. Taking into account the push-forward identity for Lie brackets, e.g. \cite[Lemma~1.22]{one83s}, we have that
%--------------------------------------------------------------------
\begin{equation*}
\d\hat{\pi} ([\hor{X},\hor{Y}]) = [\d\hat{\pi}(X),\d\hat{\pi}(Y)] = [X,Y] ,
\end{equation*}
%--------------------------------------------------------------------
and it follows that
%--------------------------------------------------------------------
\begin{equation*}
[\hor{X},\hor{Y}] = (X(\omega(Y)) - Y(\omega(X)) , [X,Y]) .
\end{equation*}
%--------------------------------------------------------------------
Since $\d\omega(X,Y) = X(\omega(Y)) - Y(\omega(X)) - \omega([X,Y])$, we get \eqref{br2}; namely,
%--------------------------------------------------------------------
\begin{align*}
[\hor{X},\hor{Y}] &= (\d\omega(X,Y),0) + (\omega([X,Y]),[X,Y]) \\
&= \d\omega(X,Y)\partial_t + \hor{[X,Y]} .
\end{align*}
%--------------------------------------------------------------------
\end{proof}
%--------------------------------------------------------------------

Let $T$ and $A$ be the fundamental tensors for the submersion \eqref{submMap}; see \cite[\S~2]{one66f} or \cite[\S~2]{gray67p} for the definitions.

%--------------------------------------------------------------------
\begin{lemma}
\label{AXY}
For all $X,Y \in \X (M)$,
%--------------------------------------------------------------------
\begin{subequations}
\begin{equation}
T = 0 ;
\end{equation}
\begin{equation}
\label{Atensor}
A_{\hor{X}} \hor{Y} = \frac{1}{2} \d\omega(X,Y) \partial_t .
\end{equation}
\end{subequations}
%--------------------------------------------------------------------
\end{lemma}
%--------------------------------------------------------------------
%--------------------------------------------------------------------
\begin{proof}
Let $\hnab$ denote the Levi-Civita connection for $(\R\times M,\hat{g})$. By the properties of the tensor $T$, it suffices to prove that $T_{\partial_t} \partial_t = 0$. To verify this, we will actually show that $\hnab_{\partial_t} \partial_t$ vanishes. On one hand,
%--------------------------------------------------------------------
\begin{equation*}
\hat{g}(\hnab_{\partial_t} \partial_t,\partial_t) = \frac{1}{2} \partial_t (\hat{g}(\partial_t,\partial_t)) = 0 ,
\end{equation*}
%--------------------------------------------------------------------
because $\hat{g}(\partial_t,\partial_t) = -1$. So the vertical part of $\hnab_{\partial_t} \partial_t$ is zero. On the other hand, for any $X \in \X(M)$,
%--------------------------------------------------------------------
\begin{align*}
2\hat{g}(\hnab_{\partial_t} \partial_t,{\hor{X}}) =& 2\partial_t (\underbrace{\hat{g}(\partial_t,\hor{X})}_{0}) - \underbrace{\hor{X} (\hat{g}(\partial_t,\partial_t))}_{0} \\
&+ \hat{g}(\underbrace{[\partial_t,\partial_t]}_{0},\hor{X}) -2\hat{g}(\underbrace{[\partial_t,\hor{X}]}_{0},\partial_t) = 0 .
\end{align*}
%--------------------------------------------------------------------
Hence, $\hnab_{\partial_t} \partial_t$ is indeed zero.

For the tensor $A$ on horizontal vector fields, if $X,Y \in \X(M)$, then
%--------------------------------------------------------------------
\begin{align*}
2\hat{g}(\hnab_{\hor{X}} \hor{Y},\partial_t) =& \hor{X} (\underbrace{\hat{g}(\hor{Y},\partial_t)}_{0}) + \hor{Y} (\underbrace{\hat{g}(\hor{X},\partial_t)}_{0}) - \partial_t (\hat{g}(\hor{X},\hor{Y})) \\
&+ \hat{g}([\hor{X},\hor{Y}],\partial_t) - \hat{g}(\underbrace{[\hor{X},\partial_t]}_{0},\hor{Y}) - \hat{g}(\underbrace{[\hor{Y},\partial_t]}_{0},\hor{X}) \\
=& - \underbrace{\partial_t (g(X,Y))}_{0} + \hat{g}(\d\omega(X,Y)\partial_t + \hor{[X,Y]},\partial_t) \\
=& \d\omega(X,Y)\underbrace{\hat{g}(\partial_t,\partial_t)}_{-1} + \underbrace{\hat{g}(\hor{[X,Y]},\partial_t)}_{0} \\
=& - \d\omega(X,Y) .
\end{align*}
%--------------------------------------------------------------------
Therefore,
%--------------------------------------------------------------------
\begin{equation*}
A_{\hor{X}} \hor{Y} = \frac{\hat{g}(\hnab_{\hor{X}} \hor{Y},\partial_t)}{\hat{g}(\partial_t,\partial_t)}\partial_t = \frac{1}{2} \d\omega(X,Y) \partial_t .
\end{equation*}
%--------------------------------------------------------------------
\end{proof}
%--------------------------------------------------------------------

%--------------------------------------------------------------------
\begin{remark}
For any $X,Y \in \X(M)$, it is straightforward to compute that the horizontal part of $\hnab_{\hor{X}} \hor{Y}$ is exactly $\hor{\nabla_X Y}$. Together with \eqref{Atensor}, it implies
%--------------------------------------------------------------------
\begin{equation}
\label{nabHorVectors}
\hnab_{\hor{X}} \hor{Y} = \frac{1}{2} \d\omega(X,Y) \partial_t + \hor{\nabla_X Y} .
\end{equation}
%--------------------------------------------------------------------
\end{remark}
%--------------------------------------------------------------------

%--------------------------------------------------------------------
\begin{lemma}
\label{lemmaAxi}
For any $X \in \X(M)$,
%--------------------------------------------------------------------
\begin{equation}
\label{Axi}
\hnab_{\partial_t} \hor{X} = A_{\hor{X}} \partial_t = \frac{1}{2}\hor{\d\omega^{\sharp} (X)} ,
\end{equation}
%--------------------------------------------------------------------
where $\d\omega^{\sharp} (X)$ is determined by $g(\d\omega^{\sharp}(X),Y) = \d\omega(X,Y)$, $\forall Y \in \X(M)$.
\end{lemma}
%--------------------------------------------------------------------
%--------------------------------------------------------------------
\begin{proof}
First consider that $\hat{g}(\hnab_{\hor{X}} \partial_t,\partial_t) = \frac{1}{2} \hor{X} (\hat{g}(\partial_t,\partial_t)) = 0$, because $\hat{g}(\partial_t,\partial_t) = -1$. So $\hnab_{\hor{X}}\partial_t$ is horizontal and, hence, $\hnab_{\hor{X}} \partial_t = A_{\hor{X}} \partial_t$. By torsion-freeness and \eqref{br1}, we have $\hnab_{\partial_t} \hor{X} = \hnab_{\hor{X}} \partial_t$, and we conclude the first equality in \eqref{Axi}. Now, by the anti-symmetry of $A$, for any $Y \in \X(M)$,
%--------------------------------------------------------------------
\begin{align*}
\hat{g}(A_{\hor{X}} \partial_t,\hor{Y}) &= -\hat{g}(\partial_t,A_{\hor{X}} \hor{Y}) \overset{\eqref{Atensor}}{=} -\frac{1}{2} \hat{g}(\partial_t,\d\omega(X,Y)\partial_t) \\
&= \frac{1}{2}\d\omega(X,Y) = \frac{1}{2} g(\d\omega^{\sharp}(X),Y) = \frac{1}{2} \hat{g}(\hor{\d\omega^{\sharp}(X)},\hor{Y}) ,
\end{align*}
%--------------------------------------------------------------------
which implies the second equality in \eqref{Axi}.
\end{proof}
%--------------------------------------------------------------------

Recall the O'Neill formulas for the curvature tensor of a Lorentzian submersion; see \cite[Theorems~2 and~3]{one66f} or \cite[Theorem~3.3]{gray67p}, taking into account our opposite convention on curvature sign. To write them, let $\hat{R}$ be the curvature tensor of $(\R\times M)$, $R$ the curvature tensor of $(M,g)$, and denote by $\hor{R}$ the horizontal lift of $R$. For all $X,Y,Z,H \in \X(M)$,
%--------------------------------------------------------------------
\begin{subequations}
\begin{multline*}
\hat{g}(\hat{R}(\hor{X},\hor{Y}) \hor{Z}, \hor{H}) = \hat{g}(\hor{R}(\hor{X},\hor{Y}) \hor{Z}, \hor{H}) + 2\hat{g}(A_{\hor{X}} \hor{Y}, A_{\hor{Z}} \hor{H}) \\
- \hat{g}(A_{\hor{Y}} \hor{Z}, A_{\hor{X}} \hor{H}) - \hat{g}(A_{\hor{Z}} \hor{X}, A_{\hor{Y}} \hor{H}) ;
\end{multline*}
\begin{multline*}
\hat{g}(\hat{R}(\hor{X},\partial_t) \partial_t, \hor{Y}) = \hat{g}((\hnab_{\hor{X}} T)_{\partial_t} \partial_t, \hor{Y}) + \hat{g}((\hnab_{\partial_t} A)_{\hor{X}} \hor{Y}, \partial_t) \\
- \hat{g}(T_{\partial_t} \hor{X}, T_{\partial_t} \hor{Y}) + \hat{g}(A_{\hor{X}} \partial_t, A_{\hor{Y}} \partial_t) ;
\end{multline*}
\begin{multline*}
\hat{g}(\hat{R}(\hor{X},\hor{Y}) \hor{Z}, \partial_t) = \hat{g}(A_{\hor{Y}} \hor{Z}, T_{\partial_t} \hor{X}) + \hat{g}(A_{\hor{Z}} \hor{X}, T_{\partial_t} \hor{Y}) \\
- \hat{g}((\hnab_{\hor{Z}} A)_{\hor{X}} \hor{Y}, \partial_t) - \hat{g}(A_{\hor{X}} \hor{Y}, T_{\partial_t} \hor{Z}) .
\end{multline*}
\end{subequations}
%--------------------------------------------------------------------
In our case, $T=0$, and they simplify to:
%--------------------------------------------------------------------
\begin{subequations}
\begin{multline}
\label{firstFund}
\hat{g}(\hat{R}(\hor{X},\hor{Y}) \hor{Z}, \hor{H}) = \hat{g}(\hor{R}(\hor{X},\hor{Y}) \hor{Z}, \hor{H}) + 2\hat{g}(A_{\hor{X}} \hor{Y}, A_{\hor{Z}} \hor{H}) \\
- \hat{g}(A_{\hor{Y}} \hor{Z}, A_{\hor{X}} \hor{H}) - \hat{g}(A_{\hor{Z}} \hor{X}, A_{\hor{Y}} \hor{H}) ;
\end{multline}
\begin{equation}
\label{secondFund}
\hat{g}(\hat{R}(\hor{X},\partial_t) \partial_t, \hor{Y}) = \hat{g}((\hnab_{\partial_t} A)_{\hor{X}} \hor{Y}, \partial_t) + \hat{g}(A_{\hor{X}} \partial_t, A_{\hor{Y}} \partial_t) ;
\end{equation}
\begin{equation}
\label{thirdFund}
\hat{g}(\hat{R}(\hor{X},\hor{Y}) \hor{Z}, \partial_t) = - \hat{g}((\hnab_{\hor{Z}} A)_{\hor{X}} \hor{Y}, \partial_t) .
\end{equation}
\end{subequations}
%--------------------------------------------------------------------

%--------------------------------------------------------------------
\begin{lemma}
\label{lemanablaA}
Given the Lorentzian submersion \eqref{submMap}, the following holds for all $ X,Y,Z \in \X(M)$:
%--------------------------------------------------------------------
\begin{subequations}
\begin{equation}
\label{nablA2}
\hat{g}((\hnab_{\partial_t} A)_{\hor{X}} \hor{Y}, \partial_t) = 0 ,
\end{equation}
\begin{equation}
\label{nablaA3}
\hat{g}((\hnab_{\hor{Z}} A)_{\hor{X}} \hor{Y}, \partial_t) = -\frac{1}{2} (\nabla_Z \d\omega)(X,Y) ,
\end{equation}
\end{subequations}
%--------------------------------------------------------------------
where $\hnab, \nabla$ are the Levi-Civita connections of $(\R\times M, \hat{g}), (M,g)$, respectively.
\end{lemma}
%--------------------------------------------------------------------
%--------------------------------------------------------------------
\begin{proof}
By Lemmas~\ref{AXY} and~\ref{lemmaAxi}, including the facts that $\hnab_{\partial_t} \partial_t = 0$ (from the proof of Lemma~\ref{AXY}) and $A$ is anti-symmetric, we compute \eqref{nablA2}:
%--------------------------------------------------------------------
\begin{align*}
(\hnab_{\partial_t} A)_{\hor{X}} \hor{Y} &= \hnab_{\partial_t} (A_{\hor{X}} \hor{Y}) - A_{(\hnab_{\partial_t} \hor{X})} \hor{Y} - A_{\hor{X}} (\hnab_{\partial_t} \hor{Y}) \\
&= \underbrace{\hnab_{\partial_t} \left( \frac{1}{2} \d\omega(X,Y) \partial_t \right)}_{0} - \frac{1}{2} A_{\hor{\d\omega^{\sharp}(X)}} \hor{Y} - \frac{1}{2} A_{\hor{X}} (\hor{\d\omega^{\sharp}(Y)}) \\
&= \frac{1}{2} \left( A_{\hor{Y}} (\hor{\d\omega^{\sharp}(X)}) - A_{\hor{X}} (\hor{\d\omega^{\sharp}(Y)}) \right) \\
&= \frac{1}{4} \left( \d\omega(Y,\d\omega^{\sharp}(X)) - \d\omega(X,\d\omega^{\sharp}(Y)) \right) \partial_t \\
&= \frac{1}{4} \left( g(\d\omega^{\sharp}(Y),\d\omega^{\sharp}(X)) - g(\d\omega^{\sharp}(X),\d\omega^{\sharp}(Y))\right) \partial_t = 0 .
\end{align*}
%--------------------------------------------------------------------
Lemmas~\ref{AXY} and~\ref{lemmaAxi}, together with the identities $\hnab_{\hor{Z}} \hor{X} = \hor{\nabla_Z X} + A_{\hor{Z}} \hor{X}$ and $A_{\hnab_{\hor{Z}} \hor{X}} = A_{\hor{\nabla_{Z} X}}$, also imply that
%--------------------------------------------------------------------
\begin{align*}
(\hnab_{\hor{Z}} A)_{\hor{X}} \hor{Y} =& \hnab_{\hor{Z}} (A_{\hor{X}} \hor{Y}) - A_{(\hnab_{\hor{Z}} \hor{X})} \hor{Y} - A_{\hor{X}} (\hnab_{\hor{Z}} \hor{Y}) \\
=& \hnab_{\hor{Z}} \left( \frac{1}{2} \d\omega(X,Y) \partial_t \right) - A_{\hor{\nabla_{Z} X}} \hor{Y} - A_{\hor{X}} (\hor{\nabla_Z Y} + A_{\hor{Z}} \hor{Y}) \\
=& \frac{1}{2} \hor{Z} (\d\omega (X,Y)) \partial_t + \frac{1}{2} \d\omega(X,Y) \hnab_{\hor{Z}} \partial_t \\
&- \frac{1}{2} \d\omega(\nabla_Z X,Y)\partial_t - \frac{1}{2} \d\omega (X,\nabla_Z Y) \partial_t - A_{\hor{X}} A_{\hor{Z}} \hor{Y} \\
=& \frac{1}{2} \left( Z(\d\omega (X,Y)) - \d\omega(\nabla_Z X,Y) - \d\omega (X,\nabla_Z Y) \right) \partial_t \\
&+ \frac{1}{4} \d\omega(X,Y) \hor{\d\omega^{\sharp}(Z)} - \frac{1}{4} \d\omega(Z,Y) \hor{\d\omega^{\sharp}(X)} ,
\end{align*}
%--------------------------------------------------------------------
where $A_{\hor{X}} A_{\hor{Z}} \hor{Y}$ is computed as the horizontal part of $\hnab_{\hor{X}} (A_{\hor{Z}} \hor{Y}) = \frac{1}{2} \hnab_{\hor{X}} (\d\omega (Z,Y) \partial_t)$, which is just $A_{\hor{X}} A_{\hor{Z}} \hor{Y} = \frac{1}{2} \d\omega (Z,Y) \hnab_{\hor{X}} \partial_t$. Thus, we deduce the formula \eqref{nablaA3}.
\end{proof}
%--------------------------------------------------------------------

%--------------------------------------------------------------------
\begin{proposition}
\label{CurvTen}
The curvature tensor of a stationary spacetime as in \eqref{StationaryMetric} is given by
%--------------------------------------------------------------------
\begin{multline}
\label{submCurv}
\hat{g}(\hat{R}(\mu\partial_t + \hor{X}, \hor{Y}) \hor{Y}, \mu\partial_t + \hor{X}) = g(R(X,Y)Y,X) + \frac{3}{4} (\d\omega(X,Y))^2 \\
+ \mu (\nabla_Y \d\omega)(X,Y) + \frac{\mu^2}{4} \Vert \d\omega^{\sharp} (Y) \Vert^2 ,
\end{multline}
%--------------------------------------------------------------------
for arbitrary $X,Y \in \X(M)$ and $\mu\in \R$, where $\Vert \cdot \Vert$ is the norm associated with the Riemannian metric $g$.
\end{proposition}
%--------------------------------------------------------------------
%--------------------------------------------------------------------
\begin{proof}
First, using the multilinearity and symmetries of the curvature tensor, it follows that
%--------------------------------------------------------------------
\begin{multline*}
\hat{g}(\hat{R}(\mu\partial_t + \hor{X}, \hor{Y}) \hor{Y}, \mu\partial_t + \hor{X}) = \\
= \hat{g}(\hat{R}(\hor{X},\hor{Y})\hor{Y},\hor{X}) + 2\mu \hat{g}(\hat{R}(\hor{X},\hor{Y}) \hor{Y}, \partial_t) +
\mu^2 \hat{g}(\hat{R}(\hor{Y},\partial_t) \partial_t, \hor{Y}) .
\end{multline*}
%--------------------------------------------------------------------
Moreover, using O'Neill formulas for the Lorentzian submersion \eqref{submMap} as in \eqref{firstFund}--\eqref{thirdFund} and Lemmas~\ref{AXY}--\ref{lemanablaA}, we get
%--------------------------------------------------------------------
\begin{subequations}
\begin{align*}
\hat{g}(\hat{R}(\hor{H},\hor{Z}) \hor{Z}, \hor{H}) =& \hat{g}(\hor{R}(\hor{H},\hor{Z}) \hor{Z}, \hor{H}) + 2\hat{g}(A_{\hor{H}} \hor{Z}, A_{\hor{Z}} \hor{H}) \\
&- \hat{g}(\underbrace{A_{\hor{Z}} \hor{Z}}_{0}, \underbrace{A_{\hor{H}} \hor{H}}_{0}) - \hat{g}(A_{\hor{Z}} \hor{H}, A_{\hor{Z}} \hor{H}) \\
=& \hat{g}(\hor{R}(\hor{H},\hor{Z}) \hor{Z}, \hor{H}) - 3\hat{g}(A_{\hor{H}} \hor{Z}, A_{\hor{H}} \hor{Z}) \\
=& g(R(H,Z)Z,H) + \frac{3}{4} (\d\omega(H,Z))^2,
\end{align*}
\begin{equation*}
\hat{g}(\hat{R}(\hor{H},\hor{Z}) \hor{Z}, \partial_t) = - \hat{g}((\hnab_{\hor{Z}} A)_{\hor{H}} \hor{Z}, \partial_t) = \frac{1}{2} (\nabla_Z \d\omega) (H,Z) ,
\end{equation*}
\begin{align*}
\hat{g}(\hat{R}(\hor{H},\partial_t) \partial_t, \hor{Z}) &= \underbrace{\hat{g}((\hnab_{\partial_t} A)_{\hor{H}} \hor{Z}, \partial_t)}_{0} + \hat{g}(A_{\hor{H}} \partial_t, A_{\hor{Z}} \partial_t) \\
&= \frac{1}{4} \hat{g}(\hor{\d\omega^{\sharp} (H)}, \hor{\d\omega^{\sharp} (Z)}) = \frac{1}{4} g(\d\omega^{\sharp} (H), \d\omega^{\sharp} (Z)) ,
\end{align*}
\end{subequations}
%--------------------------------------------------------------------
for any $Z,H\in\mathfrak{X}(M)$. Substituting these formulas into the first identity of the proof, we obtain the desired equation for the curvature tensor.
\end{proof}
%--------------------------------------------------------------------
%--------------------------------------------------------------------
\section{Curvature of the stationary spacetime in standard form}
\label{stdstationarysection}
%--------------------------------------------------------------------

The {\em standard stationary spacetime} consists of a product manifold $\R\times M$ endowed with the Lorentzian metric $\til{g}$ as in \eqref{stdstationary}, i.e.
%--------------------------------------------------------------------
\begin{equation*}
\til{g}((\tau_1,v_1),(\tau_2,v_2)) = g_0(v_1,v_2) + g_0(v_1,\delta)\tau_2 + g_0(v_2,\delta)\tau_1 - \beta\tau_1\tau_2 ,
\end{equation*}
%--------------------------------------------------------------------
where $(M,g_0)$ is a Riemannian manifold, $\delta \in \X(M)$, and $\beta: M \to \R$ is a positive smooth function. We denote by $\tnab, \nabla$ the Levi-Civita connections of $\til{g}, g_0$, respectively, and by $\nabla\beta$ the gradient of $\beta$ with respect to $g_0$. We will also denote with overtilde the lift of $X \in \X(M)$ to $\R\times M$, namely, $\til{X} = (0,X) \in \X(\R\times M)$. Furthermore, we define the auxiliary function
%--------------------------------------------------------------------
\begin{equation}
\label{Lambda}
\Lambda = - \frac{1}{\beta + g_0(\delta,\delta)} ,
\end{equation}
%--------------------------------------------------------------------
and, for $X,Y \in \X(M)$, the operators
%--------------------------------------------------------------------
\begin{equation}
\label{Omega}
\begin{aligned}
\Omega_{+} (X,Y) &= g_0(\nabla_X \delta, Y) + g_0(\nabla_Y \delta, X) = g_0(\Omega_{+}^{\sharp}(X), Y) , \\
\Omega_{-} (X,Y) &= g_0(\nabla_X \delta, Y) - g_0(\nabla_Y \delta, X) = g_0(\Omega_{-}^{\sharp}(X), Y) .
\end{aligned}
\end{equation}
%--------------------------------------------------------------------

Let $\til{R}$ be the Riemann curvature tensor for $\til{g}$. Our goal is to compute $\til{g}(\til{R}(\mu\partial_t + \til{X}, \til{Y}) \til{Y}, \mu\partial_t + \til{X})$, for $X,Y \in \X(M)$, $\mu \in \R$, and $\partial_t = (1,0)$. We begin with the following identities:
%--------------------------------------------------------------------
\begin{subequations}
\begin{equation}
\label{firstlev}
\tnab_{\partial_t} \partial_t = \frac{1}{2} \left( -\Lambda \delta(\beta) (\partial_t - \til{\delta}) + \til{\nabla\beta} \right) ,
\end{equation}
\begin{equation}
\label{secondlev}
\tnab_{\til{X}} \til{Y} = \frac{\Lambda}{2} \Omega_{+}(X,Y) (\partial_t - \til{\delta}) + \til{\nabla_X Y} ,
\end{equation}
\begin{equation}
\label{thirdlev}
\tnab_{\til{X}} \partial_t = \tnab_{\partial_t} \til{X} = -\frac{\Lambda}{2} \left( X(\beta) + \Omega_{-}(X,\delta) \right) (\partial_t - \til{\delta}) + \frac{1}{2}\til{\Omega_{-}^{\sharp}(X)} ;
\end{equation}
\end{subequations}
%--------------------------------------------------------------------
see \cite[Theorem~1]{flo03g}.

%--------------------------------------------------------------------
\begin{lemma}
For all $X \in \X(M)$,
%--------------------------------------------------------------------
\begin{subequations}
\begin{equation}
\label{fourthlev}
\tnab_{\partial_t} (\partial_t - \til{\delta}) = \frac{1}{2} (\til{\nabla\beta} - \til{\Omega_{-}^{\sharp}(\delta)}) ;
\end{equation}
\begin{equation}
\label{fifthlev}
\tnab_{\til{X}} (\partial_t - \til{\delta}) = - \frac{X(\Lambda)}{2\Lambda} (\partial_t - \til{\delta}) - \frac{1}{2} \til{\Omega_{+}^{\sharp}(X)} .
\end{equation}
\end{subequations}
\end{lemma}
%--------------------------------------------------------------------
%--------------------------------------------------------------------
\begin{proof}
Equation~\eqref{fourthlev} follows directly from \eqref{firstlev} and \eqref{thirdlev}, as $\Omega_{-}(X,X) = 0$, $\forall X\in \X(M)$. Similarly, \eqref{fifthlev} follows from \eqref{secondlev} and \eqref{thirdlev}, using the identities $\d\Lambda(X) = \Lambda^2 (X(\beta) + 2g_0(\nabla_X \delta, \delta))$ and $\Omega_{+}^{\sharp}(X) + \Omega_{-}^{\sharp}(X) = 2\nabla_X \delta$, $\forall X \in \X(M)$.
\end{proof}
%--------------------------------------------------------------------

%--------------------------------------------------------------------
\begin{lemma}
The Riemann curvature tensor $\til{R}$ for the standard stationary Lorentzian metric $\til{g}$ as in \eqref{stdstationary} satisfies, for all $X,Y \in \X(M)$,
%--------------------------------------------------------------------
\begin{subequations}
\begin{multline}
\label{curvXYY}
\til{R}(\til{X},\til{Y}) \til{Y} = \left( \frac{1}{4} \big( X(\Lambda) \Omega_{+}(Y,Y) - Y(\Lambda) \Omega_{+}(X,Y) \big) \right. \\
\left. + \frac{\Lambda}{2} \big( (\nabla_X \Omega_{+})(Y,Y) - (\nabla_Y \Omega_{+})(X,Y) \big) \right) (\partial_t - \til{\delta}) \\
- \frac{\Lambda}{4} \left( \Omega_{+}(Y,Y) \til{\Omega_{+}^{\sharp} (X)} - \Omega_{+}(X,Y) \til{\Omega_{+}^{\sharp} (Y)} \right) + \til{R(X,Y)Y} ,
\end{multline}
\begin{multline}
\label{curvtYY}
\til{R}(\partial_t,\til{Y}) \til{Y} = \left( \frac{Y(\Lambda)}{4} \big( Y(\beta) + \Omega_{-}(Y,\delta) \big) - \frac{\Lambda}{4} \Omega_{+}(Y,\Omega_{-}^{\sharp}(Y)) \right. \\
\left. + \frac{\Lambda}{2} \big( (\nabla_Y \d\beta)(Y) + (\nabla_Y\Omega_{-})(Y,\delta) + \Omega_{-}(Y,\nabla_Y \delta) \big) \right) (\partial_t - \til{\delta}) \\
+ \frac{\Lambda}{4} \left( \Omega_{+}(Y,Y) \big( \til{\nabla\beta} - \til{\Omega_{-}^{\sharp}(\delta)} \big) - \big( Y(\beta) + \Omega_{-}(Y,\delta) \big) \til{\Omega_{+}^{\sharp}(Y)} \right) \\
- \frac{1}{2} \til{(\nabla_Y \Omega_{-}^{\sharp}) (Y)} ,
\end{multline}
\end{subequations}
%--------------------------------------------------------------------
where $\nabla$ and $R$ denote, respectively, the Levi-Civita connection and the Riemann curvature tensor for $(M,g_0)$.
\end{lemma}
%--------------------------------------------------------------------
%--------------------------------------------------------------------
\begin{proof}
Let us choose $X,Y \in \X(M)$ such that, at a certain point $p\in M$, they satisfy $[X,Y]\vert_p = \nabla_X Y\vert_p = \nabla_Y X\vert_p = \nabla_X X\vert_p = \nabla_Y Y\vert_p = 0$, which is not restrictive. Then, using \eqref{secondlev},
%--------------------------------------------------------------------
\begin{align*}
\tnab_{\til{X}} (\til{\nabla_Y Y})\vert_{(t,p)} &= \til{\nabla_X \nabla_Y Y}\vert_{(t,p)} , \\
\tnab_{\til{Y}} (\til{\nabla_X Y})\vert_{(t,p)} &= \til{\nabla_Y \nabla_X Y}\vert_{(t,p)} .
\end{align*}
%--------------------------------------------------------------------
It now follows from \eqref{secondlev}, \eqref{fifthlev}, and the above identities that
%--------------------------------------------------------------------
\begin{subequations}
\begin{align*}
\tnab_{\til{X}} \tnab_{\til{Y}} \til{Y}\vert_{(t,p)} =& \frac{1}{2} \big( X(\Lambda) \Omega_{+}(Y,Y) + \Lambda (\nabla_X \Omega_{+})(Y,Y) \big) (\partial_t - \til{\delta})\vert_{(t,p)} \\
&+ \frac{\Lambda}{2} \Omega_{+}(Y,Y) \tnab_{\til{X}} (\partial_t - \til{\delta})\vert_{(t,p)} + \tnab_{\til{X}} (\til{\nabla_Y Y})\vert_{(t,p)} \\
=& \left( \frac{X(\Lambda)}{4} \Omega_{+}(Y,Y) + \frac{\Lambda}{2} (\nabla_X \Omega_{+})(Y,Y) \right) (\partial_t - \til{\delta})\vert_{(t,p)} \\
&- \frac{\Lambda}{4} \Omega_{+}(Y,Y) \til{\Omega_{+}^{\sharp}(X)}\vert_{(t,p)} + \til{\nabla_X \nabla_Y Y}\vert_{(t,p)} ,
\end{align*}
\begin{align*}
\tnab_{\til{Y}} \tnab_{\til{X}} \til{Y}\vert_{(t,p)} =& \frac{1}{2} \big( Y(\Lambda) \Omega_{+}(X,Y) + \Lambda (\nabla_Y \Omega_{+})(X,Y) \big) (\partial_t - \til{\delta})\vert_{(t,p)} \\
&+ \frac{\Lambda}{2} \Omega_{+}(X,Y) \tnab_{\til{Y}} (\partial_t - \til{\delta})\vert_{(t,p)} + \tnab_{\til{Y}} (\til{\nabla_X Y})\vert_{(t,p)} \\
=& \left( \frac{Y(\Lambda)}{4} \Omega_{+}(X,Y) + \frac{\Lambda}{2} (\nabla_Y \Omega_{+})(X,Y) \right) (\partial_t - \til{\delta})\vert_{(t,p)} \\
&- \frac{\Lambda}{4} \Omega_{+}(X,Y) \til{\Omega_{+}^{\sharp}(Y)}\vert_{(t,p)} + \til{\nabla_Y \nabla_X Y}\vert_{(t,p)} .
\end{align*}
\end{subequations}
%--------------------------------------------------------------------
Taking into account that $[\til{X},\til{Y}] = \til{[X,Y]}$, and using \eqref{secondlev} again, we have that
%--------------------------------------------------------------------
\begin{equation*}
\tnab_{[\til{X},\til{Y}]} \til{Y}\vert_{(t,p)} = \tnab_{\til{[X,Y]}} \til{Y}\vert_{(t,p)} = 0 .
\end{equation*}
%--------------------------------------------------------------------
Substituting the above equalities into the definition
%--------------------------------------------------------------------
\begin{equation*}
\til{R}(\til{X},\til{Y}) \til{Y} = \tnab_{\til{X}} \tnab_{\til{Y}} \til{Y} - \tnab_{\til{Y}} \tnab_{\til{X}} \til{Y} - \tnab_{[\til{X},\til{Y}]} \til{Y} ,
\end{equation*}
%--------------------------------------------------------------------
we conclude \eqref{curvXYY} at $(t,p)$, which is an arbitrary point.

The proof of \eqref{curvtYY} is going to be similar. Under the same assumptions on $X,Y \in \X(M)$, we first obtain from \eqref{thirdlev} that
%--------------------------------------------------------------------
\begin{equation*}
\tnab_{\partial_t} \til{\nabla_Y Y}\vert_{(t,p)} = 0 .
\end{equation*}
%--------------------------------------------------------------------
Then, using \eqref{firstlev}, \eqref{fourthlev}, and the fact that $\Lambda$ and $\Omega_{+}$ have no dependence on $t$,
%--------------------------------------------------------------------
\begin{align*}
\tnab_{\partial_t} \tnab_{\til{Y}} \til{Y}\vert_{(t,p)} &= \frac{\Lambda}{2} \Omega_{+}(Y,Y) \tnab_{\partial_t} (\partial_t - \til{\delta})\vert_{(t,p)} + \tnab_{\partial_t} \til{\nabla_Y Y}\vert_{(t,p)} \\
&= \frac{\Lambda}{4} \Omega_{+}(Y,Y) \big( \til{\nabla\beta} - \til{\Omega_{-}^{\sharp}(\delta)} \big)\vert_{(t,p)} .
\end{align*}
%--------------------------------------------------------------------
Next, \eqref{secondlev}, \eqref{thirdlev}, \eqref{fifthlev}, and $\beta, \delta, \Lambda, \Omega_{-}$ being independent of $t$ imply that
%--------------------------------------------------------------------
\begin{align*}
\tnab_{\til{Y}} \tnab_{\partial_t} &\til{Y}\vert_{(t,p)} = \\
=& -\frac{Y(\Lambda)}{2} \big( Y(\beta) + \Omega_{-}(Y,\delta) \big) (\partial_t - \til{\delta})\vert_{(t,p)} \\
&- \frac{\Lambda}{2} \big( (\nabla_Y \d\beta)(Y) + (\nabla_Y \Omega_{-})(Y,\delta) + \Omega_{-}(Y,\nabla_Y \delta) \big) (\partial_t - \til{\delta})\vert_{(t,p)} \\
&- \frac{\Lambda}{2} (Y(\beta) + \Omega_{-}(Y,\delta)) \left.\left( -\frac{Y(\Lambda)}{2\Lambda} (\partial_t - \til{\delta}) - \frac{1}{2}\til{\Omega_{+}^{\sharp}(Y)} \right)\right\vert_{(t,p)} \\
&+ \frac{\Lambda}{4}\Omega_{+}(Y,\Omega_{-}^{\sharp}(Y)) (\partial_t - \til{\delta})\vert_{(t,p)} + \frac{1}{2} \til{\nabla_Y (\Omega_{-}^{\sharp}(Y))}\vert_{(t,p)} \\
=& \left( -\frac{Y(\Lambda)}{4} (Y(\beta) + \Omega_{-}(Y,\delta)) + \frac{\Lambda}{4} \Omega_{+}(Y,\Omega_{-}^{\sharp}(Y))
\right. \\
&\left. -\frac{\Lambda}{2} \big( (\nabla_Y \d\beta)(Y) + (\nabla_Y \Omega_{-})(Y,\delta) + \Omega_{-}(Y,\nabla_Y \delta) \big) \right) (\partial_t - \til{\delta})\vert_{(t,p)} \\
&+ \frac{\Lambda}{4} \big( Y(\beta) + \Omega_{-}(Y,\delta) \big) \til{\Omega_{+}^{\sharp}(Y)}\vert_{(t,p)} + \frac{1}{2} \til{(\nabla_Y \Omega_{-}^{\sharp})(Y)}\vert_{(t,p)} .
\end{align*}
%--------------------------------------------------------------------
Finally, $\tnab_{[\partial_t,\til{Y}]} \til{Y} = 0$, because $[\partial_t,\til{Y}] = 0$. Putting all this together in the definition of the curvature tensor yields \eqref{curvtYY} at $(t,p)$, as required.
\end{proof}
%--------------------------------------------------------------------

Now, take into account that
%--------------------------------------------------------------------
\begin{equation}
\label{gg0}
\til{g}(\til{X},\til{Y}) = g_0(X,Y) , \; \forall X,Y \in \X(M) ,
\end{equation}
%--------------------------------------------------------------------
and that
%--------------------------------------------------------------------
\begin{subequations}
\begin{equation}
\label{ortvNorm}
\til{g} (\partial_t - \til{\delta}, \partial_t - \til{\delta}) = \frac{1}{\Lambda} ;
\end{equation}
\begin{equation}
\label{ortv}
\til{g} (\partial_t - \til{\delta}, \til{X}) = 0 , \; \forall X \in \X(M) .
\end{equation}
\end{subequations}
%--------------------------------------------------------------------

%--------------------------------------------------------------------
\begin{proposition}
Let $\til{g}$ be the stationary Lorentzian metric defined in \eqref{stdstationary}. Then its curvature tensor $\til{R}$ satisfies, for any $X,Y \in \X(M)$ and $\mu\in\R$,
%--------------------------------------------------------------------
\begin{multline}
\label{finalresult}
\til{g}(\til{R}(\mu\partial_t + \til{X}, \til{Y}) \til{Y}, \mu\partial_t + \til{X}) = \\
= \frac{\Lambda}{4} \Omega_{+}(Y,Y) \big( \mu ((X + \mu\delta)(\beta) + \Omega_{-}(X,\delta)) - \Omega_{+}(X,X + \mu\delta) \big) \\
- \frac{\Lambda}{4} \Omega_{+}(Y,X + \mu\delta) \big( \mu (Y(\beta) + \Omega_{-}(Y,\delta)) - \Omega_{+}(X,Y) \big) \\
+ \frac{\mu}{4\Lambda} \Big( X(\Lambda)\Omega_{+}(Y,Y) + Y(\Lambda) \big( \mu (Y(\beta) + \Omega_{-}(Y,\delta)) - \Omega_{+}(X,Y) \big) \Big) \\
+ \frac{\mu}{2} \big( (\nabla_X \Omega_{+})(Y,Y) - (\nabla_Y \Omega_{+})(X,Y) \big) \\
+ \frac{\mu}{2} \big( (\nabla_{X + \mu\delta} \Omega_{+})(Y,Y) - (\nabla_Y \Omega_{+})(X + \mu\delta,Y) \big) \\
+ \frac{\mu^2}{2} \big( (\nabla_Y \d\beta)(Y) + (\nabla_Y \Omega_{-})(Y,\delta) + \Omega_{-}(Y,\nabla_Y \delta) \big) \\
- \frac{\mu^2}{4} \Omega_{+}(Y,\Omega_{-}^{\sharp}(Y)) + g_0(R(X + \mu\delta, Y) Y, X + \mu\delta) ,
\end{multline}
%--------------------------------------------------------------------
where $\nabla$ and $R$ are the Levi-Civita connection and the curvature tensor of $g_0$, respectively.
\end{proposition}
%--------------------------------------------------------------------
%--------------------------------------------------------------------
\begin{proof}
By the multilinearity of the curvature tensor $\til{R}$,
%--------------------------------------------------------------------
\begin{equation}
\begin{aligned}
\label{finalcount}
&\til{g}(\til{R}(\mu\partial_t + \til{X}, \til{Y}) \til{Y}, \mu\partial_t + \til{X}) = \\
&= \mu \til{g}(\til{R}(\mu\partial_t + \til{X}, \til{Y}) \til{Y}, \partial_t - \til{\delta}) + \til{g}(\til{R}(\mu\partial_t + \til{X}, \til{Y}) \til{Y}, \til{X} + \mu\til{\delta}) .
\end{aligned}
\end{equation}
%--------------------------------------------------------------------

From \eqref{curvXYY}, \eqref{ortvNorm} and \eqref{ortv}, it follows that
%--------------------------------------------------------------------
\begin{multline}
\label{partial_tdelta}
\til{g}(\til{R}(\til{X},\til{Y}) \til{Y}, \partial_t - \til{\delta}) = \frac{1}{4\Lambda} \big( X(\Lambda) \Omega_{+}(Y,Y) - Y(\Lambda) \Omega_{+}(X,Y) \big)\\
+ \frac{1}{2} \big( (\nabla_X \Omega_{+})(Y,Y) - (\nabla_Y \Omega_{+})(X,Y) \big) ,
\end{multline}
%--------------------------------------------------------------------
and, using \eqref{gg0},
%--------------------------------------------------------------------
\begin{multline}
\label{tildeZ}
\til{g}(\til{R}(\til{X},\til{Y}) \til{Y}, \til{Z}) = - \frac{\Lambda}{4} \big( \Omega_{+}(Y,Y) \Omega_{+}(X,Z) - \Omega_{+}(X,Y) \Omega_{+}(Y,Z) \big) \\
+ g_0(R(X,Y)Y, Z) ,
\end{multline}
%--------------------------------------------------------------------
for any $Z \in \X(M)$. Similarly, by \eqref{curvtYY}, \eqref{ortvNorm} and \eqref{ortv}, we obtain
%--------------------------------------------------------------------
\begin{multline}
\label{gdelta}
\tilde g(\tilde R(\partial_t,\til{Y})\til{Y},\partial_t-\til{\delta}) = \frac{Y(\Lambda)}{4\Lambda} \big( Y(\beta) + \Omega_{-}(Y,\delta) \big) - \frac{1}{4} \Omega_{+}(Y,\Omega_{-}^{\sharp}(Y)) \\
+ \frac{1}{2} \big( (\nabla_Y \d\beta)(Y) + (\nabla_Y\Omega_{-})(Y,\delta) + \Omega_{-}(Y,\nabla_Y \delta) \big) ,
\end{multline}
%--------------------------------------------------------------------
and for any $Z \in\X(M)$,
%--------------------------------------------------------------------
\begin{multline}
\label{tildeZ2}
\tilde g(\tilde R(\partial_t,\til{Y})\til{ Y},\til{Z}) = \frac{\Lambda}{4} \Omega_{+}(Y,Y) \big( Z(\beta) + \Omega_{-}(Z,\delta) \big) \\
- \frac{\Lambda}{4} \Omega_{+}(Y,Z) \big( Y(\beta) + \Omega_{-}(Y,\delta) \big) - \frac{1}{2} (\nabla_Y \Omega_{-})(Y,Z) ,
\end{multline}
%--------------------------------------------------------------------
because of \eqref{gg0}, the antisymmetry of the operator $\Omega_{-}$, and the fact that $g_0((\nabla_X \Omega_{-}^{\sharp}) (Y), Z) = (\nabla_X \Omega_{-})(Y,Z)$, $\forall X,Y,Z \in \X(M)$.

Then, putting together \eqref{partial_tdelta} and \eqref{gdelta}, we compute 
%--------------------------------------------------------------------
\begin{multline*}
\til{g}(\til{R}(\mu\partial_t + \til{X}, \til{Y}) \til{Y}, \partial_t - \til{\delta}) = \\
= \frac{1}{4\Lambda} \Big( X(\Lambda) \Omega_{+}(Y,Y) + Y(\Lambda) \big( \mu (Y(\beta) + \Omega_{-}(Y,\delta)) - \Omega_{+}(X,Y) \big) \Big) \\
+ \frac{1}{2} \big( (\nabla_X \Omega_{+})(Y,Y) - (\nabla_Y \Omega_{+})(X,Y) \big) \\
+ \frac{\mu}{2} \big( (\nabla_Y \d\beta)(Y) + (\nabla_Y \Omega_{-})(Y,\delta) + \Omega_{-}(Y,\nabla_Y \delta) \big) - \frac{\mu}{4} \Omega_{+}(Y,\Omega_{-}^{\sharp}(Y)) .
\end{multline*}
%--------------------------------------------------------------------
Likewise, combining \eqref{tildeZ} and \eqref{tildeZ2} for $Z = X + \mu\delta$ gives
%--------------------------------------------------------------------
\begin{multline*}
\til{g}(\til{R}(\mu\partial_t + \til{X}, \til{Y}) \til{Y}, \til{X} + \mu\til{\delta}) = \\
\frac{\Lambda}{4} \Omega_{+}(Y,Y) \Big( \mu \big( (X + \mu\delta)(\beta) + \Omega_{-}(X,\delta) \big) - \Omega_{+}(X,X + \mu\delta) \Big) \\
- \frac{\Lambda}{4} \Omega_{+}(Y,X + \mu\delta) \Big( \mu \big( Y(\beta) + \Omega_{-}(Y,\delta) \big) - \Omega_{+}(X,Y) \Big) \\
- \frac{\mu}{2} (\nabla_Y \Omega_{-})(Y,X + \mu\delta) + g_0(R(X,Y) Y, X + \mu\delta) ,
\end{multline*}
%--------------------------------------------------------------------
as $\Omega_{-}(X + \mu\delta,\delta) = \Omega_{-}(X,\delta)$, by bilinearity and antisymmetry of $\Omega_{-}$. Substituting the last two identities into \eqref{finalcount}, we get to
%--------------------------------------------------------------------
\begin{multline}
\label{prelimfinalresult}
\til{g}(\til{R}(\mu\partial_t + \til{X}, \til{Y}) \til{Y}, \mu\partial_t + \til{X}) = \\
= \frac{\Lambda}{4} \Omega_{+}(Y,Y) \big( \mu ((X + \mu\delta)(\beta) + \Omega_{-}(X,\delta)) - \Omega_{+}(X,X + \mu\delta) \big) \\
- \frac{\Lambda}{4} \Omega_{+}(Y,X + \mu\delta) \big( \mu (Y(\beta) + \Omega_{-}(Y,\delta)) - \Omega_{+}(X,Y) \big) \\
+ \frac{\mu}{4\Lambda} \Big( X(\Lambda)\Omega_{+}(Y,Y) + Y(\Lambda) \big( \mu (Y(\beta) + \Omega_{-}(Y,\delta)) - \Omega_{+}(X,Y) \big) \Big) \\
+ \frac{\mu}{2} \big( (\nabla_X \Omega_{+})(Y,Y) - (\nabla_Y \Omega_{+})(X,Y) - (\nabla_Y \Omega_{-})(Y,X + \mu\delta) \big) \\
+ \frac{\mu^2}{2} \big( (\nabla_Y \d\beta)(Y) + (\nabla_Y \Omega_{-})(Y,\delta) + \Omega_{-}(Y,\nabla_Y \delta) \big) \\
- \frac{\mu^2}{4} \Omega_{+}(Y,\Omega_{-}^{\sharp}(Y)) + g_0(R(X,Y) Y, X + \mu\delta) .
\end{multline}
%--------------------------------------------------------------------

To conclude \eqref{finalresult}, assume (without loss of generality) that $Z,Y \in \X(M)$ satisfy $[Z,Y]\vert_p = \nabla_Z Y\vert_p = \nabla_Y Z\vert_p = \nabla_Y Y\vert_p = 0$, at a certain point $p\in M$. As $\Omega_{+} (Z,Y) + \Omega_{-} (Z,Y) = 2 g_0(\nabla_Z \delta, Y)$, $\forall Z,Y \in \X(M)$, we compute
%--------------------------------------------------------------------
\begin{align*}
(\nabla_Y &\Omega_{-})(Z,Y)\vert_p = Y(\Omega_{-}(Z,Y))\vert_p \\
&= Y( 2 g_0(\nabla_Z \delta, Y) - \Omega_{+} (Z,Y))\vert_p \\
&= 2 g_0(\nabla_Y \nabla_Z \delta, Y)\vert_p - (\nabla_Y \Omega_{+})(Z,Y)\vert_p \\
&= 2 g_0(R(Y,Z)\delta, Y)\vert_p + 2 g_0(\nabla_Z \nabla_Y \delta, Y)\vert_p - (\nabla_Y \Omega_{+})(Z,Y)\vert_p \\
&= 2 g_0(R(\delta, Y)Y,Z)\vert_p + 2 Z(g_0(\nabla_Y \delta, Y))\vert_p - (\nabla_Y \Omega_{+})(Z,Y)\vert_p \\
&= 2 g_0(R(\delta, Y)Y,Z)\vert_p + Z(\Omega_{+}(Y,Y))\vert_p - (\nabla_Y \Omega_{+})(Z,Y)\vert_p \\
&= ( 2 g_0(R(\delta, Y)Y,Z) + (\nabla_Z \Omega_{+})(Y,Y) - (\nabla_Y \Omega_{+})(Z,Y) )\vert_p ,
\end{align*}
%--------------------------------------------------------------------
using the symmetries of the curvature tensor $R$. Since $p \in M$ is arbitrary and $\nabla_Y \Omega_{-}$ is antisymmetric, we have, for any $Z\in \X(M)$,
%--------------------------------------------------------------------
\begin{equation*}
(\nabla_Y \Omega_{-})(Y, Z) = - 2 g_0(R(\delta, Y) Y, Z) - (\nabla_Z \Omega_{+})(Y,Y) + (\nabla_Y \Omega_{+})(Z, Y) .
\end{equation*}
%--------------------------------------------------------------------
Plugging the above identity for $Z = X + \mu\delta$ into \eqref{prelimfinalresult} finishes the proof.
\end{proof}
%--------------------------------------------------------------------

\section*{Acknowledgments}
The authors warmly thank Professor Marcos Alexandrino for several helpful comments and stimulating discussions on the contents of this paper. MAJ  is partially supported by the grant PID2025-167696NB-I00, funded by MICIU/AEI/10.13039/501100011033, Spain; in part by FAPESP (Processo 2024/22850-0);  and in part by Fundación Séneca-Agencia de Ciencia y Tecnología de la Región de Murcia, under project ``Ayudas a proyectos para el desarrollo de investigación científica y técnica por grupos competitivos (Comunidad Autónoma de la Región de Murcia)'', included in ``Programa Regional de Fomento de la Investigación Científica y Técnica (Plan de Actuación 2022)'', REF. 21899/PI/22. PM is supported by FAPESP (Processo 2024/03446-4).

%  BIBLIOGRAPHY
\bibliography{bib}

@article {bao04z,
  title={{Zermelo navigation on Riemannian manifolds}},
  author={Bao, David and Robles, Colleen and Shen, Zhongmin},
  journal={Journal of Differential Geometry},
  volume={66},
  number={3},
  pages={377--435},
  year={2004}
}

@article {cap11i,
  title={{On the interplay between Lorentzian causality and Finsler metrics of Randers type}},
  author={Caponio, Erasmo and Javaloyes, Miguel \'Angel and S\'anchez, Miguel},
  journal={Revista Matem\'atica Iberoamericana},
  volume={27},
  number={3},
  pages={919--952},
  year={2011}
  }

@article {cap24w,
  title={{Wind Finslerian structures: From Zermelo’s navigation to the causality of spacetimes}},
  author={Caponio, Erasmo and Javaloyes, Miguel \'Angel and S\'anchez, Miguel},
  journal={Memoirs of the American Mathematical Society},
  volume={300},
  number={1501},
  pages={121 pp},
  year={2024}
}

@article {che12n,
  title={{A note on Randers metrics of scalar flag curvature}},
  author={Chen, Bin and Zhao, Lili},
  journal={Canadian Mathematical Bulletin. Bulletin Canadien de Math\'ematiques},
  volume={55},
  number={3},
  pages={474--486},
  year={2012}
}

@article {chen20c,
  title={{A class of Randers metrics of scalar flag curvature}},
  author={Cheng, Xinyue and Yin, Li and Li, Tingting},
  journal={International Journal of Mathematics},
  volume={31},
  number={13},
  pages={2050114},
  year={2020}
}

@article {gib09s,
  title={{Stationary metrics and optical Zermelo-Randers-Finsler geometry}},
  author={Gibbons, G. W. and Herdeiro, C. A. R. and Warnick, C. M. and Werner, M. C.},
  journal={Physical Review D. Particles, Fields, Gravitation, and Cosmology},
  volume={79},
  number={4},
  pages={044022},
  year={2009}
}

@article{jav20c,
  title={{Curvature computations in Finsler geometry using a distinguished class of anisotropic connections}},
  author={Javaloyes, Miguel {\'A}ngel},
  journal={Mediterranean Journal of Mathematics},
  volume={17},
  number={4},
  pages={123},
  year={2020}
}

@article{har82t,
  title={{A triangle comparison theorem for Lorentz manifolds}},
  author={Harris, Stacey G},
  journal={Indiana University Mathematics Journal},
  volume={31},
  number={3},
  pages={289--308},
  year={1982}
}

@article{jav20a,
  title={{Anisotropic conformal invariance of lightlike geodesics in pseudo-Finsler manifolds}},
  author={Javaloyes, Miguel {\'A}ngel and Soares, Bruno Learth},
  journal={Classical and Quantum Gravity},
  volume={38},
  number={2},
  pages={025002},
  year={2021}
}

@book{bao12i,
  title={{An introduction to Riemann-Finsler geometry}},
  author={Bao, David and Chern, Shiing-Shen and Shen, Zhongmin},
  series={Graduate Texts in Mathematics},
  volume={200},
  year={2000},
  publisher={Springer}
}

@article{per06f,
  title={{Fermat principle in Finsler spacetimes}},
  author={Perlick, Volker},
  journal={General Relativity and Gravitation},
  volume={38},
  number={2},
  pages={365--380},
  year={2006}
}

@article{jav20d,
  title={{On the definition and examples of cones and Finsler spacetimes}},
  author={Javaloyes, Miguel {\'A}ngel and S{\'a}nchez, Miguel},
  journal={Revista de la Real Academia de Ciencias Exactas, F{\'\i}sicas y Naturales. Serie A. Matem{\'a}ticas},
  volume={114},
  number={1},
  pages={30},
  year={2020}
}

@article{jav14c,
  title={{Chern connection of a pseudo-Finsler metric as a family of affine connections}},
  author={Javaloyes, Miguel {\'A}ngel},
  journal={Publicationes Mathematicae Debrecen},
  volume={84},
  number={1-2},
  pages={29--43},
  year={2014}
}

@article{jav19a,
  title={{Anisotropic tensor calculus}},
  author={Javaloyes, Miguel {\'A}ngel},
  journal={International Journal of Geometric Methods in Modern Physics},
  volume={16},
  number={supp02},
  pages={1941001},
  year={2019}
}

@article{jav14d,
  title={{On the definition and examples of Finsler metrics}},
  author={Javaloyes, Miguel {\'A}ngel and S{\'a}nchez, Miguel},
  journal={Annali della Scuola Normale Superiore di Pisa. Classe di Scienze},
  volume={13},
  number={3},
  pages={813--858},
  year={2014}
}

@article {lis09r,
  title={{On Randers metrics of quadratic Riemann curvature}},
  author={Li, Benling and Shen, Zhongmin},
  journal={International Journal of Mathematics},
  volume={20},
  number={3},
  pages={369--376},
  year={2009}
}

@book{asa12f,
  title={{Finsler Geometry, Relativity and Gauge Theories}},
  author={Asanov, Gennadii S},
  year={1985},
  series={Fundamental Theories of Physics},
  volume={12},
  publisher={D. Reidel Publishing Company / Springer}
}

@book{kob69f,
  title={{Foundations of differential geometry, Vol. II}},
  author={Kobayashi, Shoshichi and Nomizu, Katsumi},
  series={Interscience Tracts in Pure and Applied Mathematics},
  volume={15},
  publisher={Interscience / John Wiley \& Sons},
  year={1969}
}

@article{pfe12f,
  title={{Finsler geometric extension of Einstein gravity}},
  author={Pfeifer, Christian and Wohlfarth, Mattias N R},
  journal={Physical Review D},
  volume={85},
  number={6},
  pages={064009},
  year={2012}
}

@book{she01d,
  title={{Differential Geometry of Spray and Finsler Spaces}},
  author={Shen, Zhongmin},
  year={2001},
  publisher={Springer Dordrecht}
}

@book{mat80z,
  title={{Zwei Verallgemeinerungen eines Satzes von Gromoll und Meyer}},
  author={Matthias, Hans-Helmut},
  year={1980},
  publisher={Bonn Mathematical Publications},
  series={Bonner Mathematischen Schriften},
  number={126}
}

@book{one83s,
  title={{Semi-Riemannian Geometry with applications to Relativity}},
  author={O'Neill, Barrett},
  series={Pure and Applied Mathematics},
  volume={103},
  year={1983},
  publisher={Academic Press}
}

@article{alb10g,
  title={{A geometrical interpretation of the null sectional curvature}},
  author={Albujer, Alma L and Haesen, Stefan},
  journal={Journal of Geometry and Physics},
  volume={60},
  number={3},
  pages={471--476},
  year={2010}
}

@article{pal07f,
  title={{The fibre bundle of degenerate tangent planes of a Lorentzian manifold and the smoothness of the null sectional curvature}},
  author={Palomo, Francisco J},
  journal={Differential Geometry and its Applications},
  volume={25},
  number={6},
  pages={667--673},
  year={2007}
}

@incollection{flo03g,
  title={{Geodesics in stationary spacetimes. Application to Kerr spacetime}},
  author={Flores, Jos{\'e} Luis and S{\'a}nchez, Miguel},
  booktitle={Theoretical Physics 2002, Part 2},
  editor={George, Thomas F. and Arnoldus, Henk F.},
  publisher={Nova Science Publishers},
  year={2003},
  pages={141--158}
}

@article{she13c,
  title={{A class of Randers metrics of scalar flag curvature}},
  author={Shen, Zhongmin and Xia, Qiaoling},
  journal={International Journal of Mathematics},
  volume={24},
  number={07},
  pages={1350055},
  year={2013}
}

@article{cap11e,
  title={{On the energy functional on Finsler manifolds and applications to stationary spacetimes}},
  author={Caponio, Erasmo and Javaloyes, Miguel {\'A}ngel and Masiello, Antonio},
  journal={Mathematische Annalen},
  volume={351},
  number={2},
  pages={365--392},
  year={2011}
}

@book{bes87e,
  title={{Einstein Manifolds}},
  author={Besse, Arthur L.},
  series={Classics in Mathematics},
  year={1987},
  publisher={Springer}
}

@article{gray67p,
  title={{Pseudo-Riemannian almost product manifolds and submersions}},
  author={Gray, Alfred},
  journal={Journal of Mathematics and Mechanics},
  volume={16},
  number={7},
  pages={715--737},
  year={1967}
}

@article{one66f,
  title={{The fundamental equations of a submersion}},
  author={O'Neill, Barrett},
  journal={Michigan Mathematical Journal},
  volume={13},
  number={4},
  pages={459--469},
  year={1966}
}

@article{mak18d,
  title={{Differential geometric aspects of causal structures}},
  author={Makhmali, Omid},
  journal={SIGMA. Symmetry, Integrability and Geometry. Methods and Applications},
  volume={14},
  number={080},
  pages={50 pp},
  year={2018}
}

@article{min19c,
  title={{Causality theory for closed cone structures with applications}},
  author={Minguzzi, Ettore},
  journal={Reviews in Mathematical Physics},
  volume={31},
  number={5},
  pages={1930001},
  year={2019}
}

@article{pfe19f,
  title={{Finsler spacetime geometry in physics}},
  author={Pfeifer, Christian},
  journal={International Journal of Geometric Methods in Modern Physics},
  volume={16},
  pages={1941004},
  year={2019}
}

@article{hee27e,
  title={{Einstein-Kropina metrics and their application in Finsler gravity}},
  author={Heefer, Sjors and Villase\~nor, Fidel F. and Fuster, Andrea},
  journal={Journal of Geometry and Physics},
  volume={231},
  pages={105970},
  year={2027}
}

@article{edw18r,
  title={{Riemann-Finsler geometry and Lorentz-violating scalar fields}},
  author={Edwards, Benjamin R. and Kosteleck\'y, V. Alan},
  journal={Physics Letters. B. Particle Physics, Nuclear Physics and Cosmology},
  volume={783},
  pages={319--326},
  year={2018}
}

@article {hua15f,
  title={{On the flag curvature of a class of Finsler metrics produced by the navigation problem}},
  author={Huang, Libing and Mo, Xiaohuan},
  journal={Pacific Journal of Mathematics},
  volume={277},
  number={1},
  pages={149--168},
  year={2015}
}

@article{xia13f,
  title={{On the flag curvature of a class of Randers metric generated from the navigation problem}},
  author={Xia, Qiaoling},
  journal={Journal of Mathematical Analysis and Applications},
  volume={397},
  number={1},
  pages={415-427},
  year={2013}
}

@article{mo07o,
  title={{On curvature decreasing property of a class of navigation problems}},
  author={Mo, Xiaohuan and Huang, Libing},
  journal={Publicationes Mathematicae Debrecen},
  volume={71},
  number={1-2},
  pages={141--163},
  year={2007}
}

%--------------------------------------------------------------------
\end{document}